\documentclass[11pt,reqno]{amsart}

\usepackage{amsmath}\allowdisplaybreaks  
\usepackage{amssymb}        
\usepackage{amsfonts}       
\usepackage{amsthm}         
\usepackage{mathrsfs}       
\usepackage{cases}          
\usepackage{siunitx}        

\usepackage{array}
\usepackage{booktabs}       
\usepackage{longtable}      
\usepackage{multirow}
\usepackage{makecell}       
\usepackage{tabularx}       
\usepackage{nicematrix}     
\usepackage{threeparttable}

\usepackage{tikz}\usetikzlibrary{tikzmark, arrows.meta}   
\usetikzlibrary{positioning}
\usepackage{float} 
\usepackage{xcolor}           
\usepackage{graphicx}         

\usepackage[linesnumbered, ruled, vlined]{algorithm2e}  

\usepackage[colorlinks=true, linkcolor=blue, citecolor=blue, urlcolor=blue]{hyperref} 
\usepackage[nameinlink]{cleveref}	

\usepackage[left=2cm, right=2cm, top=2cm, bottom=2cm]{geometry}

\usepackage{ulem}           
\usepackage{lineno}         

\makeatletter
\newcommand\figcaption{\def\@captype{figure}\caption} 
\newcommand\tabcaption{\def\@captype{table}\caption}
\makeatother

\newtheorem{theorem}{Theorem}[section]          
\newtheorem{lemma}[theorem]{Lemma}              
\newtheorem{proposition}[theorem]{Proposition}  
\newtheorem{definition}[theorem]{Definition}    
\newtheorem{example}[theorem]{Example}          
\newtheorem{remark}[theorem]{Remark}            

\begin{document}
	\title[Neural network-enhanced AFEM]{Neural network-enhanced $hr$-adaptive finite element algorithm for parabolic equations}
	\author[J.~Hao, Y.~Huang, N.~Yi, P.~Yin]{Jiaxiong Hao$^{\dag}$, Yunqing Huang$^\S$, Nianyu Yi$^\dag$, Peimeng Yin$^{\ddag}$}
	\address{$\dag$ Hunan Key Laboratory for Computation and Simulation in Science and Engineering, School of Mathematics and Computational Science, Xiangtan University, Xiangtan 411105, Hunan, P.R.China} \email{moodbear@qq.com (J. Hao);\  yinianyu@xtu.edu.cn (N. Yi)}
	\address{$\S$ National Center for Applied Mathematics in Hunan, Key Laboratory of Intelligent Computing \& Information Processing of Ministry of Education, Xiangtan University, Xiangtan 411105, Hunan, P.R.China} \email{huangyq@xtu.edu.cn}
	\address{$^\ddag$ Department of Mathematical Sciences, University of Texas at El Paso,  El Paso, Texas 79968, USA}\email{pyin@utep.edu}
	
	\begin{abstract}
		In this paper, we propose a novel $hr$-adaptive finite element method, enhanced by neural networks, for parabolic equations. The main challenge of the conventional $h$-adaptive finite element method is interpolating the finite element solution from the previous step in the updated mesh. Recomputing this mesh-dependent interpolation at each iteration is both computationally expensive and demands considerable implementation expertise. The new approach addresses this challenge by introducing a neural network to construct a mesh-free surrogate of the finite element solution of the previous step. Since the neural network is mesh-free, it only requires training once per time step, with its parameters initialized using the minimizer of the previous time step. This approach effectively overcomes the interpolation challenges associated with non-nested meshes in computation, making node insertion and movement more convenient and efficient. The new algorithm computes the local mesh density using a mesh size field (SIZE) and then generates a new adapted mesh based on the density (GENERATE), allowing each refinement to roughly double the number of mesh nodes of the previous iteration and then redistribute them to form a new mesh that effectively captures the singularities. It significantly reduces the time required for repeated refinement of the conventional methods and achieves the desired accuracy in no more than seven space-adaptive iterations per time step. Numerical experiments confirm the efficiency of the proposed algorithm in capturing dynamic changes of singularities. The code is made publicly available on GitHub.\footnotemark
	\end{abstract}
	\footnotetext{GitHub: \url{https://github.com/FEMmaster/NN-enhanced-hr-AFEA-for-parabolic-equations}.}
	
	
	
	\maketitle
	\section{Introduction}
	The adaptive finite element method (AFEM) is one of the effective numerical approaches for solving partial differential equations, particularly showing strong adaptability for problems involving singularities or multiscale features. 
	The fundamental idea behind adaptive finite element methods is to achieve an optimal asymptotic error decay rate by refining the finite element space. 
	Generally, there are three primary approaches to refining a finite element space: the $h$-adaptive method, which improves the accuracy of finite element approximation by refining the mesh; the $p$-adaptive method, which uses a fixed mesh but increases the polynomial degree of the shape functions to improve its accuracy; and $r$-adaptivity (or moving mesh methods), which relocates mesh nodes to concentrate them in regions of interest without changing the total number of nodes. 
	Naturally, these strategies can be synergistically combined to achieve greater efficiency and robustness than any single approach alone. 
	Among them, the two most prominent combinations are the $hp$- and $hr$-adaptive methods, coupling $h$-adaptivity with $p$- and $r$-adaptivity, respectively.
	The $hp$-adaptive method is renowned for achieving an exponential convergence rate in the energy norm \cite{Babuka1987, Melenk2002}. 
	However, its practical implementation, especially in high dimensions, is notoriously challenging due to the complexity of managing both mesh refinement and variable polynomial degrees. 
	Given these challenges, this work focuses on the $hr$-adaptive framework.
	
	Irrespective of the specific adaptive strategy, the key mechanism of an AFEM is to drive the mesh adaptation process (refinement or coarsen) based on an a posteriori error indicator. 
	The a posteriori error estimator, which is a computable quantity depending on the numerical solution, provides information about the distribution of the error of the finite element approximation \cite{Ainsworth1997, Verfrth2013}. In particular, for the $hr$-AFEM, the success of such an algorithm hinges on a reliable a posteriori error estimator to guide both mesh refinement and node distribution to achieve better accuracy with minimum degrees of freedom. 
	A posteriori error analysis, along with convergence and optimality theories for adaptive finite element algorithms, are well established for elliptic equations \cite{Bonito2024, Mekchay2005, Morin2002}. 
	Corresponding theories are available for parabolic equations \cite{Eriksson1991, Lakkis2011, Picasso1998}. 
	However, translating this theoretical framework into a robust computational procedure is nontrivial, as applying AFEM to time-dependent problems presents unique challenges not encountered in stationary elliptic settings. 
	The sequential nature of time-stepping and the evolution of solution features exacerbate the computational cost of mesh adaptation and the error propagation from frequent data transfer. 
	Consequently, conventional AFEMs suffer from the following inherent fundamental bottlenecks when applied to such problems:
	\begin{itemize}
		\item Inefficient refinement: The refinement process in each time layer can require up to $30-40$ iterations \cite{Picasso1998, Xiao2026}. 
		\item Interpolation challenge: The interpolation of the solution from the previous time step onto each newly refined mesh is algorithmically demanding and a significant computational bottleneck.   
		\item Redundant initial mesh choice at each time step: At each new time step, the adaptive cycle begins with the final mesh from the previous step, resulting in a linear system with an excessively large number of degrees of freedom. Here, an additional and intricate mesh coarsening step (COARSEN) must be enforced. 
	\end{itemize}
	
	In pursuit of an efficient and straightforward adaptive finite element method for parabolic problems, this work proposes a novel neural network-enhanced $hr$-adaptive finite element algorithm. 
	Its key innovation lies in employing a data-driven neural network as our primary tool, which serves as a continuous function representation. 
	The choice of this tool is natural, given the proven capacity of neural networks to act as universal function approximators and to learn efficient representations from data. 
	Nowadays, neural networks have achieved remarkable success in artificial intelligence, enabling complex tasks in domains such as speech recognition, natural language processing, computer vision, and autonomous vehicles through data-trained models \cite{He2016, Khurana2022}. 
	In scientific computing, neural networks leverage their dual strengths—as universal approximators \cite{Barron1993, Hornik1991, Hornik1989, Siegel2020, Yarotsky2017} and for mitigating high-dimensional complexity \cite{E2021}—underpin their fundamental capability to handle complex, data- and physics-driven computational tasks. 
	This capability is vividly evidenced by the growing number of studies that directly apply neural networks to solving partial differential equations \cite{Cai2022, Chen2022, E2018, He2020, Lagaris1998, Raissi2019, Zeng2022, Xu2020}. 
	
	
	Despite this rapid expansion, several key limitations have become apparent when deep learning is used as a stand-alone PDE solver, including high training costs, limited interpretability, and accuracy that often remains below that of well-established numerical methods.
		To overcome these limitations while leveraging the strengths of both artificial intelligence and classical numerical methods,
		a growing trend is to employ artificial intelligence as a tool for augmentation rather than replacement, shifting the role of neural networks from stand-alone PDE solvers toward specialized components embedded within mature numerical frameworks.
		In this vein, neural networks have been incorporated into finite element methods in various ways. Examples include error-informed adaptive FEM enhanced by neural networks \cite{Oishi2020} and $h$-adaptive finite element–interpolated neural networks designed to resolve localized features \cite{Badia2025}. Further related developments and applications can be found in \cite{Aballay2026, Caboussat2024, Han2023, Xiong2025} and the references therein. These works illustrate how neural networks can be strategically integrated to complement, rather than replace, classical numerical methods.

		
		Following the strategy of neural networks as specialized modules, our proposed method effectively addresses the three aforementioned limitations. 
		This method is more efficient in terms of adaptive performance compared to the conventional refinement strategies discussed earlier. 
		Specifically, the advantages of the proposed neural network-enhanced $hr$-AFEM in addressing the three limitations of the conventional AFEMs are as follows:
		\begin{itemize} 
			\item Efficient refinement: By controlling the refinement of mesh nodes at each step, the number of iterations required for convergence is reduced to a maximum of $7$ iterations.
			\item Interpolation free: By training a neural network to replace the prior finite element solution, the method bypasses the complex and costly steps of interpolation. With only a single fine-tuning required per time step, it then enables direct function evaluation at any spatial point.
			\item Cheap initial mesh choices per time step: By employing a combined error indicator that incorporates information from two adjacent time steps, the method ensures reliable re-initializing of the mesh to a prescribed, coarse uniform state at each time step.
		\end{itemize}
		Together, these mechanisms—drastic iteration reduction, elimination of interpolation costs, and avoidance of redundant initial computations—constitute a comprehensive efficiency gain.
		The rest of the paper is outlined as follows: 
		Section \ref{hafem} introduces the model problem and the fully discrete scheme for the parabolic equation.
		Then, we recap the conventional adaptive algorithm, offering a comprehensive analysis of its methodology while critically examining its key limitations and potential drawbacks.
		In response, Section \ref{hrafealg} discusses the development of an efficient neural network-enhanced $hr$-adaptive finite element algorithm, specifically highlighting the main advantages over conventional methods.
		In Section \ref{numExp}, several numerical experiments are presented to demonstrate the robustness and effectiveness of the proposed method.

		\section{Standard adaptive finite element algorithm}\label{hafem}
		\subsection{Model problem}
		Let $\Omega$ be a polygonal or polyhedral domain in $\mathbb{R}^d$ (with $d=2,3$), and $T>0$. Consider the following parabolic equation:
		\begin{equation} \label{model}
			\left\{\begin{aligned}
				u_t - \nabla\cdot(a\nabla u) &= f, & \text{in } &  \Omega \times (0, T], \\
				u(\mathbf{x}, t) &= g, & \text{on } &  \partial\Omega \times (0, T], \\
				u(\mathbf{x}, 0) &= u_0, & \text{in } &  \Omega,
			\end{aligned}\right.
		\end{equation}
		where $f\in L^2(0, T; L^2(\Omega))$, $u_0\in L^2(\Omega)$, and $a\in L^{\infty}(\Omega)$ are uniformly bounded and positive in the sense that there exist positive constants $\lambda_{\min}$ and $\lambda_{\max}$ satisfying \[\lambda_{\min}\leq a(\mathbf{x}) \leq\lambda_{\max},\quad \forall \mathbf{x}\in\bar{\Omega}.\]
		This type of equation plays a critical role in applications such as heat conduction models \cite{Rogolino2018}, image denoising \cite{Aboulaich2008}, material science \cite{Wang2014}, and financial mathematics \cite{BARLES1995}.
		The solution to the parabolic equation is generally smooth in most regions. 
		However, discontinuities in initial conditions, nonuniform material properties, or interactions between multiple physical processes can give rise to singularities that emerge and gradually evolve over time.
		This requirement arises from the fact that the smoothness of solutions to parabolic equations may vary considerably across different regions, necessitating finer spatial discretization and more sophisticated numerical techniques to resolve the solution behavior accurately.
		
		Throughout this work, we use the standard notation $W^{m,q}(D)$ to refer to Sobolev spaces on an open set $D\subset \mathbb{R}^d$, with the norm $\Vert\cdot\Vert_{W^{m,q}(D)}$ and the semi-norm $\vert\cdot\vert_{W^{m,q}(D)}$. The space $W^{m,q}_0(D)$ is defined as $\{w \in W^{m,q}(D): w|_{\partial D} = 0\}$. Furthermore, $W^{m,2}(D)$ and $W^{m,2}_0(\Omega)$ are denoted by $H^m(D)$ and $H^m_0(D)$, respectively. 
		We also denote by $L^{p}(0,T;H^k(D))$ the space of all $L^p$ integrable functions from $(0,T)$ into $H^k(D)$ with norm $\|v\|_{L^{p}(0,T;H^k(D))}=\left(\int_0^T \|v\|_{k,D}^pdt\right)^{\frac{1}{p}}$ for $p\in [1,\infty)$. 
		
		Given $f\in L^2(0, T; L^2(\Omega))$ and $u_0\in L^2(\Omega)$, the variational formulation of \eqref{model} reads: Find $u \in L^2(0, T; H^1(\Omega))\bigcap H^{1}(0, T; H^{-1}(\Omega))$ such that
		\begin{equation}	\label{weak form}
			\left\{\begin{aligned}
				(u_t,v)+(a\nabla u,\nabla v) & =(f,v),  & \forall v \in H_0^1(\Omega),\quad  t\in(0,T], \\
				u(\mathbf{x},0) & =u_0.
			\end{aligned}\right.
		\end{equation}
		
		\subsection{Fully discrete FEM}
		Let the time step size be defined as $\tau = T/S$, which partitions the interval $[0, T]$ uniformly into $S$ subintervals and produces the discrete time levels $t_n = n\tau$, for $n = 0, \dots, S$.
		Let $\mathcal{T}_h^n$ be a shape-regular triangulation of $\Omega$ at the time level $t_n$, whose associated linear finite element space is
		\[V_h^n=\{v\in H_0^1(\Omega):v|_K\in P_1(K),\forall K\in \mathcal{T}_h^n\},\]
		where $P_1(K)$ denotes the space of linear polynomials on $K\subseteq R^d(d=2,3)$. 
		Let $u_h^0=Pu_0$, where $P: L^2(\Omega)\rightarrow V_h^0$ is the $L^2$ projection operator onto the finite element space $V_h^0$ defined on the initial mesh $\mathcal{T}_h^0$. 
		Applying the backward Euler method for time discretization and the finite element method for spatial discretization to the weak formulation \eqref{weak form} yields the following fully discrete scheme: Find a sequence of functions $u_h^n \in V_h^n,\,n=1,2,\cdots, S$, such that
		\begin{equation}	\label{full-discrete}
			\left\{\begin{aligned}
				\left(\frac{u_h^n-\Pi_h^n u_h^{n-1}}{\tau},v\right)+(\nabla{u_h^n},\nabla v ) &= \left(f_h^{n},v\right), \quad & \forall v \in V_h^n, \\
				u_h(\mathbf{x},0) & =u_h^0,
			\end{aligned}\right.
		\end{equation} 
		where the operator $\Pi_h^n: V_h^{n-1} \to V_h^n$ denotes the nodal interpolation operator from the previous time level's finite element space onto the current one. 
		Note that at each time level $t_n$, the adaptive finite element algorithm generates a sequence of meshes ${\mathcal{T}^{n,0}_h, \mathcal{T}^{n,1}_h, \dots}$ along with their corresponding continuous piecewise linear finite element spaces ${V^{n,0}_h, V^{n,1}_h, \dots}$. 
		Hereafter, we do not distinguish between the intermediate meshes generated within the $n$-th time step and denote them collectively by $\mathcal{T}_h^n$, and similarly by $V_h^n$ for the associated finite element spaces.
		
		\subsection{Recap of the standard $h$-AFEM}
		The standard $h$-adaptive finite element method for parabolic equations proceeds sequentially in time, with multiple mesh updates possible within a single time step. 
		Given the numerical solution $u_h^{n-1}$ from the previous time level $t_{n-1}$, the iterative algorithm for computing the solution $u_h^n$ at time $t_n$ can be summarized as the following workflow:
		\begin{center}
			\begin{tikzpicture}[node distance=0.5cm]        
				\node (node1) {$u_h^{n-1}$};
				\node[right=of node1] (node2) {\textit{INTERPOLATE}};
				\node[right=of node2] (node3) {\textit{SOLVE}};
				\node[right=of node3] (node4) {\textit{ESTIMATE}};
				\node[right=of node4] (node5) {\textit{REFINE}};
				\node[right=of node5] (node6) {$u_h^n$};
				\node[right=of node6] (node7) {\textit{COARSEN}};
				
				\draw[->] (node1) -- (node2);
				\draw[->] (node2) -- (node3);
				\draw[->] (node3) -- (node4);
				\draw[->] (node4) -- (node5);
				\draw[->] (node6) -- (node7);
				
				\draw[->] (node5.south) -- ++(0,-0.3) -| (node2.south); 
				\draw[->] (node4.north) -- ++(0,0.4) -| (node6.north);      
			\end{tikzpicture}
		\end{center}
		Here, the \textit{INTERPOLATE} step computes $\Pi_h^n u_h^{n-1} \in V_h^n$ from $u_h^{n-1}$, 
		and is required whenever a mesh change occurs.
		The essential \textit{COARSEN} step is applied once, at the end of the adaptive iteration loop within a time step. 
		Its role is to coalesce elements with low estimated error. 
		Typically, these elements lie in regions that were highly refined to capture features of the previous solution $u_h^{n-1}$, but become redundant for subsequent computations once the new solution $u_h^n$ is obtained. 
		If not, these elements would accumulate, contributing to an unchecked growth in the total number of degrees of freedom and ultimately leading to computational blow-up. 
		The process of refinement and coarsening will be further illustrated in \Cref{exp2d}. 
		While mesh coarsening improves computational efficiency, it risks coarsening elements that were just refined in the previous iteration, which can degrade the accuracy of the subsequent solution. 
		
		We also recap the gradient recovery-based a posteriori error estimators, which use a certain norm of the difference between direct and post-processed gradient approximations as an indicator, have been widely used in the engineering and scientific computation community since the work of Zienkiewicz and Zhu \cite{Zienkiewicz1992}. 
		The success of recovery-based error estimators stems from their computational efficiency, ease of implementation, and general asymptotic exactness \cite{Chen2024, Huang2010, Lakkis2011, Liu2025Enhanced, Zhang2005}.
		Let $G: V_h^n\rightarrow V_h^n\times V_h^n$ denote the gradient recovery operator, and let $\mathcal{N}_h^n$ denote the set of vertices of mesh $\mathcal{T}_h^n$.
		For each node $z\in\mathcal{N}_h^n$, let $\omega_z$ be the element patch associated with $z$, 
		we define the recovered gradient $G(\nabla u_h^n(z))$ as 
		\[
		G(\nabla u_h^n(z)):=\sum_{K_i\in \omega_z}\frac{1/|K_i|}{\sum_{K_j\in \omega_z}1/|K_j|}\nabla u_h^n(c_i),
		\]
		where $c_i$ denotes the center of the element $K_i$. 
		The recovered gradient $G(\nabla u_h^n)$ over the whole domain is then obtained via interpolation:
		\[G(\nabla u_h^n)=\sum_{z\in\mathcal{N}_h^n}G(\nabla u_h^n(z))\phi_z^n,\]
		where $\phi_z^n$ is the Lagrange basis of finite element space $V_h^n$ associated with $\mathcal{T}_h^n$.
		The corresponding local and global gradient recovery-based error estimators are defined separately as
		\begin{equation}\label{estimators}
			\hat{\eta}_{h, K}^n:=\|G(\nabla u_h^n)-\nabla u_h^n\|_{0, K},\quad \hat{\eta}_{h}^n:=\|G(\nabla u_h^n)-\nabla u_h^n\|.
		\end{equation}
		
		By integrating the ideas outlined above, we introduce \Cref{alg:Parabolic1}.
		
		\normalem
		\begin{algorithm}[ht]
			\caption{Standard $h$-AFEM for the parabolic equation \eqref{model}}
			\label{alg:Parabolic1}
			\KwIn{Domain $\Omega$, source term $f$, tolerance $eTol$, initial data $u_0$, time step size $\tau$}
			\KwOut{The meshes $\{\mathcal{T}_h^n\}_{n=1}^S$ and finite element approximations $\{u_h^n\}_{n=1}^S$} 
			Generate an initial quasi-uniform mesh $\mathcal{T}_h^0$\;
			\For{$n:=0$ \KwTo $S$}{
				\While{true}{
					\If {$n==0$}{
						(PROJECT) Compute the projection $u_h^0=Pu_0$\;
					}
					\Else { 
						(SOLVE) Solve the problem \eqref{full-discrete} for $u_h^n$ on $\mathcal{T}_h^n$\;
					}
					
					(ESTIMATE) Compute the local and global error estimators $\{\hat{\eta}_{h, K}^n, \forall K \in \mathcal{T}_h^n\}$ and $\hat{\eta}_h^n$\;
					\If{$\eta_h^n \leq eTol$}{
						break\;
					}
					\Else{
						(REFINE) Refine the elements marked by the error estimator using a bisection algorithm to locally enhance the mesh resolution\;
						\If{$n>0$}{
							(INTERPOLATE) Compute the interpolation $\Pi_h^n u_h^{n-1}$\;
						}
					}
				}          
				$\mathcal{T}_h^{n+1}:=\mathcal{T}_h^{n}$, $t_{n+1}:=t_{n}+\tau$\;
				\If{$n>0$}{
					(COARSEN) Merge elements marked for coarsening to reduce the number of degrees of freedom in regions where the estimated error is sufficiently small.
				}
			}
		\end{algorithm}
		
		\begin{remark}
			\Cref{alg:Parabolic1}, as well as the following \Cref{alg:ParabolicHR}, is not restricted to any specific error estimators.
			Here, we choose the gradient recovery-based error estimator using the weighted average method primarily for its simplicity of implementation.
		\end{remark}
		
		\subsection{Limitations of $h$-AFEM}
\begin{figure}[ht]
			\centering
			\includegraphics[width=0.35\linewidth]{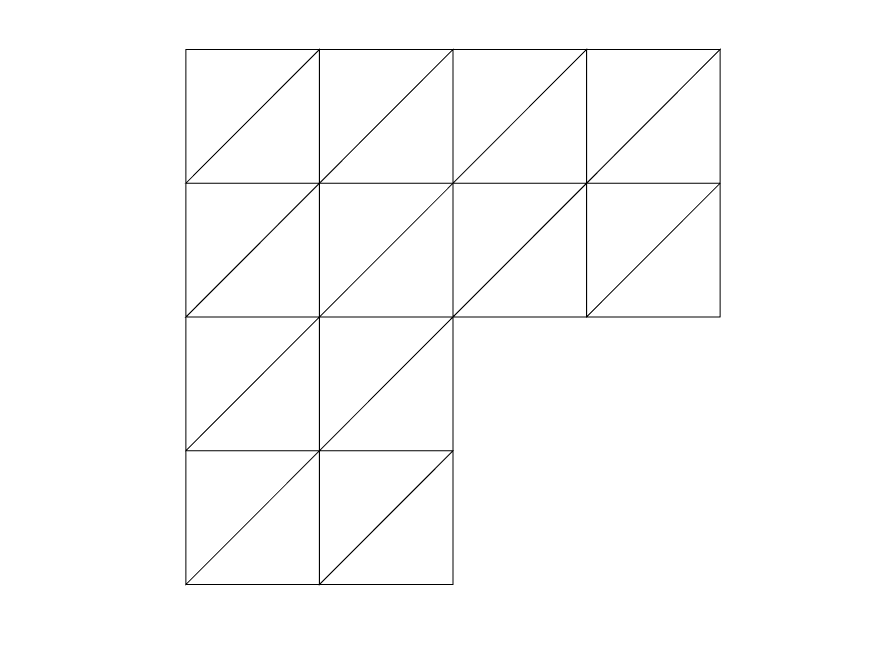}
			\includegraphics[width=0.35\linewidth]{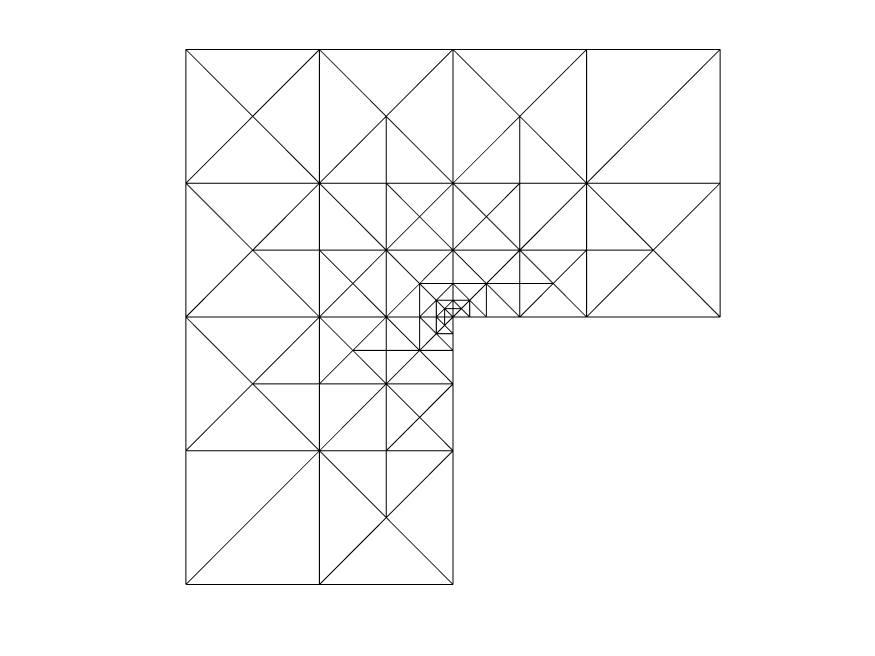}
			\caption{Left: initial uniform mesh; Right: $10$ times local refined mesh by bisection algorithm.}
			\label{Lshapemesh}
		\end{figure}

		It is important to highlight three primary limitations when applying \Cref{alg:Parabolic1} in practical computations. 
		First, the bisection algorithm is a layer-by-layer refinement method, which can lead to excessive adaptive iterations. 
		Coupled with the fact that each iteration executes the full INTERPOLATE-SOLVE–ESTIMATE–REFINE cycle, the overall computational expense becomes prohibitive. 
		A canonical example is the Laplace equation $-\Delta u = 0$ on the L‑shaped domain $[-1, 1]^2 \setminus (0,1) \times (-1,0)$ with Dirichlet boundary conditions, whose solution is  \[u(x,y)=r^{2/3}\sin(2\theta/3),\quad r=\sqrt{x^{2}+y^{2}},\quad \theta=\tan^{-1}(y/x).\]
		Applying the standard adaptive finite element method with newest‑vertex bisection (see \cite[Example 2.2]{Xiao2026}) to this problem yields the meshes shown in \Cref{Lshapemesh}. 
		The minimal element size decreases from $h_{K,\text{min}}=2^{-1}\sqrt{2}$ in the initial mesh to $h_{K,\text{min}}=2^{-5}\sqrt{2}$ after 10 adaptive steps, yet achieving a quasi‑optimal convergence rate typically requires $30‑40$ such refinements. 
		This exemplifies the high iterative cost inherent to the classical approach. 
		More numerical tests and analysis can be found in \cite{Xiao2026}.  
		
		Second, the interpolation of the solution from the previous mesh onto the newly refined mesh poses a significant implementation challenge. 
		This operation requires solving a point‑location problem by determining, for every newly created node, which element of the previous mesh contains it. 
		A straightforward implementation based on exhaustive element search is prohibitively expensive. 
		Although spatial search structures can accelerate this search, their reliable and robust integration remains algorithmically demanding. 
		Crucially, this geometrically non‑trivial step must be executed after every mesh refinement, making it a recurrent computational bottleneck. 
		
		Third, the standard practice of beginning each new time step’s adaptive process directly from the final mesh of the previous step incurs substantial and persistent overhead. 
		This forces all key computational components—matrix assembly, linear solver, error estimation, and data interpolation—to handle this inflated problem size in every iteration. 
		Moreover, to avoid carrying forward this irreversible computational baggage, the standard method is forced to insert an additional and intricate mesh coarsening step. 
		This step specifically targets and removes the fine mesh elements that were created to capture features of the old solution $u_h^{n-1}$. 
		We show the performance of \Cref{alg:Parabolic1} in \Cref{exp2d}. 
		
\begin{figure}[ht]
			\centering
			\includegraphics[width=0.22\linewidth]{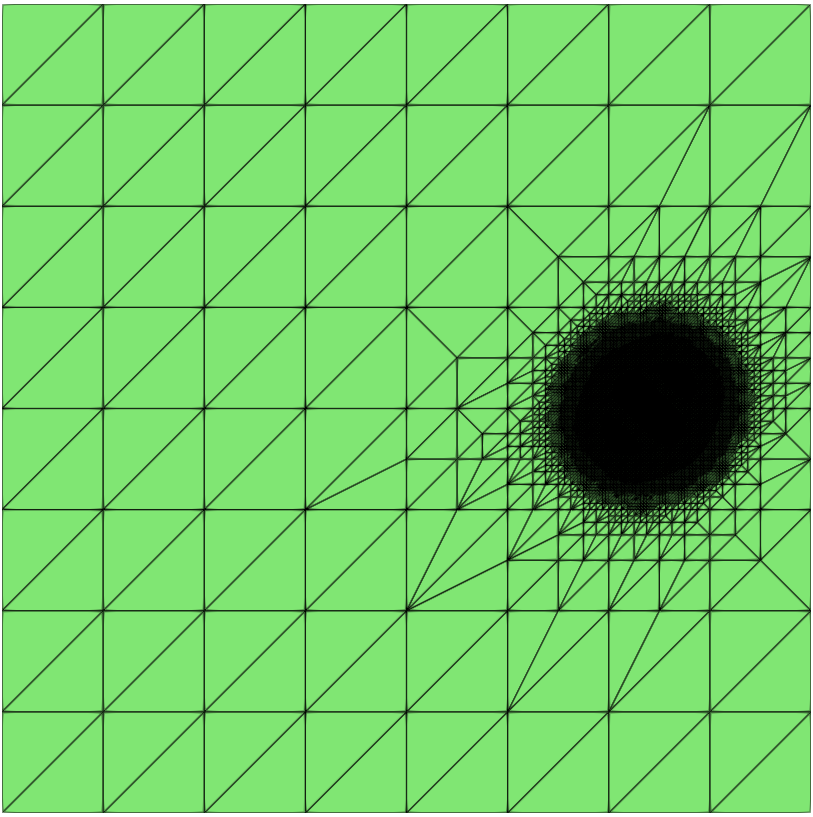}
			\includegraphics[width=0.22\linewidth]{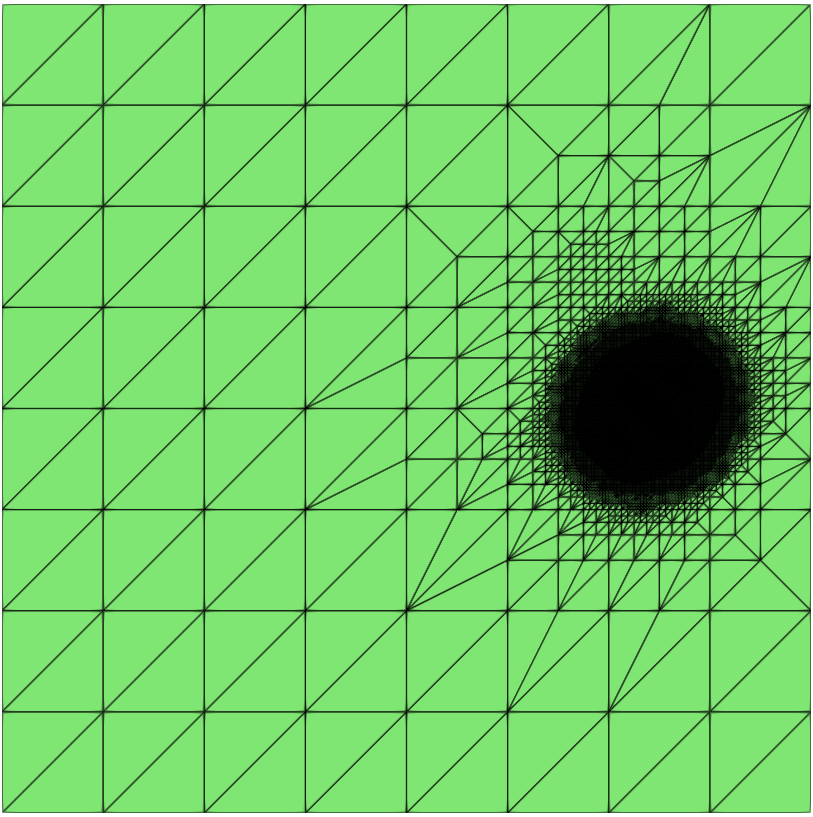}
			\includegraphics[width=0.22\linewidth]{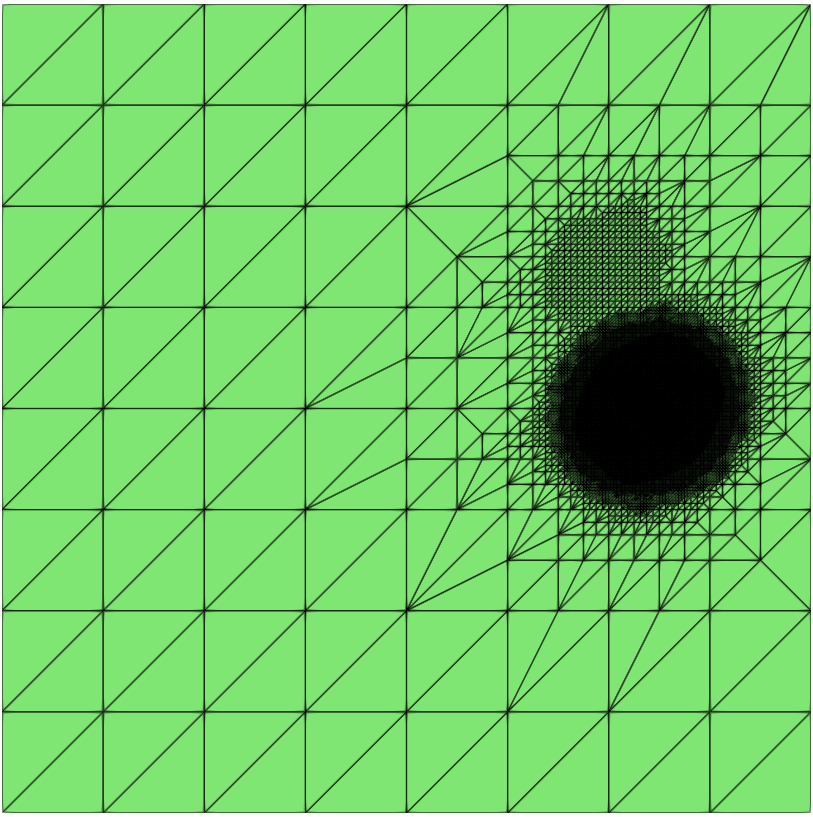}
			\includegraphics[width=0.22\linewidth]{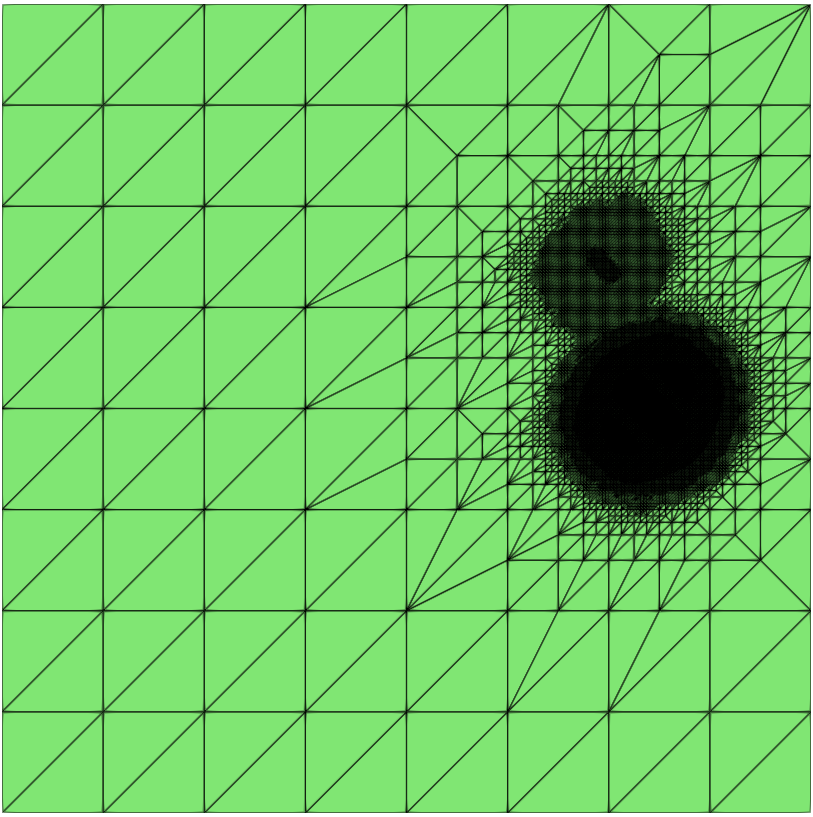}\vspace{0.1cm}
			\includegraphics[width=0.22\linewidth]{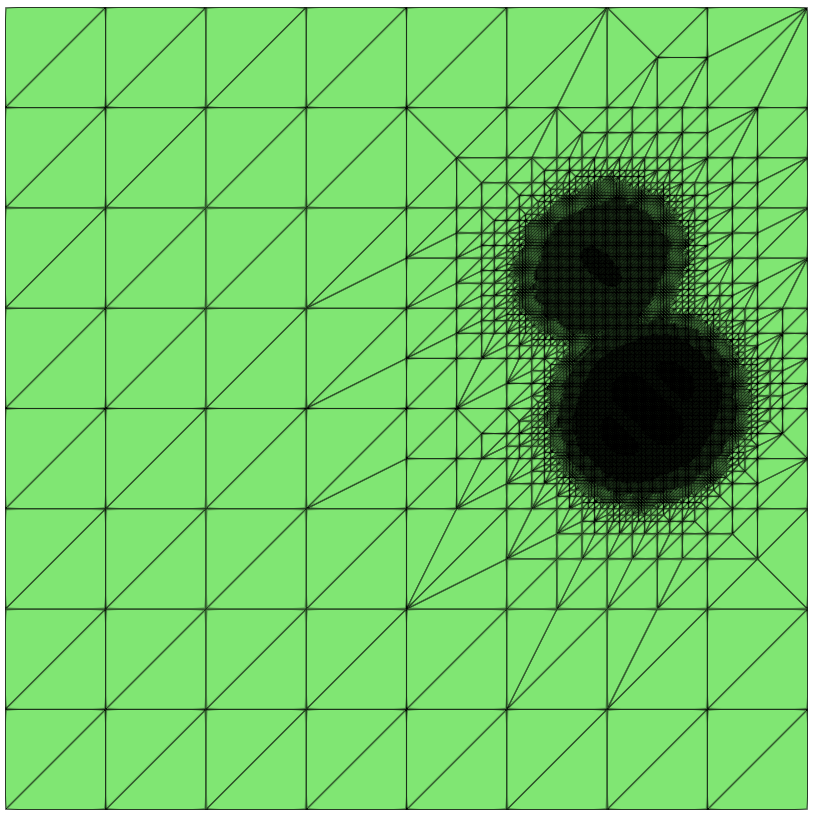}
			\includegraphics[width=0.22\linewidth]{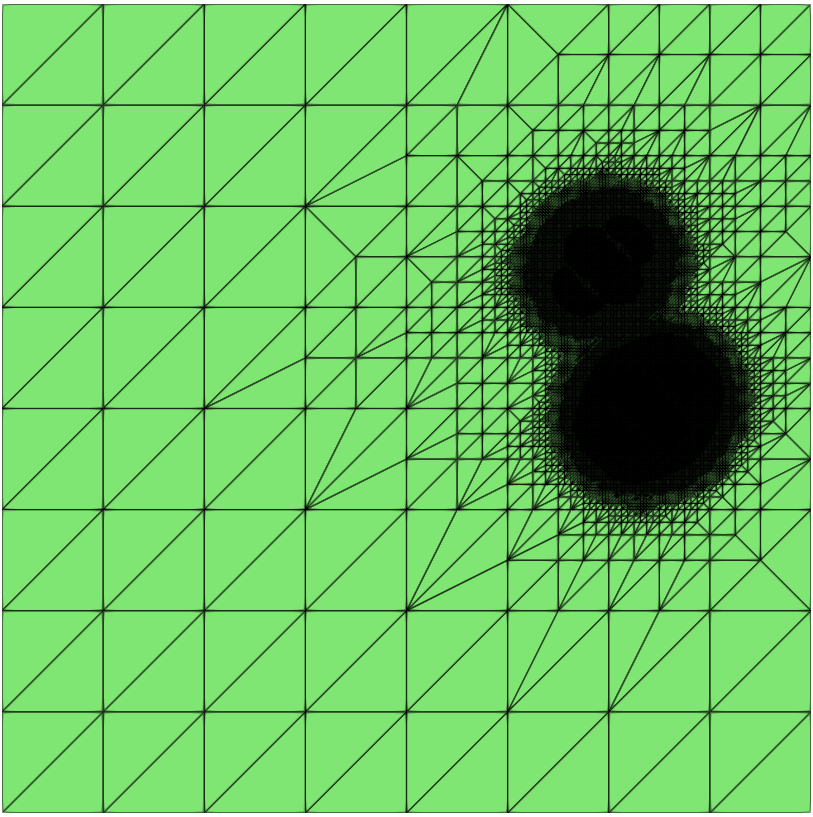}
			\includegraphics[width=0.22\linewidth]{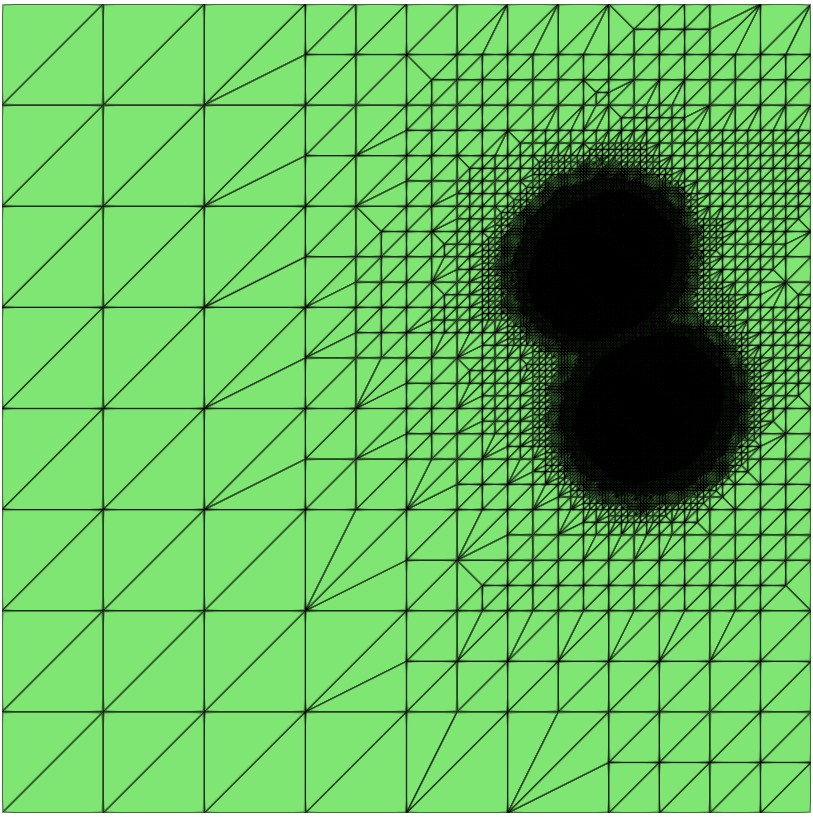}
			\includegraphics[width=0.22\linewidth]{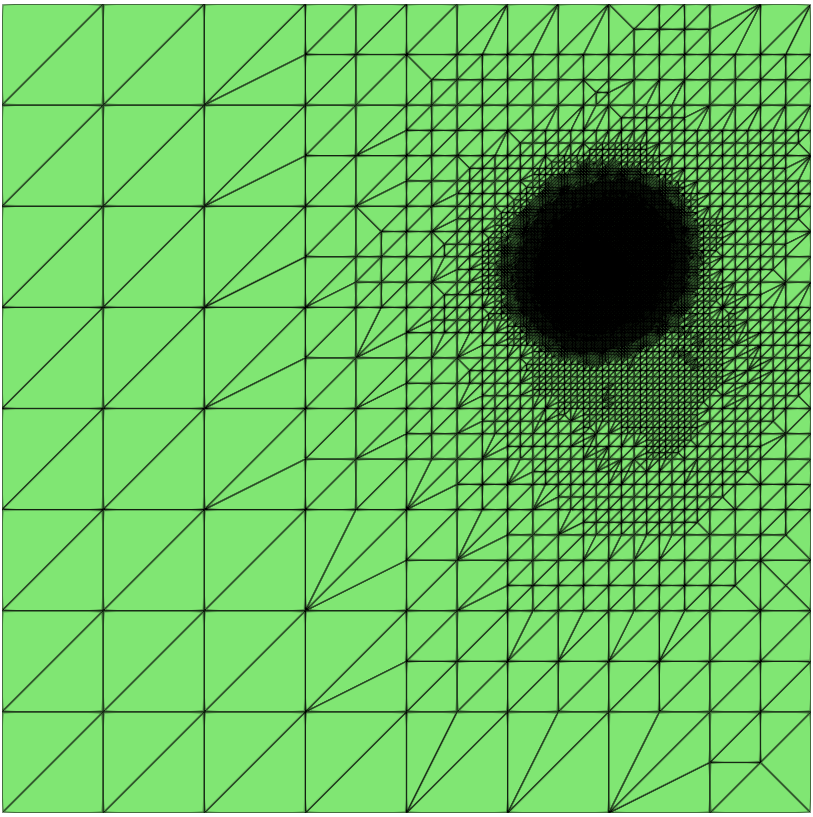}
			\caption{
				\Cref{exp2d} based on \Cref{alg:Parabolic1}: Initial mesh $\mathcal{T}^{1,0}_h$ at $t_1=0.1$ (top left); adaptively refined meshes $\mathcal{T}^{1,5}_h$, $\mathcal{T}^{1,10}_h$, $\mathcal{T}^{1,15}_h$, $\mathcal{T}^{1,20}_h$, $\mathcal{T}^{1,30}_h$, $\mathcal{T}^{1,34}_h$; and final coarsened mesh (bottom right) at $t=0.1$.}
			\label{paramesh}
		\end{figure}

		\begin{example}\label{exp2d}
			In this example, we consider the following parabolic equation 
			\begin{equation} \label{ex2d}
				\left\{\begin{aligned}
					u_t - \Delta u &= f, & \text{in } &  \Omega \times (0, 1], \\
					u &= g, & \text{on } &  \partial\Omega \times (0, 1], \\
					u(\mathbf{x}, 0) &= u_0, & \text{in } &  \Omega,
				\end{aligned}\right.
			\end{equation}
			where $\Omega=[-1,1]^2$. 
			The forcing term $f$ is chosen such that the exact solution is 
			\[u(x, y, t)=\exp(-500(x - 0.3\cos(2\pi t))^2)\exp(-500(y - 0.3\sin(2\pi t))^2).\]
			Boundary and initial conditions corresponding to this setting are directly prescribed as $g = u|_{\partial\Omega\times(0,1]}$ and $u_0 = u|_{t=0}$.
			
			We solve this problem using \Cref{alg:Parabolic1} with a tolerance $eTol = 0.01$ and a time step size $\tau = 0.1$. 
			\Cref{paramesh} presents the adaptively refined meshes at $t_1=0.1$. 
			The top left picture displays the initial mesh $\mathcal{T}^{1,0}_h$, obtained from a quasi-uniform mesh after $42$ refinements at $t_0=0.0$, with $54619$ vertices. 
			Consequently, all further adaptive refinements at this time step must start from and expand upon this already large set of degrees of freedom. 
			Subsequent frames, ordered left‑to‑right and top‑to‑bottom, show the progressively refined meshes which are $\mathcal{T}_h^{1,5}$, $\mathcal{T}_h^{1,10}$, $\mathcal{T}_h^{1,15}$, $\mathcal{T}_h^{1,20}$, $\mathcal{T}_h^{1,30}$, and $\mathcal{T}_h^{1,34}$. 
			The corresponding numbers of vertices (NOV) for these meshes are $54654$, $55452$, $57648$, $64340$, $82795$, and $108886$. 
			After solving for $u_h^1$, a final coarsening step produces the last mesh with $41599$ vertices. 
			From this example, we observe that the conventional AFEM replying on nested meshes requires large number of refinement iterations, incurs substantial computational overhead, and necessitates additional coarsening procedures. These factors collectively restrict the efficiency and flexibility of classical $h$-adaptivity. 
		\end{example}
		
		\section{Neural network-enhanced $hr$-AFEM}\label{hrafealg}
		Having identified the causes of these limitations, this section introduces a solution framework composed of two essential tools and a central strategy. 
		Specifically, we employ: (i) a mesh generation tool to perform local refinement; (ii) a data transfer tool to handle interpolation between non-matching meshes; and (iii) the mesh-reset and estimator-combine strategy to ensure a more economical number of degrees of freedom in each adaptive cycle. 
		Together, these components give rise to our main proposal: a novel neural network-enhanced $hr$-AFEM. 
		
		\subsection{Mesh size field}\label{sec:meshsize}
		The first objective is to reduce the total number of adaptive iterations. 
		In \cite{Xiao2026}, an efficient adaptive finite element method was proposed for elliptic equations by employing a refinement strategy that can control the number of mesh vertices, ensuring that the adaptive process terminates within seven iterations. 
		In this section, we extend that strategy to the adaptive finite element method for the parabolic equation \eqref{model}. 
		To materialize this adaptive refinement, a robust and controllable mesh generation tool is required. 
		Several established libraries such as Gmsh, MeshLab, TetGen, and NETGEN provide control, either directly or indirectly, to meet these requirements. 
		Among them, Gmsh \cite{Geuzaine2009} is a widely adopted and well-documented option.  
		Its key advantage is that refinement is driven by a user-supplied mesh size field, which determines the local mesh density. 
		Given its widespread acceptance and precise controllability, we employ Gmsh as the primary mesh generation tool in this work. 
		
		The mesh size field is a scalar function defined over the computational domain that prescribes the target local mesh size. 
		It serves as the direct link between the error estimator and the resulting mesh geometry: a smaller prescribed size in a region leads to a denser mesh there, and vice versa. 
		In our framework, the size field is designed to strategically concentrate mesh density in regions with high error contribution by identifying nodes where the error estimator exceeds a threshold. 
		Prior to introducing the algorithm, we present the necessary notation for the Hadamard product and its related operations as follows. 
		\begin{definition}
			Given matrices $A, B \in \mathbb{R}^{m \times n}$ with elements $A_{ij}$ and $B_{ij}$, their Hadamard product is defined as:
			\[ A\odot B \in \mathbb{R}^{m \times n}, \quad \text{with elements } (A\odot B)_{ij} = A_{ij}B_{ij}. \]
			For any integer $s \geq 1$, the Hadamard power of $A$, denoted as $A^{\circ s} \in \mathbb{R}^{m \times n}$, is given by:
			\[ A^{\circ s} = \underbrace{A \odot A \odot \cdots \odot A}_{s \text{ times}}, \quad \text{with elements } (A^{\circ s})_{ij} = A_{ij}^s. \]
			Similarly, the Hadamard division is defined as:
			\[ A\oslash B \in \mathbb{R}^{m \times n}, \quad \text{with elements } (A\oslash B)_{ij} = \frac{A_{ij}}{B_{ij}}. \]
		\end{definition}
		Then we provide \Cref{alg:SIZE} to generate such a mesh size field for Gmsh. 
		
		\normalem
		\begin{algorithm}
			\caption{Mesh size field generator for Gmsh}
			\label{alg:SIZE}
			\KwIn{Element error estimators $\{\eta_{h, K}\}$, element average edge size $\{h_e(K)\}$, iteration index $iteRO$, number of vertices $N_v$, number of elements $N_K$, dimension $d$, and mark ratio $\theta_r:=\{\text{number of marked vertices}\}/N_v$}
			\KwOut{Mesh size field $Size$}
			Construct the sparse matrix $node2cell \in \mathbb{R}^{N_v \times N_K}$, where each row corresponds to a node and each column represents an element. 
			The entry $(i, j)$ in this matrix is 1 if node $i$ belongs to element $j$, and 0 otherwise\;
			Compute the average edge length vector $h_v$ and error vector $E_v$ as follows:
			\[h_v = \left(node2cell \cdot h_e(K) \right) \oslash (node2cell \cdot \mathbf{1}_{N_K \times 1}),\]  
			\[E_v = (node2cell \cdot \eta_{h, K}) \oslash (node2cell \cdot \mathbf{1}_{N_K \times 1}),\]
			where $\mathbf{1}_{N_K \times 1}$ is ones vector\;
			Compute node density vector $\rho$ by:
			\[\rho = E_v^{\circ 2} \oslash h_v^{\circ d};\]\\
			Find the smallest set $M=\{i_1,i_2,\cdots,i_m\}$ such that:
			\[\sum_{i \in M} \rho[i] \geq \theta_r \sum_{i=1}^{N_v} \rho[i];\]\\
			Define the scaling factor vector $Scale$ as follows:
			\[
			Scale[i] =
			\begin{cases} 
				1, & i \notin M, \\
				\left(\sqrt[d]{1/2}\right)^{\log_2\left(N_v/m + 1\right)}, & i \in M;
			\end{cases}
			\]\\
			Compute the mesh size field by
			\[ Size = h_v \odot Scale^{\circ\;iteRO}. \]
		\end{algorithm}
		
		\subsection{Neural networks}\label{sec:neuralnetworks}
		The meshes generated by Gmsh with the user-specified size field are unstructured and non-nested, making the implementation of the interpolation $\Pi_h^nu_h^{n-1}$ both algorithmically complex and computationally expensive.  
		As the second key tool, we employ a neural network-based mesh‑free surrogate to construct $u_h^{n-1}$. 
		Notably, unlike the mesh‑dependent interpolation $\Pi_h^n u_h^{n-1}$ which must be recomputed after each mesh change, the neural network surrogate is mesh‑free and trained only once per time step. 
		
		\subsubsection{General neural networks}
		Consider a general $L$-layer neural network $u_{\theta}:\mathbb{R}^{d_0} \rightarrow \mathbb{R}^{d_L}$. 
		Formally, the $k$-th layer is defined as a function $\mathbf{N}^{(k)}(x)$ whose output is given by applying an element-wise activation function $\sigma$ to an affine transformation of the input. 
		Denoting the weight matrix and bias vector of this layer by $\mathbf{W}_k \in \mathbb{R}^{d_k \times d_{k-1}}$ and $\mathbf{b}_k \in \mathbb{R}^{d_k}$, respectively, the forward propagation is written as \[\mathbf{N}^{(k)}(\mathbf{x}^{(k-1)}):=\sigma(\mathbf{W}_k\mathbf{x}^{(k-1)}+\mathbf{b}_k).\]
		The activation function $\sigma$ is crucial, as it provides the necessary nonlinearity that enables the network to learn complex patterns and approximate a wide range of functions. 
		In the absence of such nonlinearity, the network would collapse into a sequence of linear transformations, severely restricting its ability to capture intricate patterns and generalize to unseen data. 
		For regression, no activation is used in the output layer to allow for unconstrained output. 
		Therefore, the general neural network is given by
		\begin{equation}\label{neuralnetwork0}
			u_{\theta}(\mathbf{x})=\mathbf{W}_L\left(\mathbf{N}^{(L-1)}\circ \cdots \circ \mathbf{N}^{(2)}\circ \mathbf{N}^{(1)}(\mathbf{x})\right)+\mathbf{b}_L. 
		\end{equation}
		
		While numerous activation functions, such as ReLU, sigmoid, and tanh, have been proposed in the literature, we employ the hyperbolic tangent function in this work,
		\[\sigma(x)=\tanh{x} = \frac{e^x-e^{-x}}{e^x+e^{-x}}.\]
		This choice is motivated by its smoothness, bounded range, and symmetry about zero, which are well suited for approximating the smooth solution profiles commonly encountered in partial differential equations.
		\begin{remark}
			Preliminary experiments suggest that, for the class of problems considered, the $\tanh$ activation function generally provides more stable and accurate approximations than ReLU. Alternative activation functions (e.g., sinusoidal activations) may achieve superior performance in certain specific scenarios. However, the primary objective of this work is to validate the proposed adaptive framework rather than to conduct a systematic comparison of activation functions. Therefore, we adopt the widely used and robust $\tanh$ activation function throughout, while noting that other choices may be incorporated to better accommodate problem-specific characteristics.
		\end{remark}
		
		In general, the training of the neural network proceeds as follows.
		Let $\theta=((\mathbf{W}_0, \mathbf{b}_0), \cdots, (\mathbf{W}_L, \mathbf{b}_L))$ represent all trainable parameters of the network. 
		Given a function $u\in C^0(\Omega)$, we train a neural network $u_{\theta}(\mathbf{x})$ that approximates $u$ using the dataset $\{(\mathbf{x}_i, u(\mathbf{x}_i)\}_{i=1}^N$. The network parameters $\theta$ are optimized by minimizing the empirical loss function $\mathcal{L}(u_{\theta})$, typically the mean squared error (MSE):
		\[\theta^*=arg \min\limits_{\theta\in \Theta}\mathcal{L}(u_{\theta}),\quad \mathcal{L}(u_{\theta}):=\frac{1}{N}\sum\limits_{i=1}^N(u_{\theta}(\mathbf{x}_i)-u(\mathbf{x}_i))^2,\]
		where $\Theta$ is the feasible set consisting of the parameters and hyperparameters of the neural network.
		
		\subsubsection{Mesh-free surrogate via neural networks}
		The objective here is to construct a neural network surrogate $u_{\theta}^{n-1}$, based on the general neural network \eqref{neuralnetwork0}, to approximate the finite element solution $u_h^{n-1}$ obtained at $t_{n-1}$. This surrogate provides a continuous, mesh-free representation of $u_h^{n-1}$ for use in the subsequent adaptive procedure at time $t_n$, thereby eliminating the need for the computationally expensive interpolation $\Pi_h^n u_h^{n-1}$ in \eqref{full-discrete}.  
		The surrogate $u_{\theta}^{n-1}$ is trained using an input–output dataset constructed from the vertices of $\mathcal{T}_h^{n-1}$ and the corresponding finite element approximations, namely
		$$
		\{ (\mathbf{x}, u_h^{n-1}(\mathbf{x})) \mid \mathbf{x} \in \mathcal{N}_h^{n-1} \},
		$$
		where $\mathcal{N}_h^{n-1}$ denotes the set of vertices of the mesh $\mathcal{T}_h^{n-1}$.
		
		Motivated by \cite{Sukumar2022}, the Dirichlet boundary condition $u|_{\partial \Omega} = g(\mathbf{x},t), \ \mathbf{x} \in \partial \Omega$ is imposed by constructing the neural network surrogate $u_{\theta}^{n-1}$ as
		\begin{equation}\label{neuralnetwork}
			u_{\theta}^{n-1}(\mathbf{x}) = d(\mathbf{x}) u_{\theta}(\mathbf{x}) + \tilde{g}(\mathbf{x},t_{n-1}),
		\end{equation}
		where $u_{\theta}(\mathbf{x})$ denotes the general neural network defined in \eqref{neuralnetwork0}. 
		The function \(d(\mathbf{x})\) is assumed to be continuous and satisfies
		\begin{equation}\label{dxcond}
			d(\mathbf{x})=0 \quad \text{for } \mathbf{x}\in\partial\Omega,
			\qquad
			d(\mathbf{x})\neq 0 \quad \text{for } \mathbf{x}\in\Omega,
		\end{equation}
		A natural candidate is to adopt the exact distance function,
		\begin{equation}\label{dxformula}
			d(\mathbf{x})
			=\min\limits_{\mathbf{y}\in\partial\Omega}|\mathbf{x}-\mathbf{y}|.
		\end{equation}
		Nevertheless, encoding the exact distance function \eqref{dxformula} in neural networks entails considerable computational overhead \cite{Sukumar2022}. In literature, an approximate distance function $d(\mathbf{x})$ is usually employed. 
		In our application, $d(\mathbf{x})$ is not restricted to the approximate distance function. Any continuous function satisfying \eqref{dxcond} is sufficient to enforce the Dirichlet boundary condition exactly through the network construction.
		The specific choice of \(d(\mathbf{x})\) utilized in the numerical examples is provided in \Cref{app:func}.
		Furthermore,  $\tilde{g}(\mathbf{x},t_{n-1}) \in H^1(\Omega) \cap C^0(\Omega)$ is an extension of the boundary data $g(\mathbf{x},t_{n-1})$ from $\partial\Omega$ to the entire domain $\Omega$, satisfying $\tilde{g}|_{\partial \Omega} = g$. 
		

		For a homogeneous Dirichlet boundary condition $u|_{\partial \Omega} = 0, \ \mathbf{x} \in \partial \Omega$, the neural network surrogate $u_{\theta}^{n-1}$ in \eqref{neuralnetwork} reduces to
		\begin{equation}\label{neuralnetwork+}
			u_{\theta}^{n-1}(\mathbf{x}) = d(\mathbf{x}) u_{\theta}(\mathbf{x}).
		\end{equation}
		As a simple illustration, consider a two-point boundary value problem on the interval $[a, b]$ with nonhomogeneous Dirichlet boundary conditions $u(a, t_{n-1})=A$ and $u(b, t_{n-1})=B$. In this case, one may choose $\tilde{g}(x, t_{n-1}) = \frac{b-x}{b-a}A + \frac{x-a}{b-a}B$, which is simply the linear interpolant of the prescribed boundary data.  
		
		
		Then, the empirical loss function is defined as MSE:
		\begin{equation}\label{lossFunction}
			\mathcal{L}\!\left(u_{\theta}^{n-1}; \mathcal{T}_h^{n-1}\right)
			= \frac{1}{\lvert \mathcal{T}_h^{n-1} \rvert}
			\sum_{\mathbf{x} \in \mathcal{N}_h^{n-1}}
			\left| u_{\theta}^{n-1}(\mathbf{x}) - u_h^{n-1}(\mathbf{x}) \right|^2,
		\end{equation}
		where $\lvert \mathcal{T}_h^{n-1} \rvert$ denotes the number of vertices in $\mathcal{T}_h^{n-1}$. 
		
		The optimal surrogate $u_{\theta^*}^{n-1}(\mathbf{x})= d(\mathbf{x}) u_{\theta^{n-1}_*}(\mathbf{x}) + \tilde{g}(\mathbf{x},t_{n-1})$ is obtained by solving
		\[
		\theta_*^{n-1} = \arg\min_{\theta} \mathcal{L}\!\left(u_{\theta}^{n-1}; \mathcal{T}_h^{n-1}\right).
		\]
		The resulting minimizer $u_{\theta^*}^{n-1}(\mathbf{x})$ provides an accurate approximation of $u_h^{n-1}(\mathbf{x})$ and therefore serves as a mesh-free surrogate for $u_h^{n-1}$ over the entire computational domain.
		
		\begin{remark} 
			The same neural network architecture is employed for the general neural network $u_\theta$ in \eqref{neuralnetwork} across all time steps.
			The time-marching strategy in \eqref{full-discrete} introduces a strict serial dependency; therefore, the neural network surrogates are trained sequentially in time to properly capture the singular structures that may arise at each time level. 
			The neural network hyperparameters in this study, such as the network depth, neurons per layer, and learning rate, were chosen based on empirical experimentation. 
			The parameters $\theta$ of $u_\theta$ are randomly initialized and thoroughly trained only for the initial surrogate $u_\theta^0$. 
			For subsequent time steps, training is warm-started from the converged parameters of the previous step and involves only brief fine-tuning, rather than reinitialization. 
			Our experiments imply that achieving sufficient accuracy at each new time step requires only a few optimization iterations (e.g., \Cref{Ex2D1-Iter}), with the loss typically converging to expected MSE (e.g., $10^{-6} $ to $ 10^{-4}$). 
			The adequacy of this loss error can be assessed from the history of the gradient error (e.g., \Cref{Ex2-order}, \Cref{Ex3-order}): once quasi-optimal convergence is observed, we conclude that the achieved loss is sufficiently small for the intended accuracy.
			For neural networks trained using gradient-based optimizers such as ADAM or L-BFGS, this error level generally represents a numerical optimization plateau, beyond which further loss reduction yields diminishing returns. Further investigation into reducing the MSE will be left for future work. 
		\end{remark}
		
		\begin{remark}
			We emphasize that the primary purpose of \(d(\mathbf{x})\) is to vanish on the boundary while remaining nonzero in the interior of the domain, as specified in \eqref{dxcond}. In practice, such functions can often be constructed explicitly. Consequently, the cost of evaluating \(d(\mathbf{x})\) at collocation points is negligible, since it reduces to the pointwise evaluation of a known analytic expression.
			Moreover, the choice of $d(\mathbf{x})$ is not unique. For simple geometries, such as squares and cubes, $d(\mathbf{x})$ can be chosen as the exact distance function
			\(
			d(\mathbf{x}) = \min_{\mathbf{y} \in \partial \Omega} |\mathbf{x} - \mathbf{y}|,
			\)
			which admits a closed-form representation. For example, on the unit square, 
			\(
			d(\mathbf{x}) = \min(x, 1-x, y, 1-y).
			\)
			For more general geometries, $d(\mathbf{x})$ can be chosen as the approximate distance function, which relies on geometry-dependent techniques, such as domain decomposition, cut-off functions, R-functions \cite{rvachev1995r}, or other approximate distance-function constructions. In particular, the L-shaped domain considered in this work, $d(\mathbf{x})$ is obtained by decomposing the domain into two rectangular subdomains and combining the corresponding functions through an R-function formulation. We refer the readers to \cite{Sukumar2022} for a comprehensive discussion of approximate distance functions on general domains.
			
		\end{remark}

		By integrating the optimized neural network surrogate $u_{\theta^*}^{n-1}$ into the fully discrete finite element scheme \eqref{full-discrete}, we obtain the following reformulated variational problem.
		\begin{proposition}
			The neural network–enhanced finite element scheme seeks $u_h^n \in V_h^n$ such that
			\begin{equation}\label{fsna}  
				\left(\frac{u_h^n-u_{\theta^*}^{n-1}}{\tau},v\right)+(\nabla{u_h^n},\nabla v ) = \left(f_h^{n},v\right), \quad  \forall v \in V_h^n.
			\end{equation}
		\end{proposition}
	
	\begin{lemma}
		For any given mesh, source term, and boundary conditions, and for any $\tau > 0$, the neural network–enhanced finite element scheme \eqref{fsna} admits a unique solution. Moreover, the following stability estimate holds:
		\begin{equation}\label{stab}
			\|u_h^n - \tilde{u}_h^n\|
			\le 
			\|u_{\theta^*}^{\,n-1} - \tilde{u}_{\theta^*}^{\,n-1}\|,
		\end{equation}
		where $u_h^n$ and $\tilde{u}_h^n$ denote the solutions corresponding to different neural network surrogates $u_{\theta^*}^{\,n-1}$ and $\tilde{u}_{\theta^*}^{\,n-1}$, respectively.
	\end{lemma}
	The stability estimate \eqref{stab} shows that the scheme is $L^2$-stable for any $\tau>0$. Consequently, if the neural network surrogate contains only a small perturbation, then this perturbation does not amplify through the time-stepping procedure, and the resulting solution remains stable. 
	Although the time step size $\tau$ does not affect the unconditional stability of \eqref{fsna}, it has a significant impact on both the solution accuracy and the quality of the neural network initialization $u_{\theta^*}^{n-1}$. A detailed numerical investigation is presented in \Cref{example1}. 
	
	
	\subsection{$hr$-adaptive finite element algorithm}
	To effectively solve the neural network–enhanced finite element scheme \eqref{fsna}, we require an adaptive algorithm equipped with a reliable a posteriori error estimator to guide mesh refinement.
	Note that the accuracy of the term $(u_{\theta^*}^{n-1},v)$ in \eqref{fsna} may also affect the accuracy of the numerical solution $u_h^n$. 
	Consequently, the mesh $\mathcal{T}_h^n$ at time $t_n$ must be carefully generated to capture the key features of both the previous approximation $u_{\theta^*}^{n-1}$ and the current solution $u_h^n$, since the locations and intensities of singularities may vary from one time step to the next.
	Therefore, the error estimator should account for contributions from both $u_h^n$ and $u_{\theta^*}^{n-1}$, since both play a role in determining an appropriate mesh. 
	In light of this, we employ a combined error indicator in \eqref{estimators} as
	\begin{equation}\label{etarec}
		\eta_{h, K}^n=\max\{\hat{\eta}_{h, K}^n,\hat{\eta}_{h, K}^{n-1}\}=\frac{\hat{\eta}_{h, K}^n+\hat{\eta}_{h, K}^{n-1}+|\hat{\eta}_{h, K}^n-\hat{\eta}_{h, K}^{n-1}|}{2}, \qquad \eta_h^n=\left(\sum\limits_{K\in\mathcal{T}_h^n}\eta_{h, K}^2\right)^{1/2},
	\end{equation}
	where $\hat{\eta}_{h, K}^n$ is defined in \eqref{estimators}. 
	The estimator \eqref{etarec} incorporates information from both time levels and can guide the new mesh to capture the singularities at $t_n$ while ensuring accurate computation of the term $(u_{\theta^*}^{n-1},v)$. 
	
	In contrast to the coarsening procedure commonly used in conventional AFEM, the final mesh at $t_{n-1}$ is not inherited in our algorithm.
	Instead, the adaptation process employs a mesh‑reset mechanism, whereby the mesh is re‑initialized from a predefined coarse mesh $\mathcal{T}_h^{n,0}$ at each time level $t_n$. 
	As a result, the number of degrees of freedom remains controlled throughout the simulation, leading to reduced computational cost.
	
	At each time level, the adaptive procedure proceeds systematically.
	Starting from the initial mesh $\mathcal{T}_h^{n,0}$, five adaptive iterations are performed. 
	Based on the reliable data set $\{(\eta_h^{n, i}, \lvert\mathcal{T}_h^{n, i}\rvert)\}_{i=2}^4$, a least-squares fitting is used to approximate the model between the error estimator and the NOV 
	\begin{equation}\label{fitModel}
		\eta_h^n \approx cN^{-p}.
	\end{equation}
	With the fitted parameters $(c, p)$, the desired NOV is predicted as 
	$$N = \left\lceil\sqrt[p]{c/eTol}\right\rceil$$ 
	to ensure $\eta_h^n \leq eTol$.
	The predicted $N$, together with the local error estimators, is then used to compute the mesh size field via \Cref{alg:SIZE}. 
	Subsequently, a new mesh $\mathcal{T}_h^{n,5}$ is generated using Gmsh according to the computed size field, and equation \eqref{fsna} is solved on this mesh. 
	Finally, a seventh adaptive iteration is performed on $\mathcal{T}_h^{n,6}$ as a safeguard to verify the prescribed tolerance. 
	
	\begin{remark}
		The rationale for the seven-step mesh adaptation strategy is as follows (see \Cref{Exremark} for illustration).
		First, the first two iterative solutions computed on the initial mesh $\mathcal{T}_h^{n,0}$ (given) and the first refined mesh $\mathcal{T}_h^{n,1}$ are excluded from data fitting to eliminate transient start-up effects. During this stage, the mesh undergoes rapid changes, resulting in an unstable and non-asymptotic relationship between the error estimator and the NOV; the error may even increase temporarily under refinement.
		Second, after the start-up phase, at least three reliable data points on $\mathcal{T}_h^{n,i}$, $i=2,3,4$, are required to robustly fit the model \eqref{fitModel}. Note that a two-point fit can be highly sensitive to noise and lacks statistical reliability.
		Third, one adaptive iteration is carried out on the mesh $\mathcal{T}_h^{n,5}$ generated by Gmsh using the sizing function derived from the fitted model. To provide a rigorous guarantee of accuracy, an additional iteration is deliberately performed on $\mathcal{T}_h^{n,6}$. Even if the fitted model predicts that the prescribed tolerance will be achieved at the sixth step, executing the seventh iteration ensures that the final solution indeed satisfies the required error bound. Therefore, the required number of adaptive iterations is 
		\[
		2 \text{ (start-up)} + 3 \text{ (fitting)} + 1 \text{ (validation)} + 1 \text{ (assurance) } = 7 \text{(adaptive iterations)}.
		\]
	\end{remark}
	
	Based on the discussion above, given the approximation $u_h^{n-1}\in V_h^{n-1}$ at the time step $t_{n-1}$ and the selected error estimators, the workflow of the new $hr$-adaptive algorithm for the parabolic equation \eqref{model} at the $n$-th time step can be summarized as follows:
	\begin{center}
		\begin{tikzpicture}[node distance=0.5cm]        
			\node (node1) {$u_h^{n-1}$};
			\node[right=of node1] (node2) {\textit{LEARN}};
			\node[right=of node2] (node3) {\textit{SOLVE}};
			\node[right=of node3] (node4) {\textit{ESTIMATE}};
			\node[right=of node4] (node5) {\textit{SIZE}};
			\node[right=of node5] (node6) {\textit{GENERATE}};
			\node[right=of node6] (node7) {$u_h^n$};
			
			\draw[->] (node1) -- (node2);
			\draw[->] (node2) -- (node3);
			\draw[->] (node3) -- (node4);
			\draw[->] (node4) -- (node5);
			\draw[->] (node5) -- (node6);
			
			\draw[->] (node6.south) -- ++(0,-0.3) -| (node3.south);
			\draw[->] (node4.north) -- ++(0,0.3) -| (node7.north);  
		\end{tikzpicture}
	\end{center} 
	From the perspective of the adaptive workflow, the proposed method introduces three main modifications compared to the conventional one. Here, \textit{LEARN} represents the training of the neural network surrogate to transfer the solution between successive meshes; \textit{SOLVE} procedure solves the neural network–enhanced finite element scheme \eqref{fsna};
	\textit{SIZE} refers to the computation of mesh size field described in \Cref{sec:meshsize}; and 
	\textit{GENERATE} stands for creating a new adapted mesh based on that field. 
	Following this workflow, we present a neural network-enhanced $hr$-adaptive algorithm for the parabolic equation \eqref{model} in \Cref{alg:ParabolicHR}.

	\normalem
	\begin{algorithm}
		\caption{Neural network-enhanced $hr$-AFEM for the parabolic equation \eqref{model}}
		\label{alg:ParabolicHR}
		\KwIn{Domain $\Omega$, source term $f$, tolerance $eTol$, initial data $u_0$, time step size $\tau$}
		\KwOut{The meshes $\{\mathcal{T}_h^n\}_{n=1}^S$ and finite element approximations $\{u_h^n\}_{n=1}^S$}
		Generate an initial uniform mesh $\mathcal{T}_h^{initial}$ using Gmsh with a given mesh diameter\;
		\For {$n=0$ \KwTo $S$ }{
			Set $\mathcal{T}_h^n=\mathcal{T}_h^{initial}$, $t_n=t_0+n\tau$\;
			\For {$k=0$ \KwTo $6$}{
				Determine the number of the vertices $\lvert \mathcal{T}_h^{n} \rvert$ in the mesh $\mathcal{T}_h^n$\;
				\If {$n==0$}{
					(PROJECT) Compute the projection $u_h^0=Pu_0$\;
				}
				\Else { 
					(SOLVE) Solve problem \eqref{fsna} for $u^n_h$ on $\mathcal{T}_h^n$ using the data $u_{\theta^*}^{n-1}$\;
				}
				
				(ESTIMATE) Compute the recovery-based error estimators $\{\eta_{h, K}^{n, k}\}$ and $\eta_h^{n, k}$\;
				\If {$\eta_h^{n, k} \leq eTol$ \textbf{or} $k==6$ }{
					Break \;
				}
				\Else {          
					\If {$k==4$}{
						Calculate the parameters $(c, p)$ in $\eta_h^n=cN^{-p}$ by least square fitting with data $\{(\eta_h^{n, i}, N_i)\}_{i=2}^4$ \;
						Calculate $N=\left\lceil\sqrt[p]{c/eTol}\right\rceil$\;
						Set $iteRO=\max\left(\left\lceil\log_2\left(\frac{N}{N_{k-1}}\right)\right\rceil, 1\right)$ \;
					}
					\Else {
						Set $iteRO=1$ \;
					}
				}
				
				(SIZE) Compute the mesh size field $F_n$ by \Cref{alg:SIZE} based on the local error estimators $\{\eta_{h,K}^{n,k}\}$ and $iteRO$\;
				(GENERATE) Update the mesh $\mathcal{T}_h^n$ using Gmsh based on the mesh size field $F_n$ \;
			}
			
			(LEARN) Train a neural network $u_{\theta}^{n}$ to learn the finite element solution $u^{n}_h$, and obtain the minimizer as $u_{\theta^*}^{n}$.
		}
	\end{algorithm} 
	
	Compared with \Cref{alg:Parabolic1}, the main changes in \Cref{alg:ParabolicHR} are as follows. The generation of the initial mesh $\mathcal{T}_h^0$ is now performed within the time-stepping loop. The adaptive process is limited to a maximum of seven iterations per time step. The \textit{SIZE} and \textit{GENERATE} steps replace the \textit{REFINE} and \textit{INTERPOLATE} steps from \Cref{alg:Parabolic1}, although their operations are fundamentally different. Similarly, the \textit{LEARN} step takes the place of \textit{COARSEN}, but functions in a distinct manner. 
	
	\begin{example}\label{ex:alg3}
		This example repeats the experiment from \Cref{exp2d} using \Cref{alg:ParabolicHR} under the same configuration. 
		Multiple characteristics and conclusions from the prior analysis are thereby validated.
		\Cref{Exremark} presents the adaptive meshes generated by \Cref{alg:ParabolicHR} at $t_1=0.1$. 
		The mesh sequence shows two distinct regions being refined simultaneously, corresponding precisely to the singularities of the finite element solutions $u_h^0$ and $u_h^1$. 
		Without $\hat{\eta}_{h, K}^{n-1}$ in \eqref{etarec}, refinement would concentrate only on regions of large error for $u_h^1$. 
		This demonstrates the effectiveness and necessity of the combined error estimator.
		The adaptive process terminates successfully after 7 iterations once the prescribed tolerance $eTol$ is met, with the NOV progressively increasing as $98$, $163$, $322$, $633$, $1277$, $12649$, and $25764$. 
		A notable jump occurs after the fifth refinement, where the NOV rises sharply from $1277$ to $12649$. 
		This sharp increase is driven by the least‑squares fit, which predicts the NOV required to meet the tolerance and is then passed to the mesh size field computation, leading to the generation of a much denser mesh and the observed jump in NOV. 
		The execution of a seventh adaptive iteration confirms that the mesh $\mathcal{T}_h^{n,5}$, generated using the NOV predicted by the least‑squares fit, did not yet satisfy the tolerance $eTol$, which unequivocally highlights the critical importance of the seventh iteration in mesh adaptation strategy.  \begin{figure}
			\begin{minipage}[c]{0.25\linewidth}
				\raggedleft
				\includegraphics[width=0.8\linewidth]{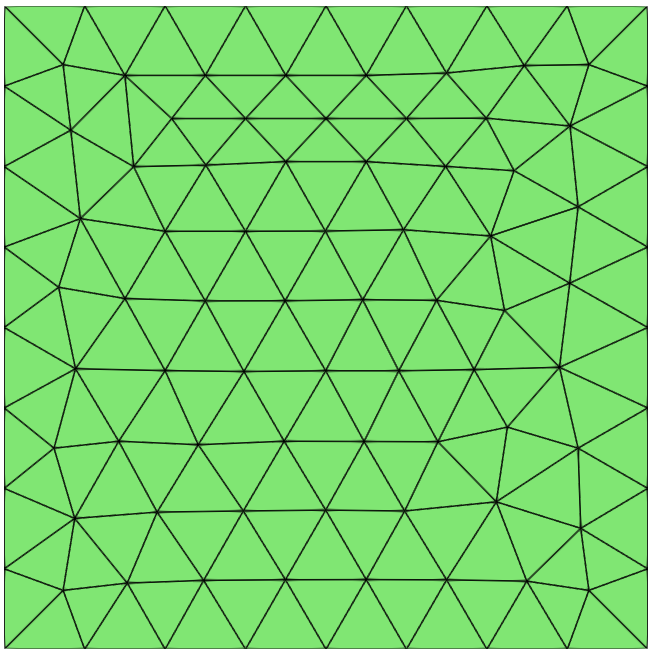}
			\end{minipage}%
			\begin{minipage}[c]{0.7\linewidth}
				\centering
				\includegraphics[width=0.3\textwidth]{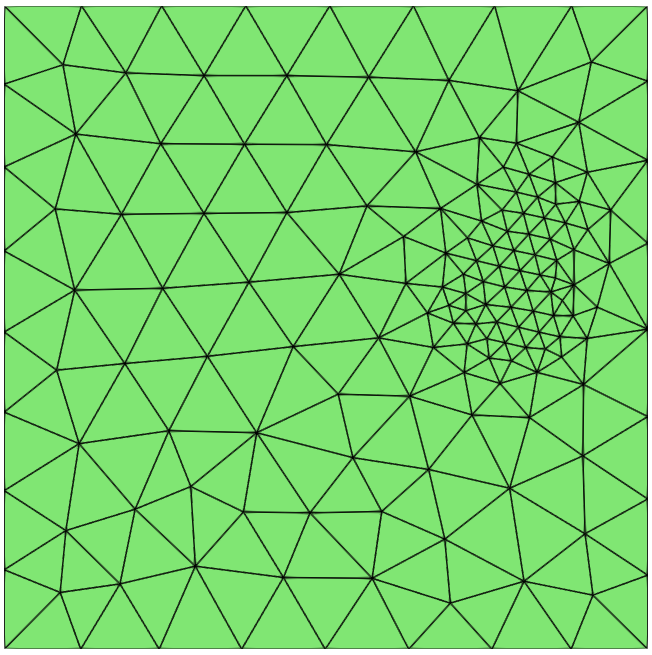}
				\includegraphics[width=0.3\textwidth]{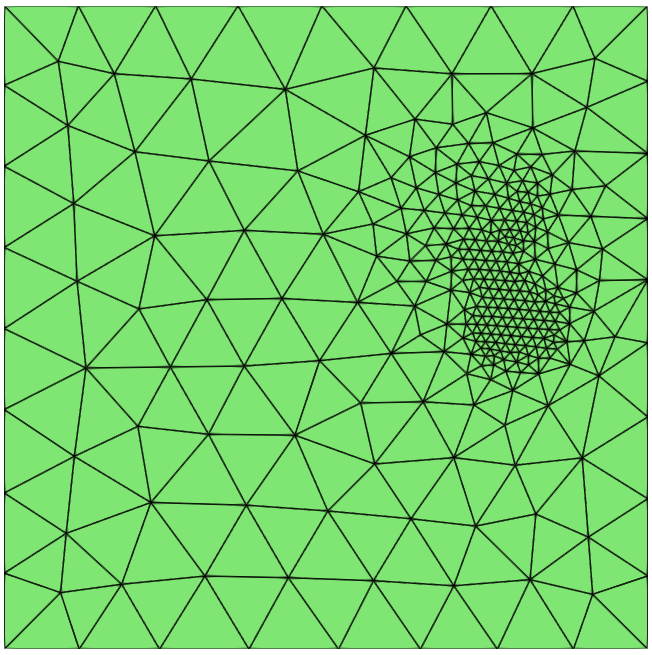}
				\includegraphics[width=0.3\textwidth]{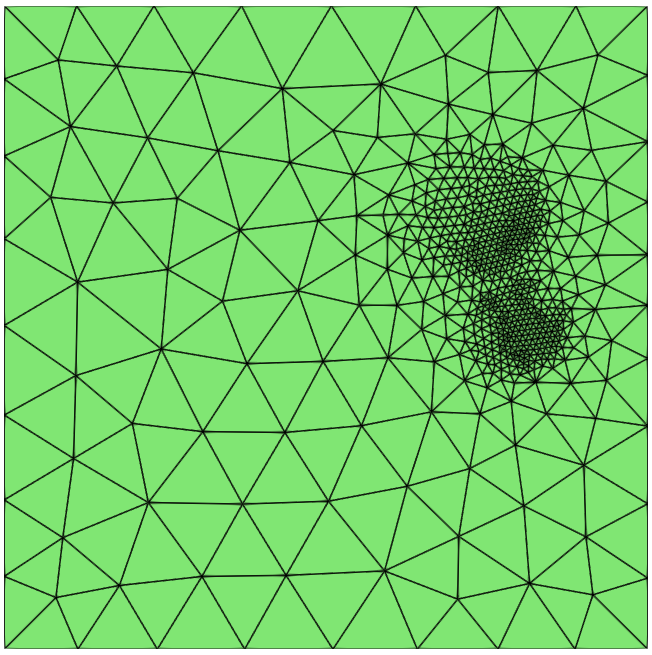}\vspace{0.1cm}
				\includegraphics[width=0.3\textwidth]{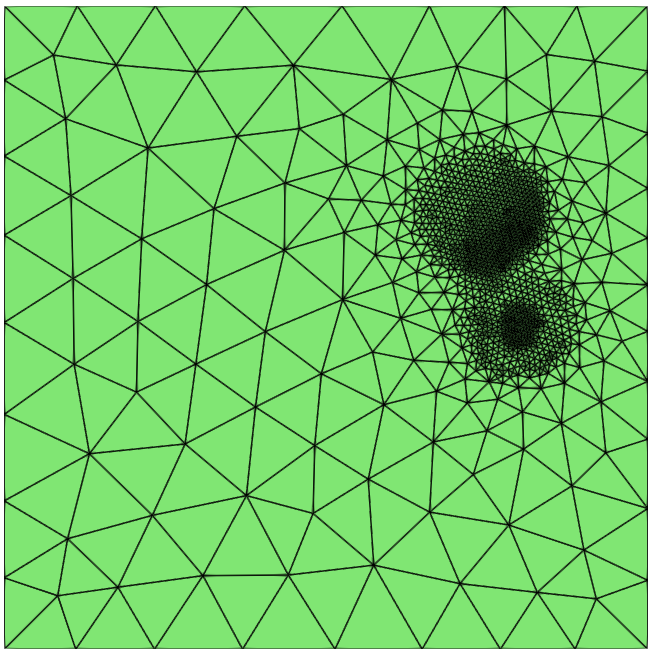}
				\includegraphics[width=0.3\textwidth]{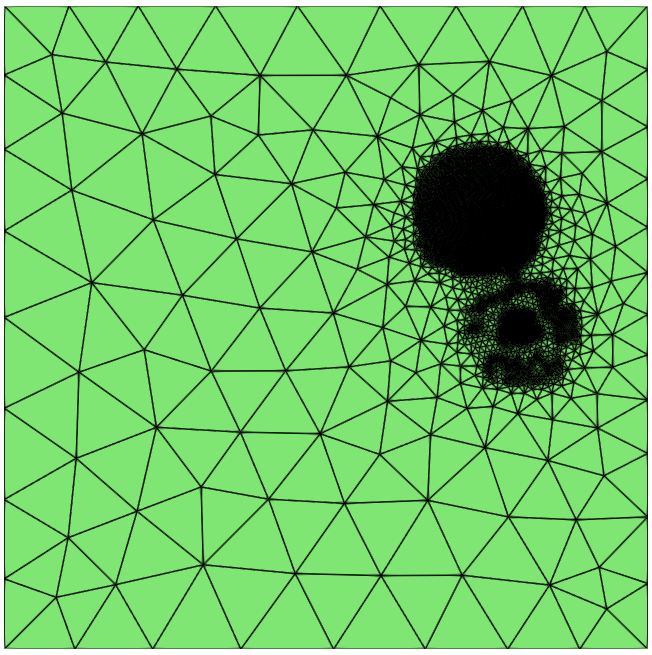}
				\includegraphics[width=0.3\textwidth]{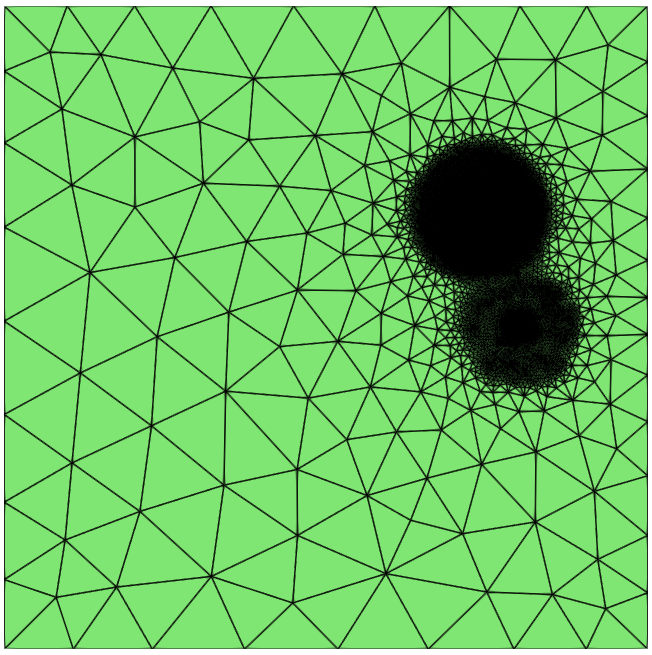}
			\end{minipage}
			\caption{
				\Cref{ex:alg3} using \Cref{alg:ParabolicHR}: Initial mesh $\mathcal{T}^{1,0}_h$ (98 nodes) at $t_1=0.1$ (left); and adaptively refined meshes $\mathcal{T}^{1,1}_h$ (163), $\mathcal{T}^{1,2}_h$ (322), $\mathcal{T}^{1,3}_h$ (633), $\mathcal{T}^{1,4}_h$ (1277), $\mathcal{T}^{1,5}_h$ (12649), $\mathcal{T}^{1,6}_h$ (25764).}
			\label{Exremark}
		\end{figure} 
	\end{example}
	
	\section{Numerical Examples}\label{numExp}
	This section presents numerical experiments to demonstrate the efficiency of the neural network-enhanced $hr$-adaptive \Cref{alg:ParabolicHR} in solving the parabolic equation in both 2D and 3D. 
	The main purpose of the tests centers on three primary test cases which feature distinct moving singularities—rotation, diffusion, and splitting—to validate the algorithm’s performance in tracking dynamic sharp gradients. 
	These 2D cases are then extended to 3D to assess scalability in resolving such singularities and show performance beyond planar geometries.
	Finally, two numerical examples, namely the Allen–Cahn equation and a parabolic problem defined on a domain with reentrant corners, are provided to further illustrate the robustness and applicability of the proposed adaptive algorithm to nonlinear parabolic equations and domains with geometric singularities.
	All examples are implemented using the \texttt{fealpy} library (version 1.2.0) \cite{fealpy}. 
	
	In the following examples, unless stated otherwise, the conforming linear finite element is used in all numerical experiments, with the time step size $\tau=0.01$. 
	The iteration tolerance $eTol$ is employed to control adaptive termination, with specific values provided in each respective example. 
	It primarily affects the computation of $iteRO$ in the sixth adaptation and determines whether a seventh adaptation is necessary. 
	For neural network training, weights are initialized via the Kaiming method with biases initialized to zero. 
	The $\tanh$ function serves as the activation function for the hidden layers, while the output layer uses a linear activation with no bias term ($b_L=0$). 
	The L-BFGS optimizer is employed to minimize the loss function \eqref{lossFunction} to ensure convergence (stagnation). The neural network hyperparameters, including the number of layers, the number of neurons per layer, and the learning rate, are selected through empirical tuning rather than a systematic grid search. 
	The distance function $d(\mathbf{x})$ and the continuous extension 
	$\tilde{g}(\mathbf{x}, t_{n-1})$ of the boundary data used in the neural 
	network \eqref{neuralnetwork} for all numerical examples are provided in 
	\Cref{app:func}. 
	\subsection{2D Example} 
	The neural network \eqref{neuralnetwork} used in this section employs an architecture $[2, 40, 40, 40, 1]$, comprising two input neurons, three hidden layers of 40 neurons each, and one output neuron. 
	It is adopted as a demonstration case and has undergone no optimization or parameter tuning.    
	\begin{figure}
		\centering
		\includegraphics[width=0.3\linewidth]{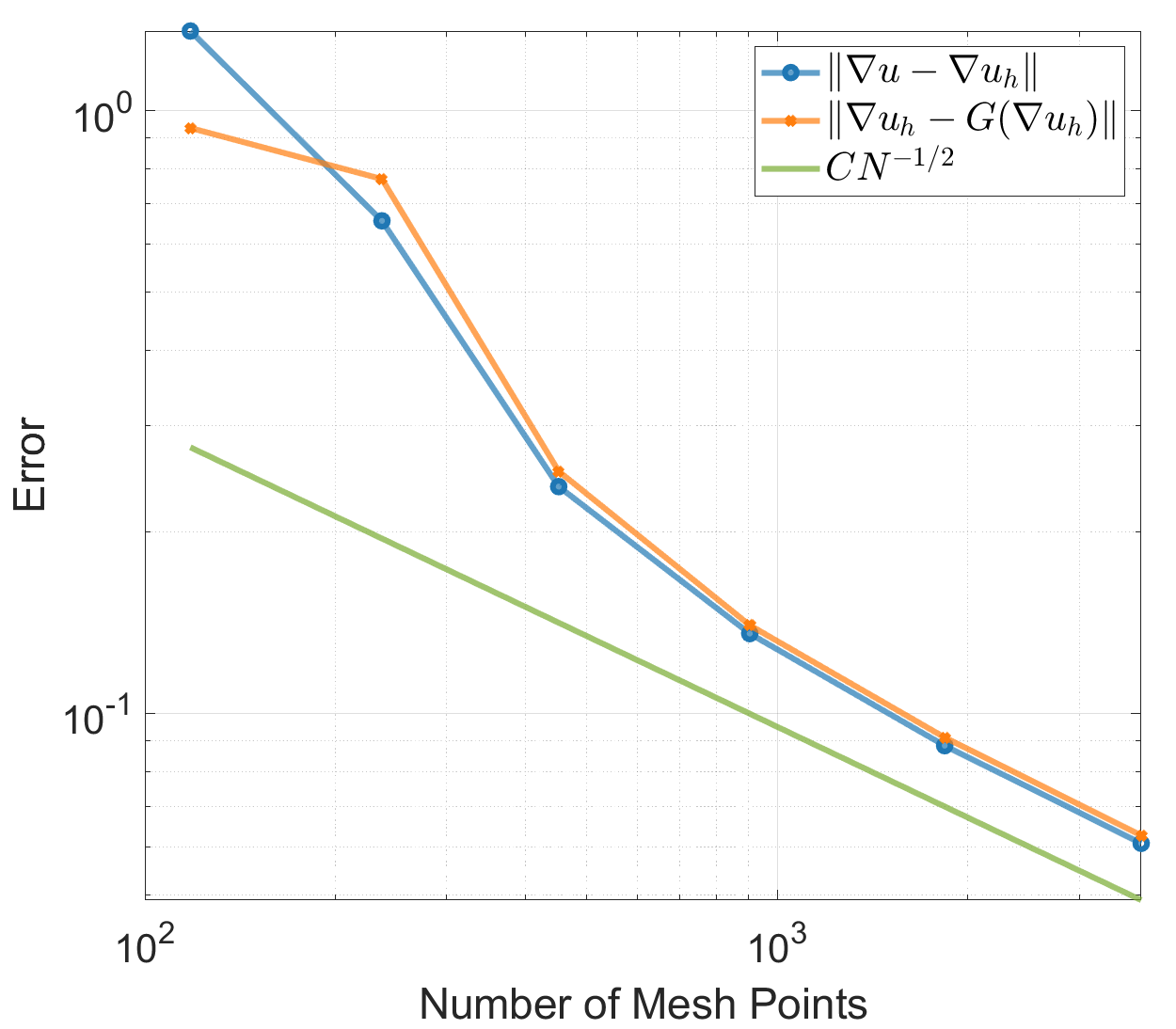}
		\includegraphics[width=0.3\linewidth]{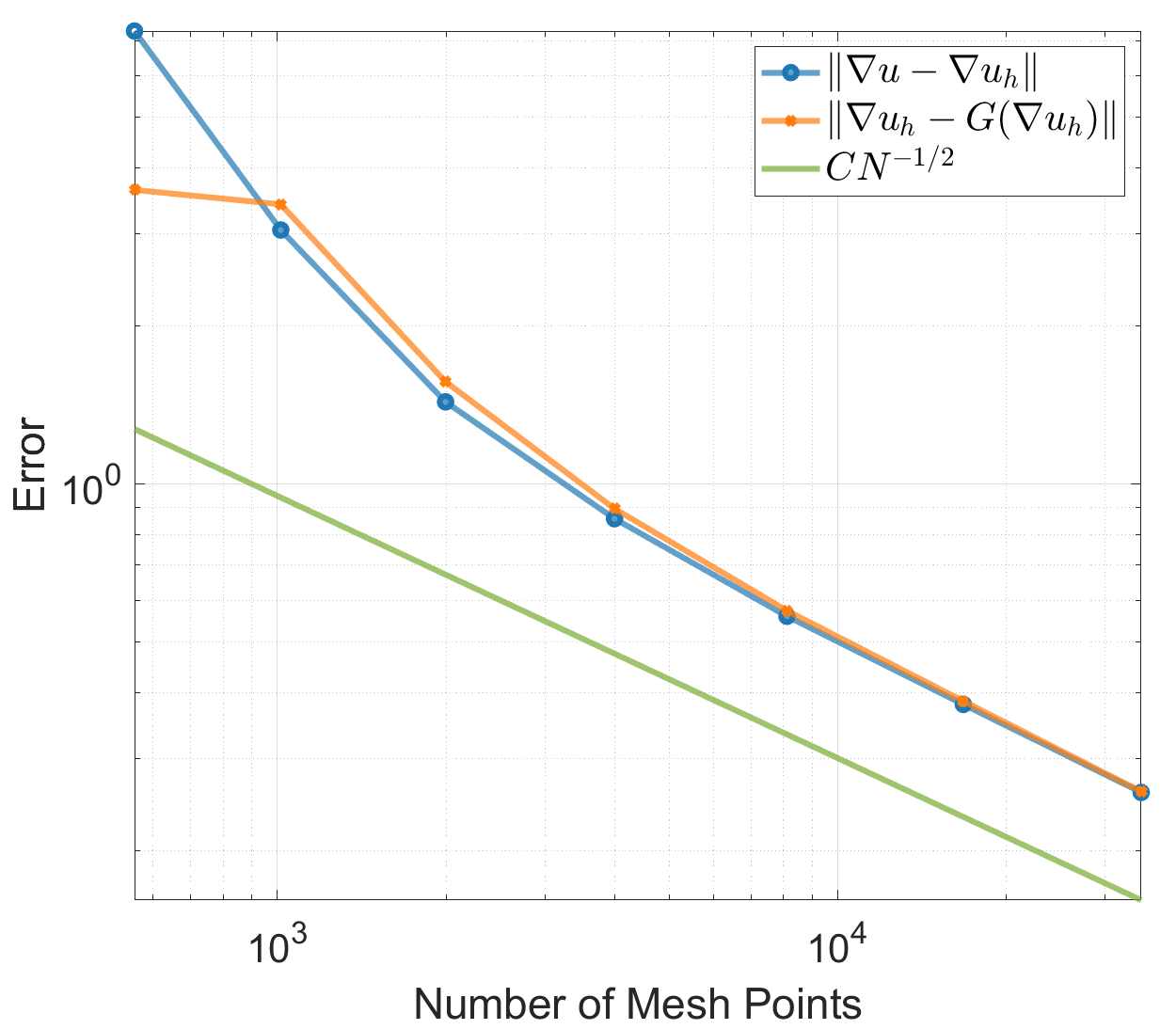}
		\includegraphics[width=0.3\linewidth]{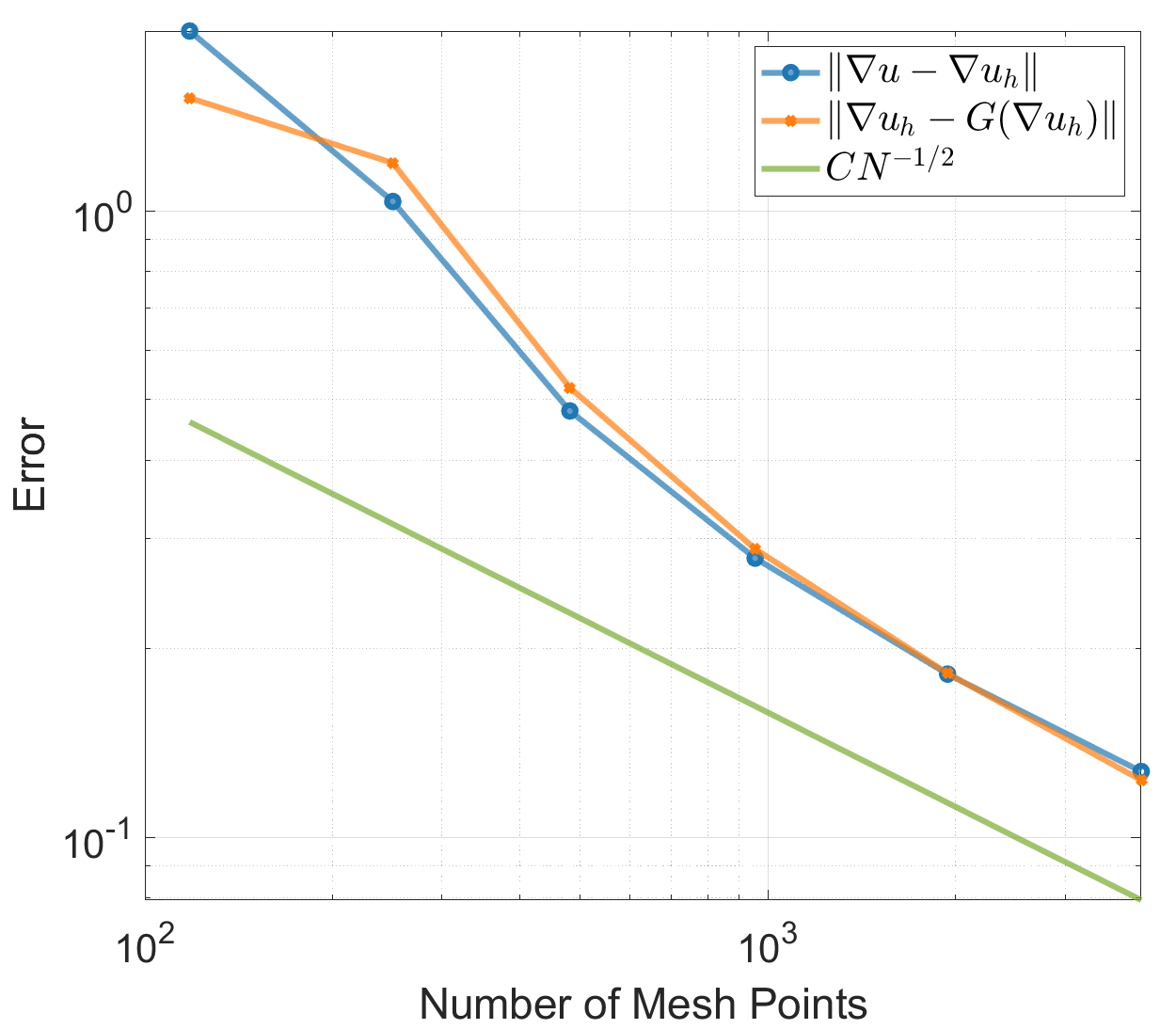}
		\caption{History of gradient error $\Vert\nabla u-\nabla u_h\Vert$ and recovered error estimator $\Vert\nabla u_h-G(\nabla u_h)\Vert$ at $T=1.0$ in \Cref{example1} (left), \Cref{example2} (middle) and \Cref{example3} (right).}
		\label{Ex2-order}
	\end{figure}
	
	\begin{example} \label{example1}
		(Rotation) Consider the 2D parabolic equation \eqref{model} in the domain \(\Omega = [-1,1]^2\), with a source function \( f \) chosen such that the exact solution \( u \) is 
		\[
		u(x,y,t) = \exp \big( -500 (x - 0.3\cos(2\pi t))^2 \big)  
		\exp \big( -500 (y - 0.3\sin(2\pi t))^2 \big).
		\]
	\end{example}
	
\begin{figure}[!ht]
		\centering
		\includegraphics[width=0.85\linewidth]{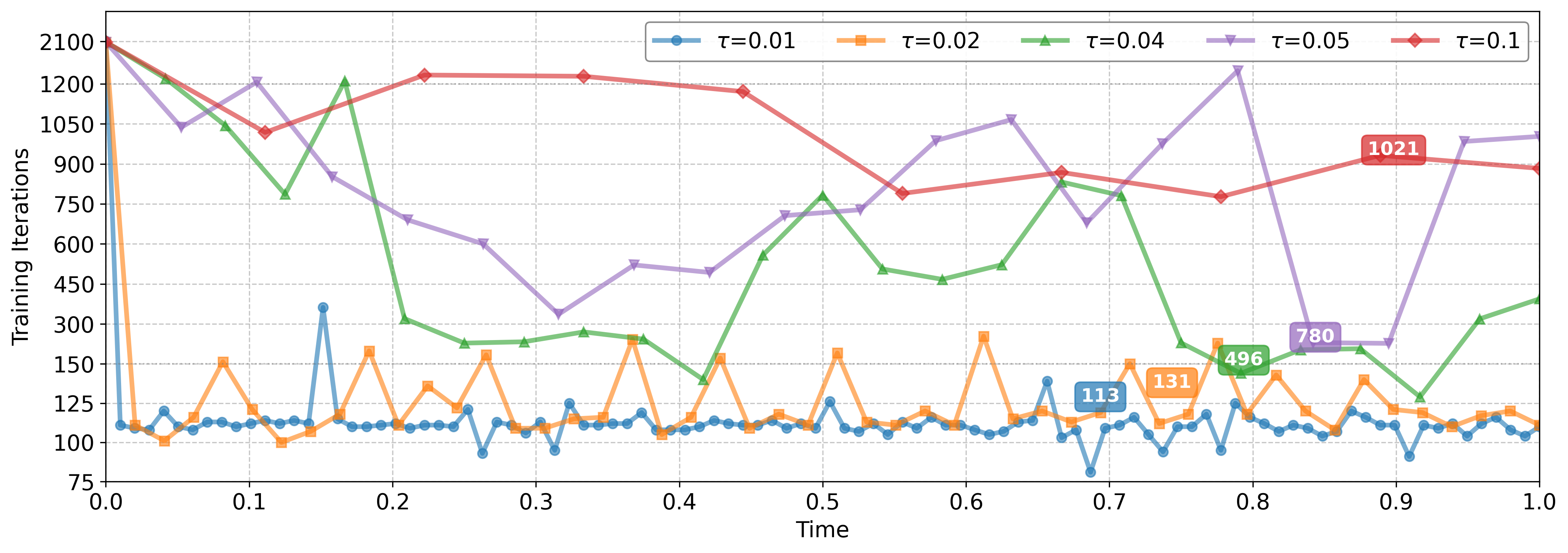}
		\caption{\Cref{example1}, training epochs for different time steps.}
		\label{Ex2D1-Iter}
	\end{figure}

	\Cref{Ex2D1-Iter} illustrates the evolution of the neural network training epochs over time for different choices of the time-step size $\tau$, where a vertical scale is adopted to improve the visibility of variations in training epochs. 
	For every tested value of $\tau$, the first time level is the most computationally demanding, requiring more than $2000$ training epochs. Once this initial network has been trained, the training effort is substantially reduced at subsequent time levels. When $\tau$ is small, the epoch counts quickly stabilize and remain nearly uniform over time. As $\tau$ increases, however, the epoch counts exhibit larger fluctuations, suggesting that larger temporal increments produce greater variations in the complexity of the learning task between successive time levels.
	The total number of training epochs (excluding the initial step) and the average number of training epochs per time step (excluding the initial step) are reported in \Cref{table:avg-iter}. 
	From these results, we observe that smaller $\tau$ reduces the discrepancy between consecutive time levels, thereby decreasing the number of training epochs required at each step. However, since a smaller $\tau$ results in more time levels overall, the total computational cost in training may increase. In summary, as $\tau$ increases, the total number of training epochs initially decreases and then increases. 
	Moreover, as the time step size $\tau$ becomes larger, the advantage of reusing the neural network parameters from the previous step gradually diminishes. Nevertheless, parameter reuse remains beneficial, as it still reduces the number of training epochs compared with retraining from scratch at every time level.
	
\begin{figure}[!ht]
		\centering
		\includegraphics[width=\linewidth]{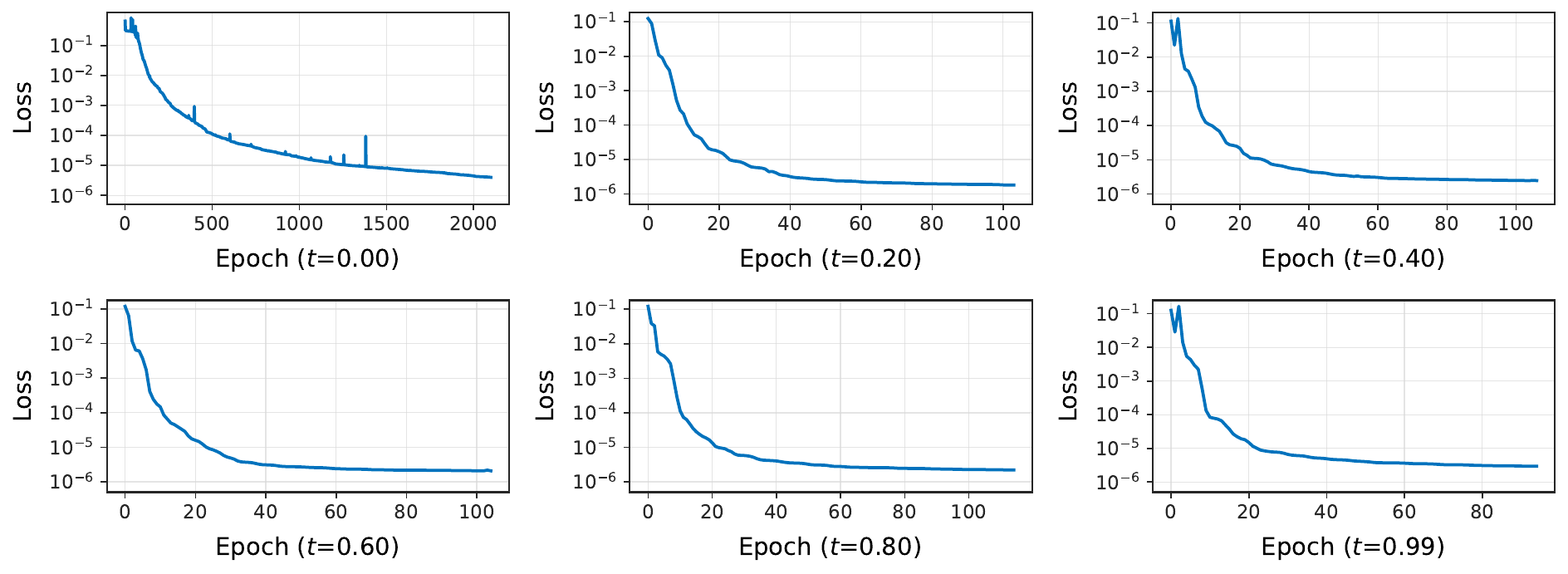}
		\caption{\Cref{example1}, training loss convergence curves at different time levels.}
		\label{Ex2D1-loss}
	\end{figure}

	
	To further examine the optimization behavior under the adopted setting, \Cref{Ex2D1-loss} reports the training loss for $\tau=0.01$ at different time levels. 
	The loss curves show that the MSE at each time level decreases to a value between $10^{-6}$ and $10^{-5}$. 
	At the first time level, approximately $2000$ training epochs are required to reach this accuracy, whereas fewer than $100$ epochs are sufficient at all subsequent time levels.

\begin{figure}[!ht]
		\includegraphics[width=0.9\linewidth]{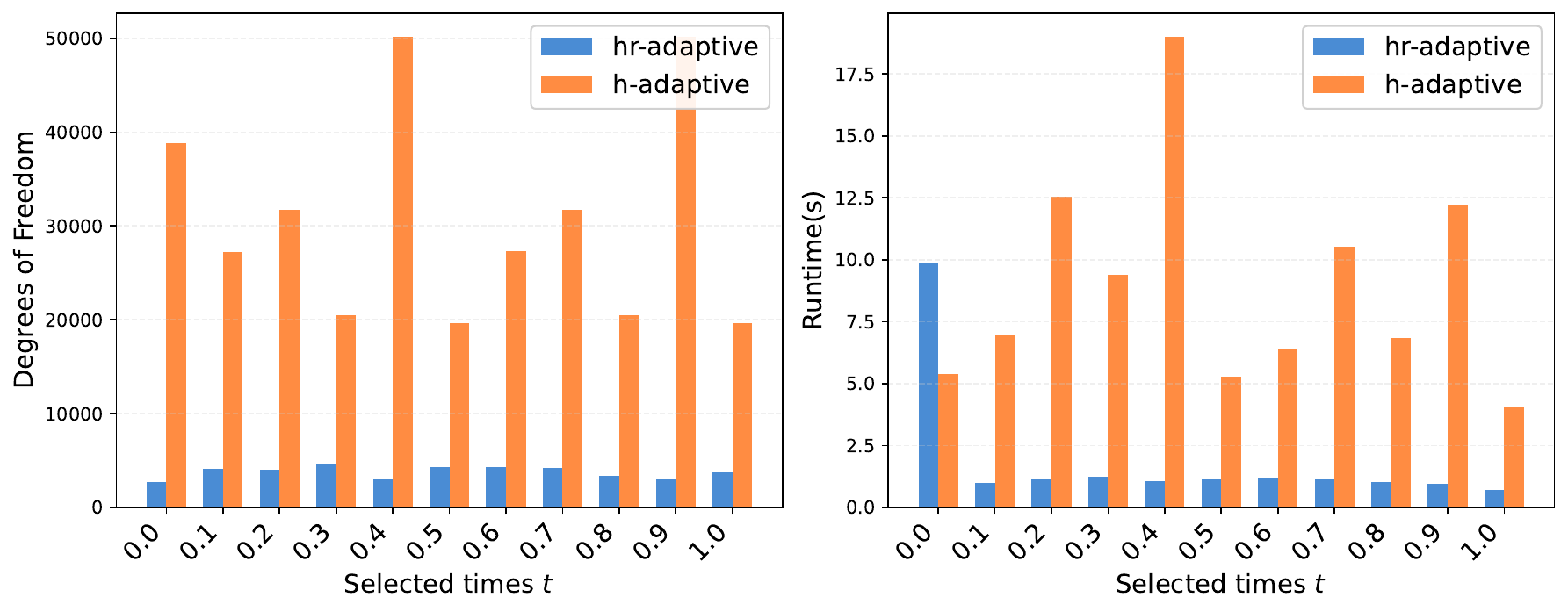}
		\caption{\Cref{example1}, comparison of $hr$-adaptive and $h$-adaptive methods in terms of degrees of freedom (left) and runtime (right).}
		\label{COMPAR}
	\end{figure}

	Next, we revisit the experiment reported in \Cref{table:avg-iter} with the time-step size fixed at $\tau=0.01$ and present the corresponding average training epochs under different network architectures and initialization strategies.  Specifically, \Cref{tab:convergence_epochs_avg_0.01} summarizes the average number of training epochs (excluding the initial time step), computed over three independent runs, required for the L-BFGS optimizer to reach its default stopping criteria. 
	The results indicate that both the network architecture and the initialization strategy have some influence on the training cost. Nevertheless, all tested configurations converge within a reasonable number of epochs, demonstrating the robustness of the proposed approach. 
	Among the network architectures considered, deeper networks generally require fewer training epochs than shallower ones, particularly for the wider architectures. This can be attributed to the greater representational capacity of deeper networks, together with the warm-start initialization from the previous time step, which facilitates rapid convergence through minor parameter updates. In contrast, shallower networks may require larger parameter updates to achieve comparable accuracy, resulting in higher training cost.
	Regarding initialization strategies, no single method consistently achieves the lowest training cost across all architectures. Nevertheless, Uniform, Xavier, and Kaiming initializations often lead to slightly faster convergence, whereas Normal initialization is less frequently the best-performing choice. Overall, the differences among initialization methods remain moderate, especially for deeper networks, suggesting that the convergence behavior of the proposed method is relatively insensitive to the choice of initialization.
	Additionally, \Cref{COMPAR} compares our neural network–enhanced $hr$-adaptive Algorithm~\ref{alg:ParabolicHR} with the conventional $h$-adaptive method in \cite{Picasso1998}. 
	The results show that our algorithm achieves a significant reduction in the degrees of freedom throughout the computation. 
	Regarding computational cost, although the initial step is more time-consuming due to the extensive neural network training required, the subsequent time steps are considerably faster. 
	As a result, the proposed method attains a substantial overall time advantage, as quantified in \Cref{tab:runtime}.

	
	\begin{table}[!ht]
		\centering
		\caption{\Cref{example1}, training epochs (excluding the initial step) for different time step $\tau$.}
		\label{table:avg-iter}
		\begin{tabular}{|c|cccccc|} 
			\hline  
			Time step size $\tau$ & $0.001$ & $0.01$ & $0.02$ & $0.04$ & $0.05$ & $0.1$ \\
			\hline
			Total training epochs & $47626$ & $11186$ & $6409$ & $11911$ & $14829$ & $9186$  \\
			\hline
			Average training epochs & $47.7$ & $113.0$ & $130.8$ & $496.3$ & $780.5$ & $1020.7$ \\
			\hline  
		\end{tabular}
	\end{table}
	
	\begin{table}[htbp]
		\centering
		\caption{Average training epochs (excluding the initial step) required for convergence for different network architectures and initialization methods with time step size $\Delta t = 0.01$. The smallest value in each row is highlighted in bold.}
		\label{tab:convergence_epochs_avg_0.01}
		\small
		\renewcommand{\arraystretch}{1.15}
		\setlength{\tabcolsep}{10pt}
		
		\begin{tabular}{llcccc}
			\toprule
			Depth & Architecture & Normal & Uniform & Xavier & Kaiming \\
			\midrule
			\multirow{2}{*}{3}
			& \texttt{[2-16-16-1]}
			& 131.0 & \textbf{55.5} & 75.4 & 150.6 \\
			& \texttt{[2-40-40-1]}
			& 165.0 & 165.2 & 302.0 & \textbf{151.5} \\
			\midrule
			
			\multirow{2}{*}{4}
			& \texttt{[2-16-16-16-1]}
			& 85.5 & 93.3 & 86.6 & \textbf{77.5} \\
			& \texttt{[2-40-40-40-1]}
			& 119.5 & 108.6 & \textbf{93.1} & 102.8 \\
			\midrule
			
			\multirow{2}{*}{5}
			& \texttt{[2-16-16-16-16-1]}
			& 74.8 & 83.2 & \textbf{70.7} & 80.4 \\
			& \texttt{[2-40-40-40-40-1]}
			& 109.2 & \textbf{96.0} & 98.5 & 101.6 \\
			\bottomrule
		\end{tabular}
	\end{table}
	
	\begin{table}[!ht]
		\centering
		\caption{Total CPU running time comparison of algorithms on three examples.}
		\label{tab:runtime}
		\begin{tabular}{lccc}
			\toprule
			Algorithm & \Cref{example1} & \Cref{example2} & \Cref{example3} \\
			\midrule
			Conventional $h$-adaptive & $818.75s$  & $959.76s$ & $555.20s$ \\
			NN‑enhanced $hr$-adaptive & $132.98s$ & $308.45s$ & $143.96s$ \\
			\bottomrule
		\end{tabular}
	\end{table}

\begin{figure}
		\centering
		\includegraphics[width=0.21\linewidth]{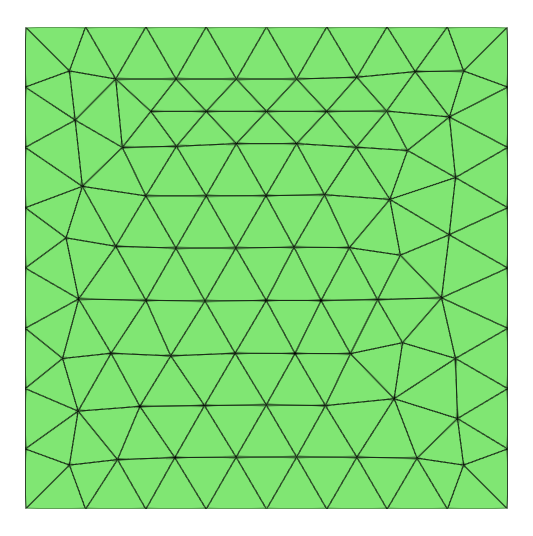}\hspace{0.3cm}
		\includegraphics[width=0.21\linewidth]{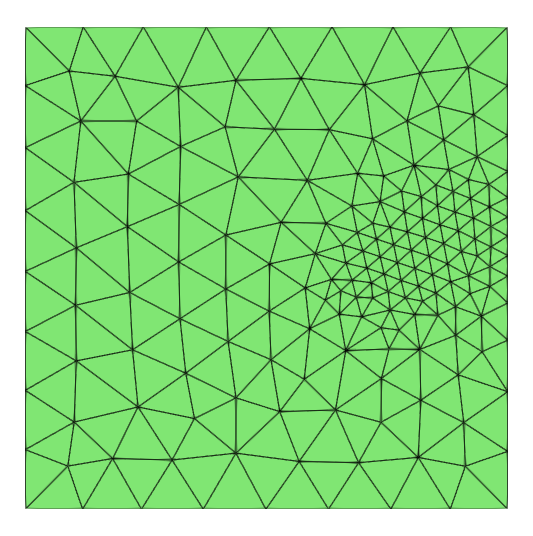}\hspace{0.3cm}
		\includegraphics[width=0.21\linewidth]{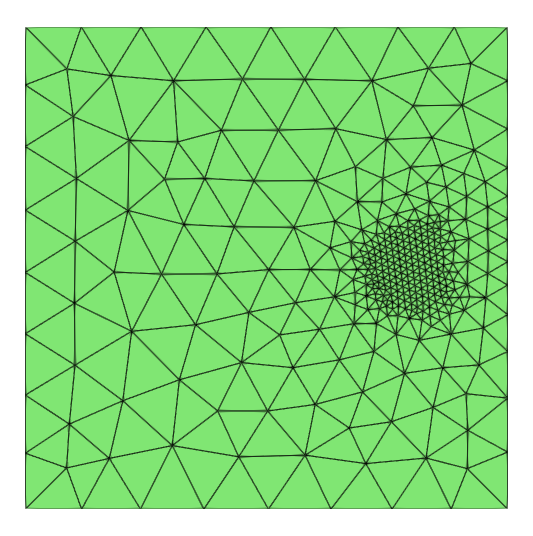} \\
		\includegraphics[width=0.21\linewidth]{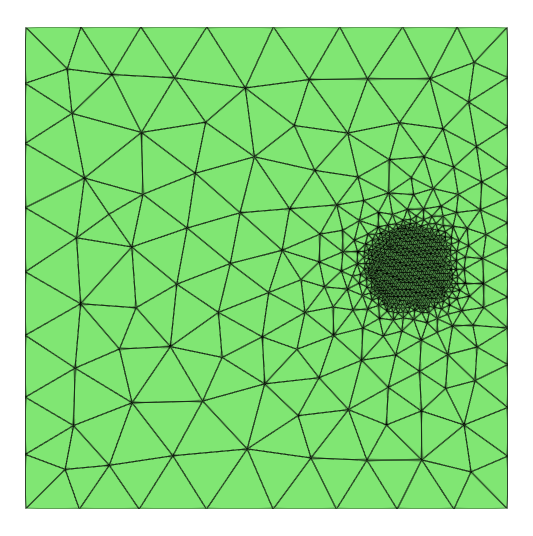}\hspace{0.3cm}
		\includegraphics[width=0.21\linewidth]{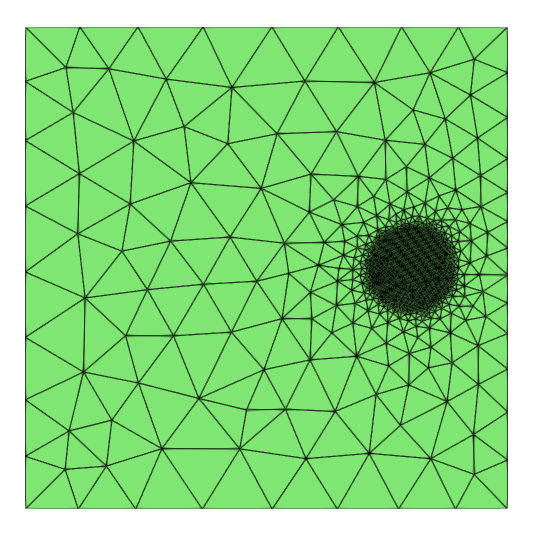}\hspace{0.3cm}
		\includegraphics[width=0.21\linewidth]{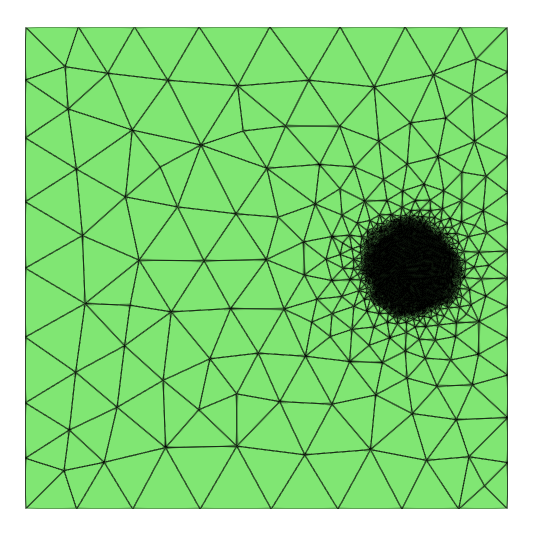}
		\caption{\Cref{example1}, initial mesh (top left) and its evolution through five adaptive refinements at $t=0.0$.}
		\label{Ex2D1-Time}
	\end{figure}

	In \Cref{Ex2D1-Time}, the top-left picture shows the initial uniform mesh $\mathcal{T}_h^{0,0}$ ($98$ nodes), which also serves as the initial mesh for adaptive iterations at each time step. 
	The next five pictures show the meshes over five successive refinements, with the NOV being $180$, $340$, $681$, $1388$, and $12818$, respectively.  
	It can be observed that the NOV approximately doubles at each step, and the mesh generated by \Cref{alg:ParabolicHR} at the fifth step is significantly larger than that of the previous mesh. 
	The process terminates before the sixth refinement as the global error estimator is smaller than the tolerance $eTol=0.01$, demonstrating that our fitted model \eqref{fitModel} indeed plays its intended role. 
	\begin{figure}
		\includegraphics[width=0.45\linewidth]{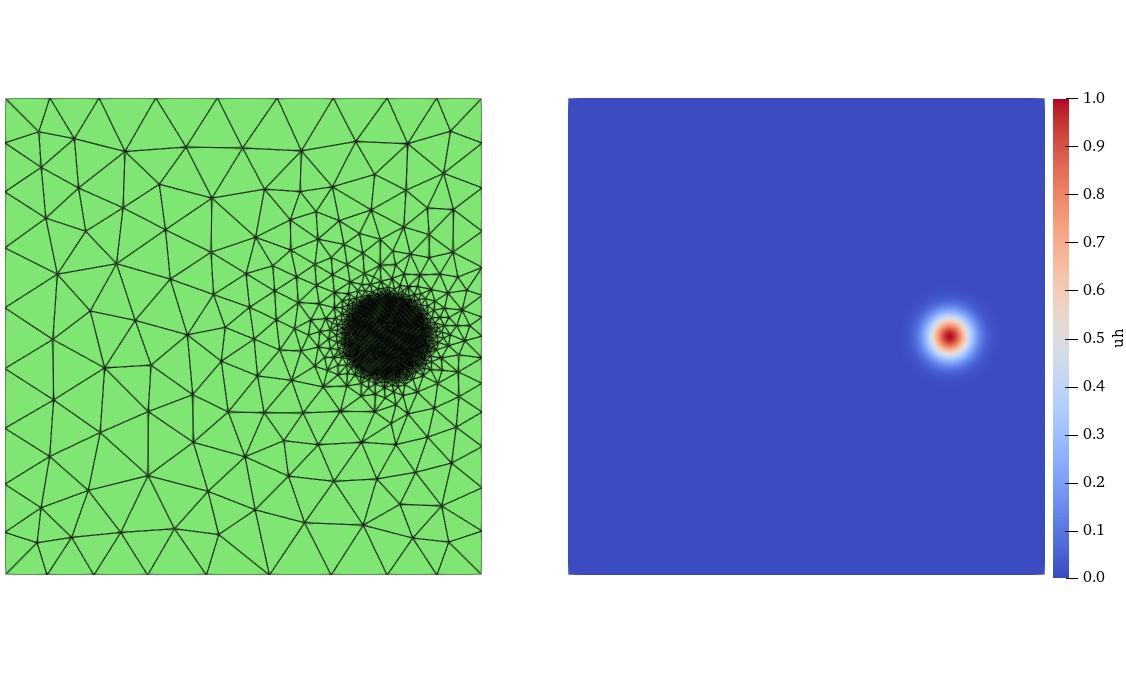}
		\includegraphics[width=0.45\linewidth]{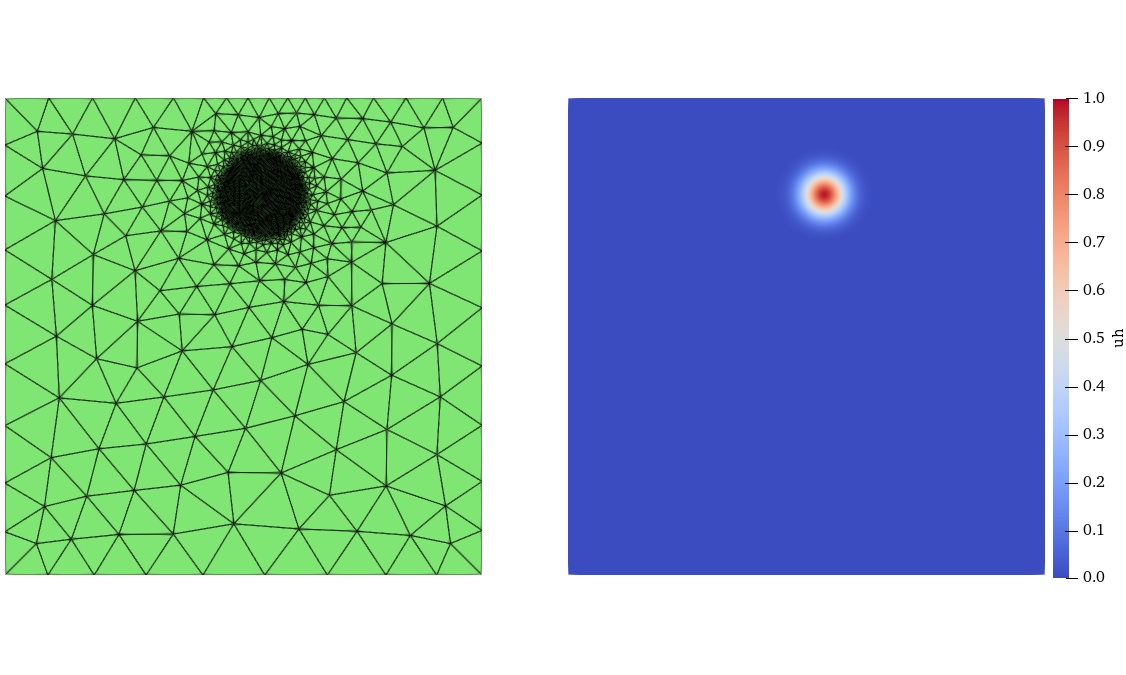}\vspace{-0.7cm}
		\includegraphics[width=0.45\linewidth]{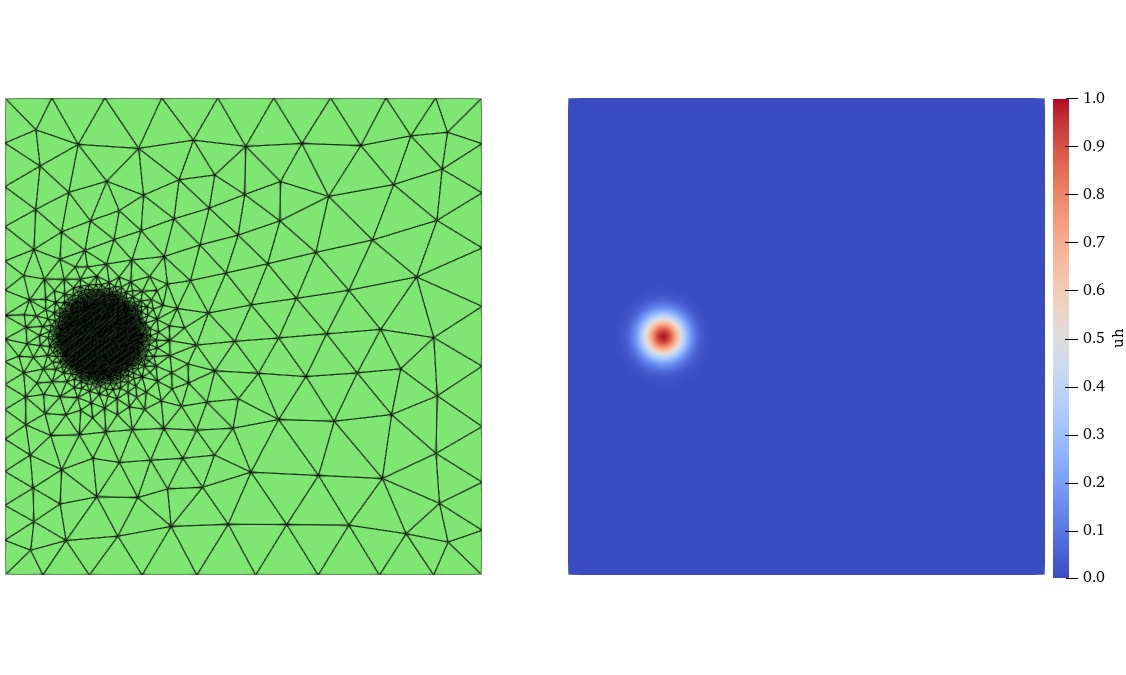}
		\includegraphics[width=0.45\linewidth]{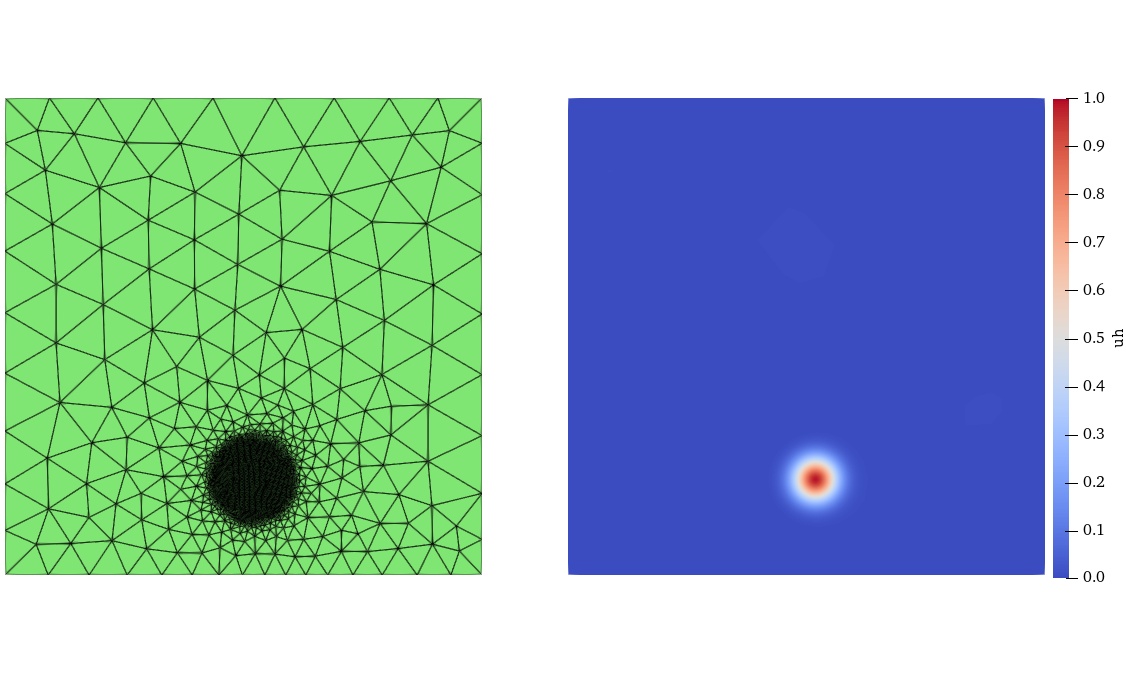}\vspace{-0.8cm}
		\caption{\Cref{example1}, snapshots of numerical solutions along with their corresponding adaptive meshes at $t=0.00, 0.25, 0.50$, and $t=0.75$.}
		\label{Ex2D1-Initial}
	\end{figure}
	
	\Cref{Ex2D1-Initial} illustrates the evolution of the numerical solutions and the corresponding mesh at different time steps.
	The mesh dynamically refines around sharp gradients or singularities, efficiently allocating resources where needed, and reveals a peak rotating around the coordinate origin.
	Notably, the grids shown here differ significantly from those in \Cref{Exremark} due to the smaller time step used, which causes the positions of two distinct singularities at adjacent time steps to nearly overlap. 
	The quasi-optimal convergence rate of the adaptive algorithm for the gradient error is demonstrated in the left panel of \Cref{Ex2-order}.

	\begin{example} \label{example2}
		(Diffusion) Consider the 2D parabolic equation \eqref{model} in the domain $\Omega=[-1,1]^2$, with the source function $f$ such that the exact solution $u$ is given by:
		\[u(x,y,t) = \exp(-5000(\sqrt{x^2+y^2}+0.3t-0.4)^2).\]
		\begin{figure}[!ht]
			\centering
			\includegraphics[width=0.95\linewidth]{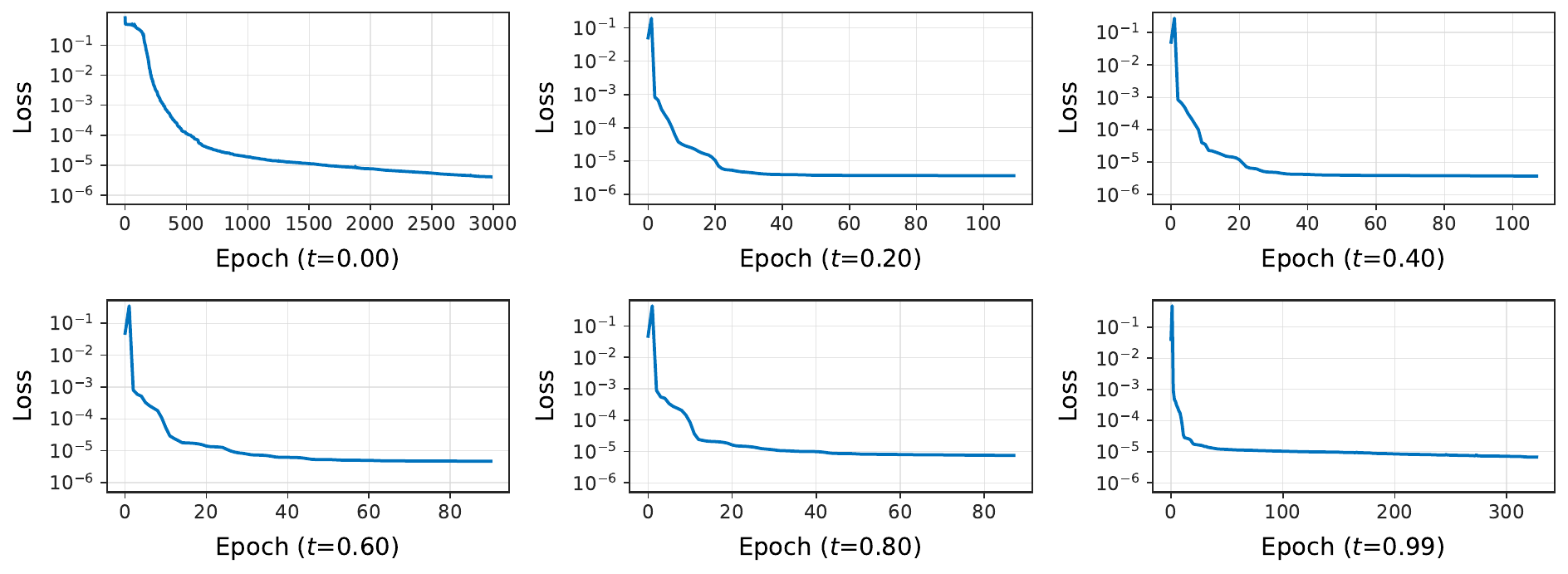}
			\caption{\Cref{example2}, training loss convergence curves at different time levels.}
			\label{Ex2D2-loss}
		\end{figure}
		\begin{figure}
			\centering
			\begin{minipage}[c]{0.24\linewidth}
				\raggedleft
				\includegraphics[width=0.73\linewidth]{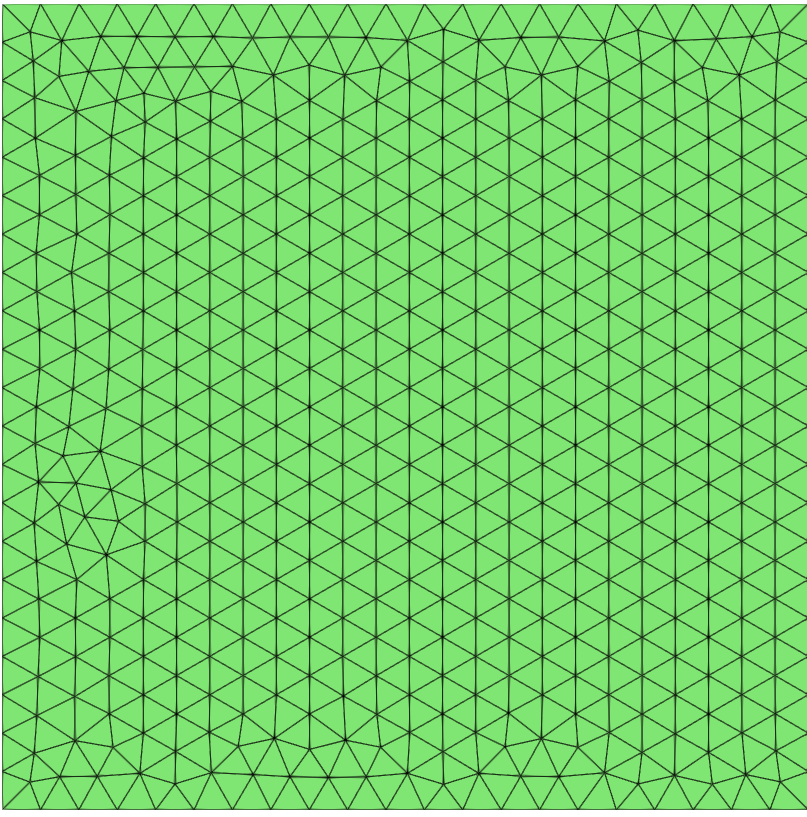}
			\end{minipage}%
			\hfill
			\begin{minipage}[c]{0.7\linewidth}
				\raggedright
				\includegraphics[width=0.28\linewidth]{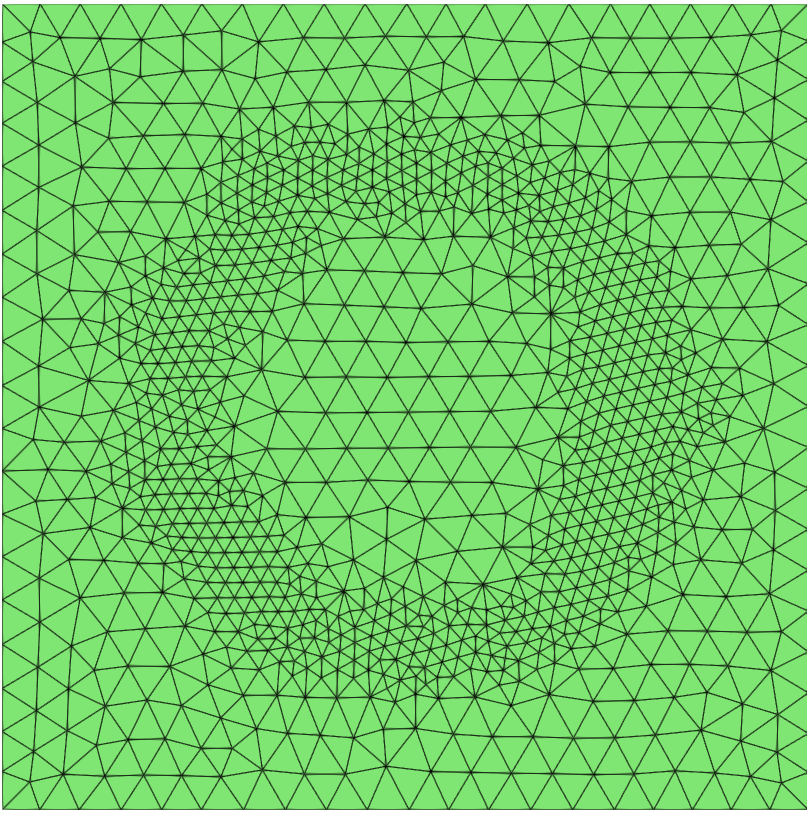}\hspace{0.2cm}
				\includegraphics[width=0.28\linewidth]{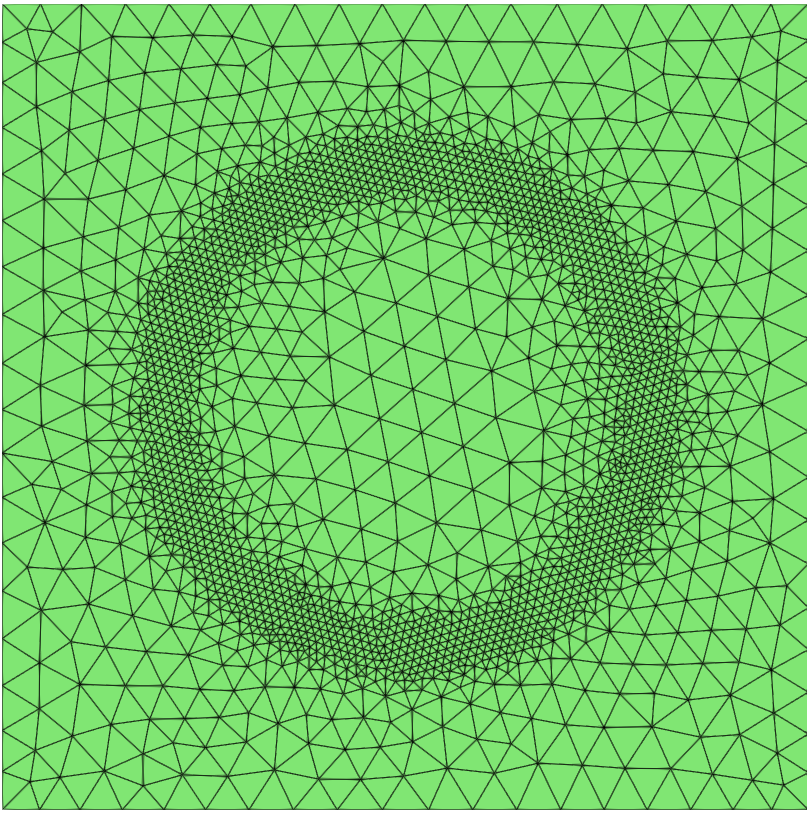}\hspace{0.2cm}
				\includegraphics[width=0.28\linewidth]{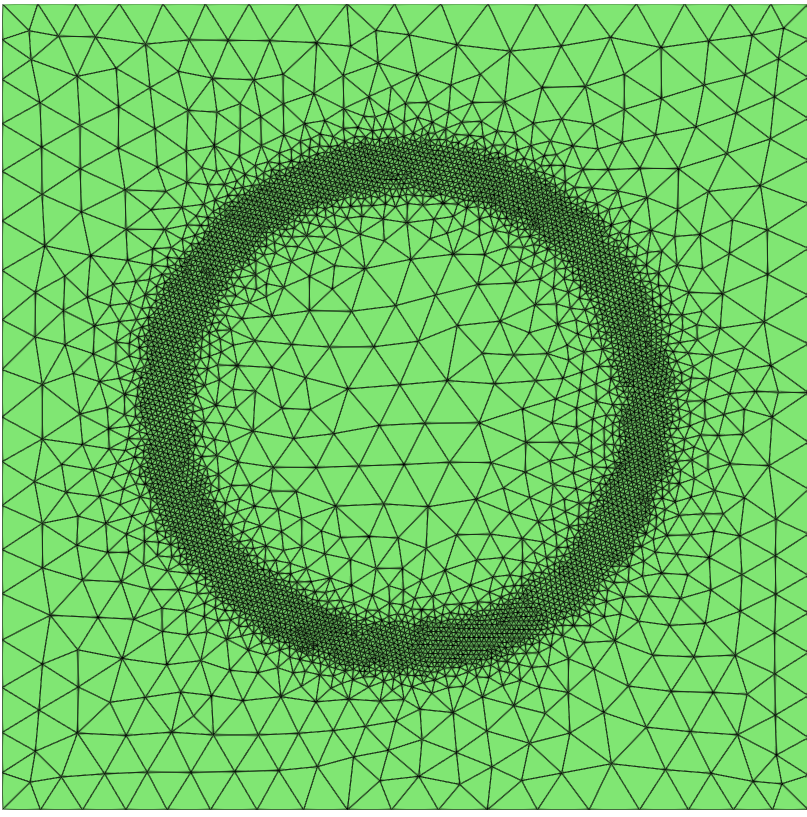}\vspace{0.1cm}\\
				\includegraphics[width=0.28\linewidth]{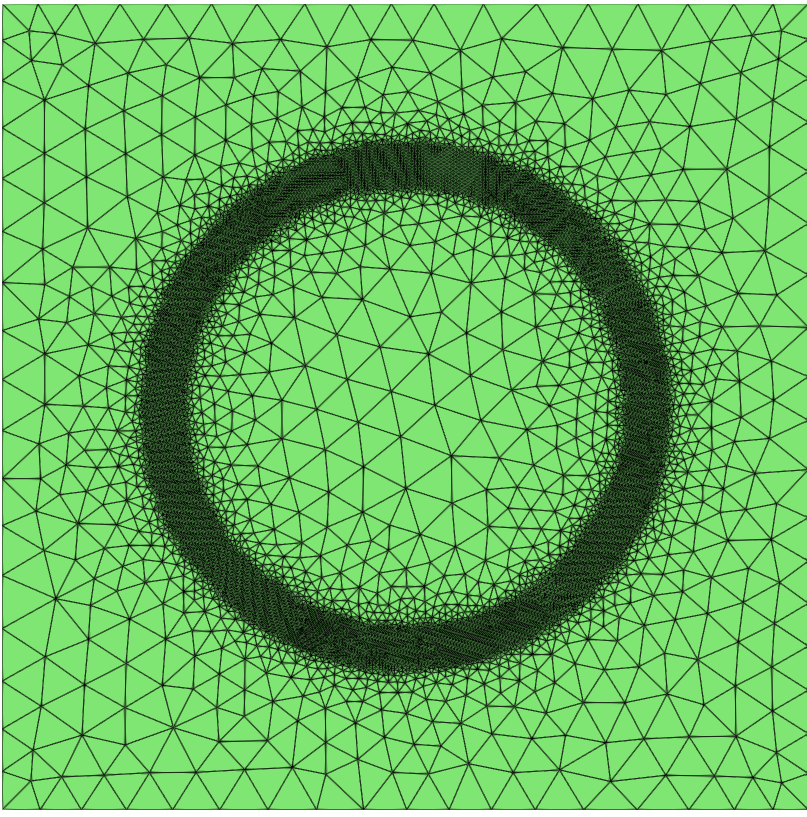}\hspace{0.2cm}
				\includegraphics[width=0.28\linewidth]{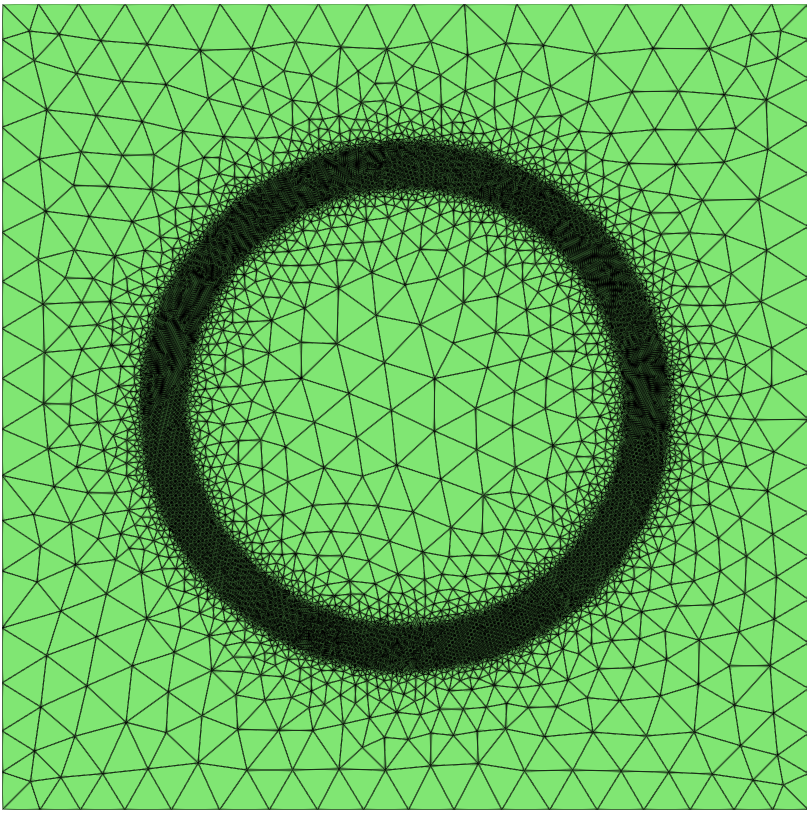}\hspace{0.2cm}
				\includegraphics[width=0.28\linewidth]{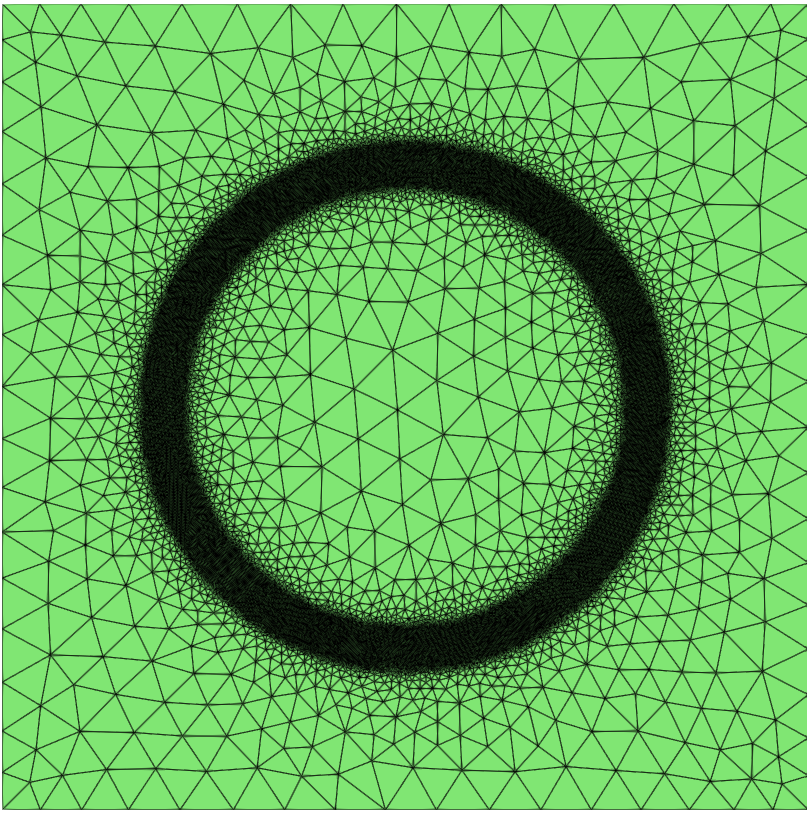}
			\end{minipage}
			\caption{\Cref{example2}, initial mesh (left) and its evolution through six refinements at $t=0.0$.}
			\label{Ex2D2-Time}
		\end{figure} 
		\begin{figure}
			\centering
			\includegraphics[width=0.85\linewidth]{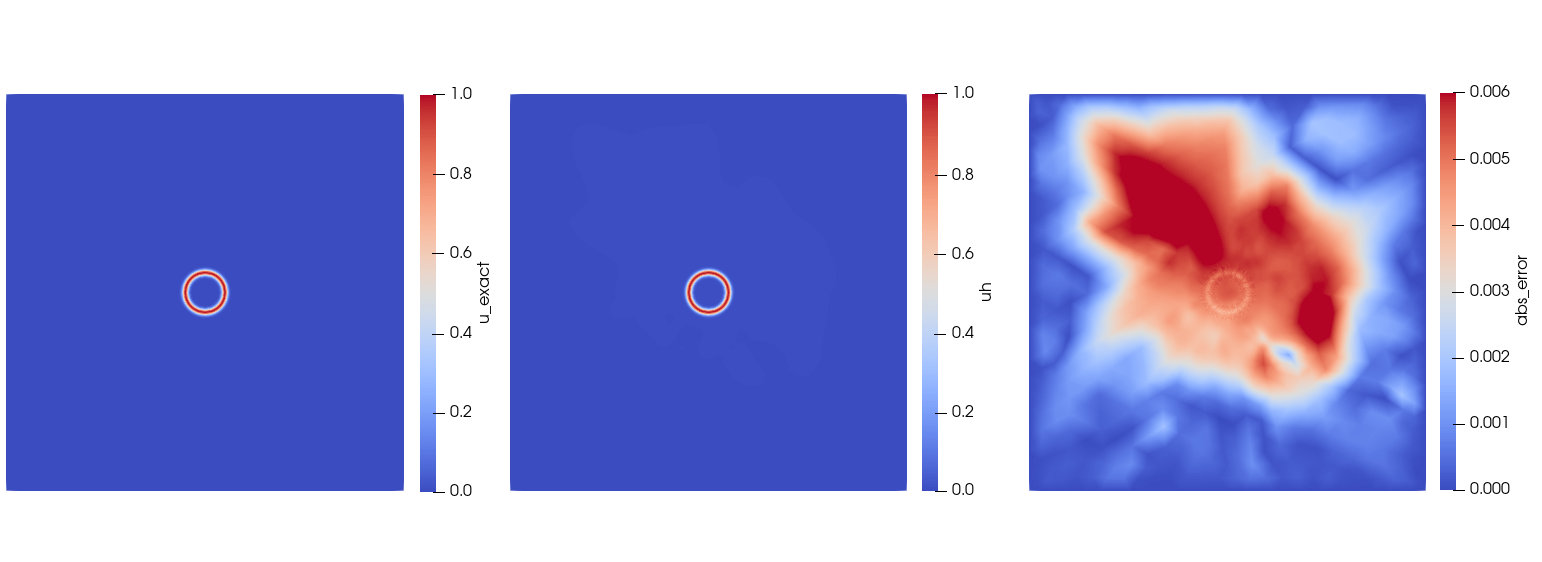}
			\caption{\Cref{example2}, comparison with the exact solution at t=1.0. (Left) Exact solution, (middle) finite element solution, and (right) pointwise error.}
			\label{Ex2D2-abs}
		\end{figure}
		\begin{figure}
			\centering
			\includegraphics[width=0.45\linewidth]{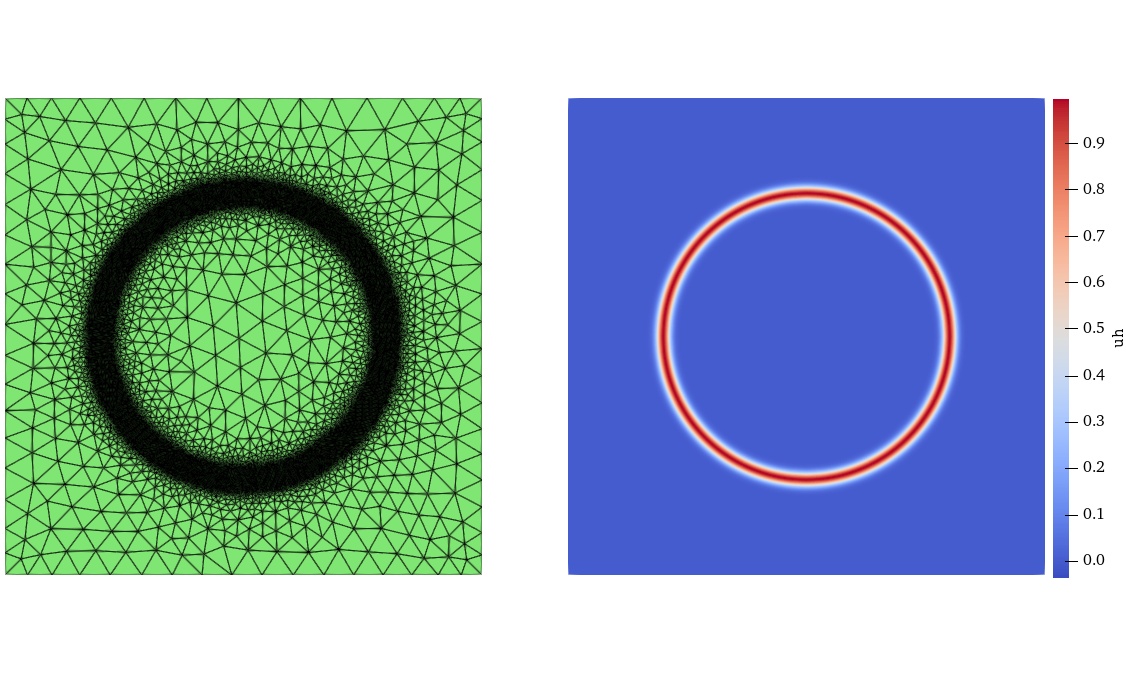}
			\includegraphics[width=0.45\linewidth]{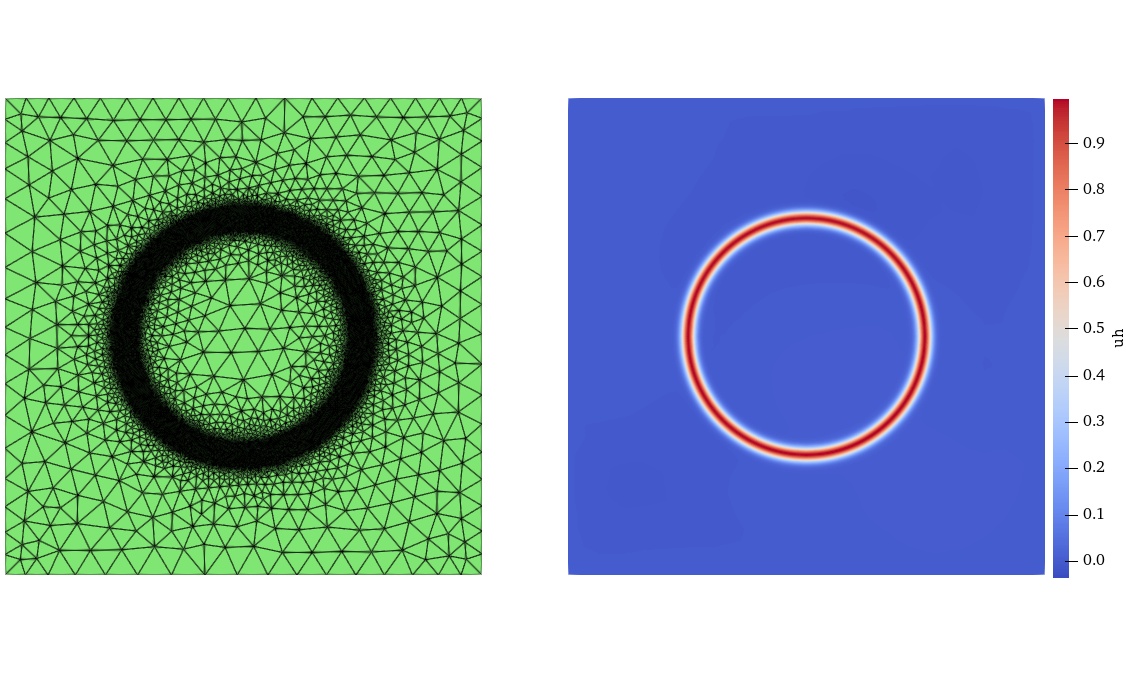}\vspace{-0.7cm}
			\includegraphics[width=0.45\linewidth]{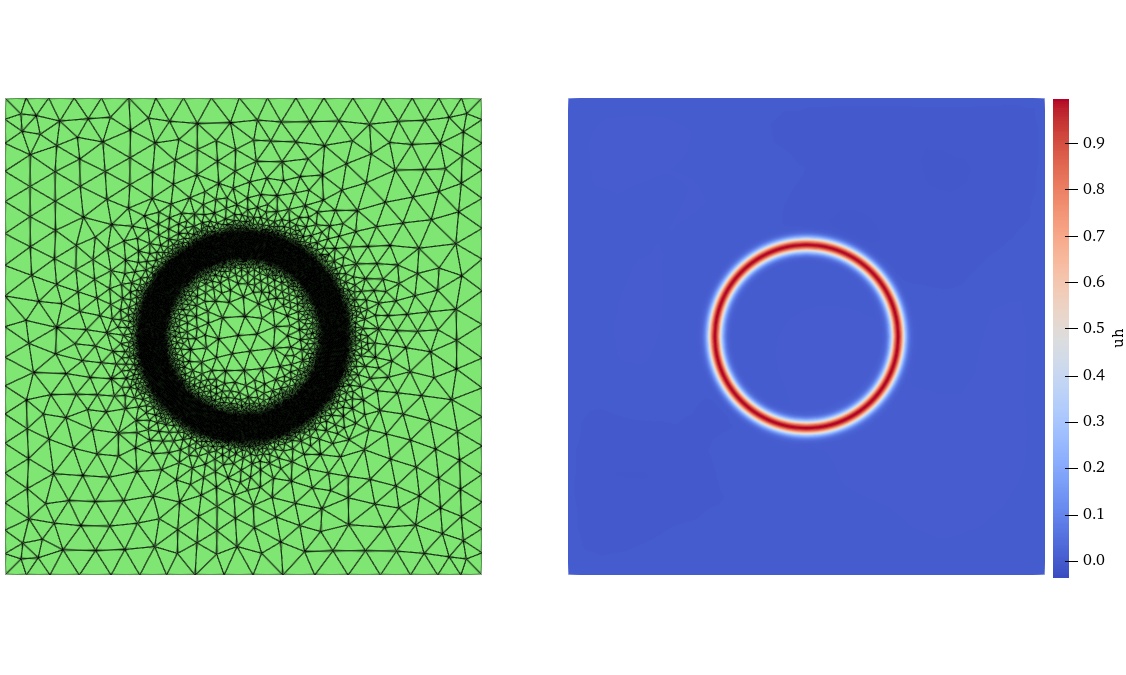}
			\includegraphics[width=0.45\linewidth]{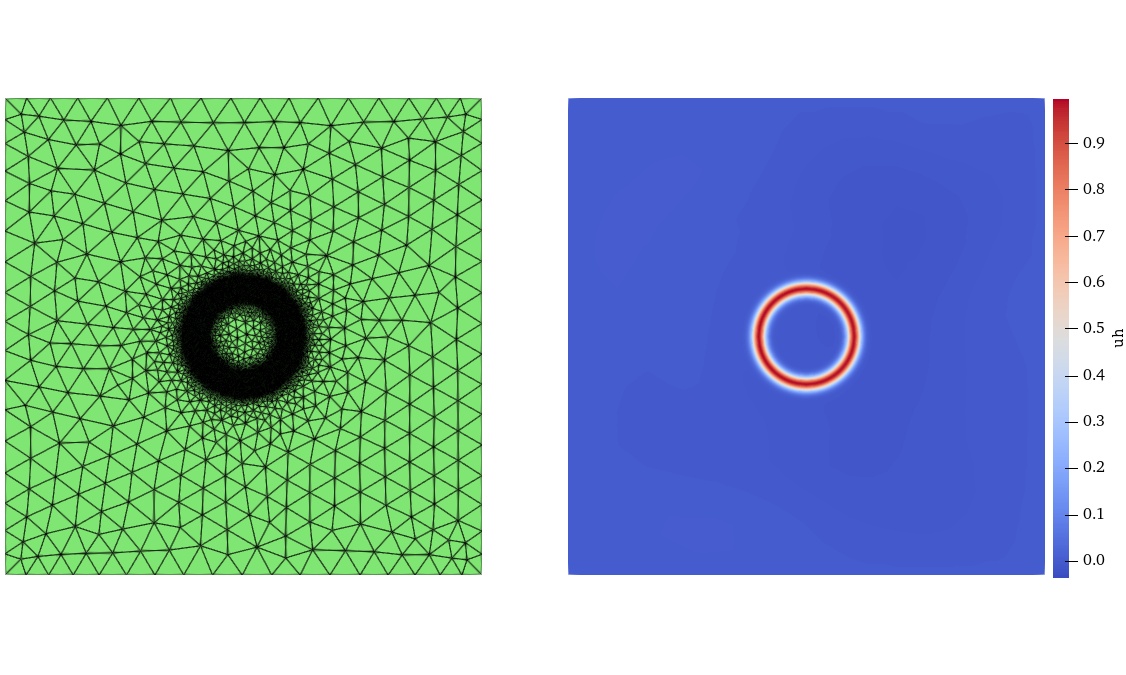}\vspace{-0.8cm}
			\caption{\Cref{example2}, snapshots of numerical solutions along with their corresponding adaptive meshes at $t=0.00, 0.25, 0.50$, and $t=0.75$.}
			\label{Ex2D2-Initial}
		\end{figure}
	\end{example}
	In this example, the settings are the same as \Cref{example1}, except for the tolerance $eTol=0.05$. 
	Training loss histories at different time levels are presented in \Cref{Ex2D2-loss}. Similar to \Cref{example1}, approximately 2000 training epochs are required at the initial time level to achieve the MSE below $10^{-5}$. In contrast, fewer than 100 epochs are sufficient at all subsequent time levels to attain comparable accuracy. 
	\Cref{Ex2D2-Time} illustrates the adaptive refinement process of the initial mesh over six iterations. The NOV at each iteration are $556$, $1088$, $2182$, $4436$, $9069$, $18678$, and $38659$, respectively. Unlike the \Cref{example1}, the NOV of the sixth mesh is approximately twice that of the fifth mesh. 
	Moreover, the need for a seventh mesh indicates that the error estimator obtained on the sixth mesh is still larger than the tolerance. This suggests that the mesh generated by the previous refinement did not meet the desired tolerance, even though the tolerance was used when calculating the NOV for the sixth mesh. This implies the necessity of one extra step after the least square fitting, namely the case $k=5$ in \Cref{alg:ParabolicHR}. 
	Four different time steps are shown in \Cref{Ex2D2-Initial}, where the mesh and numerical solutions illustrate a scenario in which a ring-shaped crater gradually decreases its radius over time. 
	The quasi-optimal convergence rate of the adaptive algorithm for the gradient error is demonstrated in the middle panel of \Cref{Ex2-order}. 
	\Cref{Ex2D2-abs} compares the exact and finite element solutions at $t=1$. The numerical solution accurately captures the localized ring-shaped profile. Moreover, the $L^\infty$-error is about $6 \times 10^{-3}$, which implies that the finite element approximation achieves satisfactory accuracy.

	\begin{example} \label{example3}
		(Splitting) Consider the 2D parabolic equation \eqref{model} in the domain $\Omega=[-1,1]^2$, with the source function $f$ such that the exact solution $u$ is given by:
		\[u(x,y,t) = \exp(-300((x-0.3t)^2+y^2)) + \exp(-300((x+0.3t)^2+y^2)).\]
		\begin{figure}[!ht]
			\includegraphics[width=0.95\linewidth]{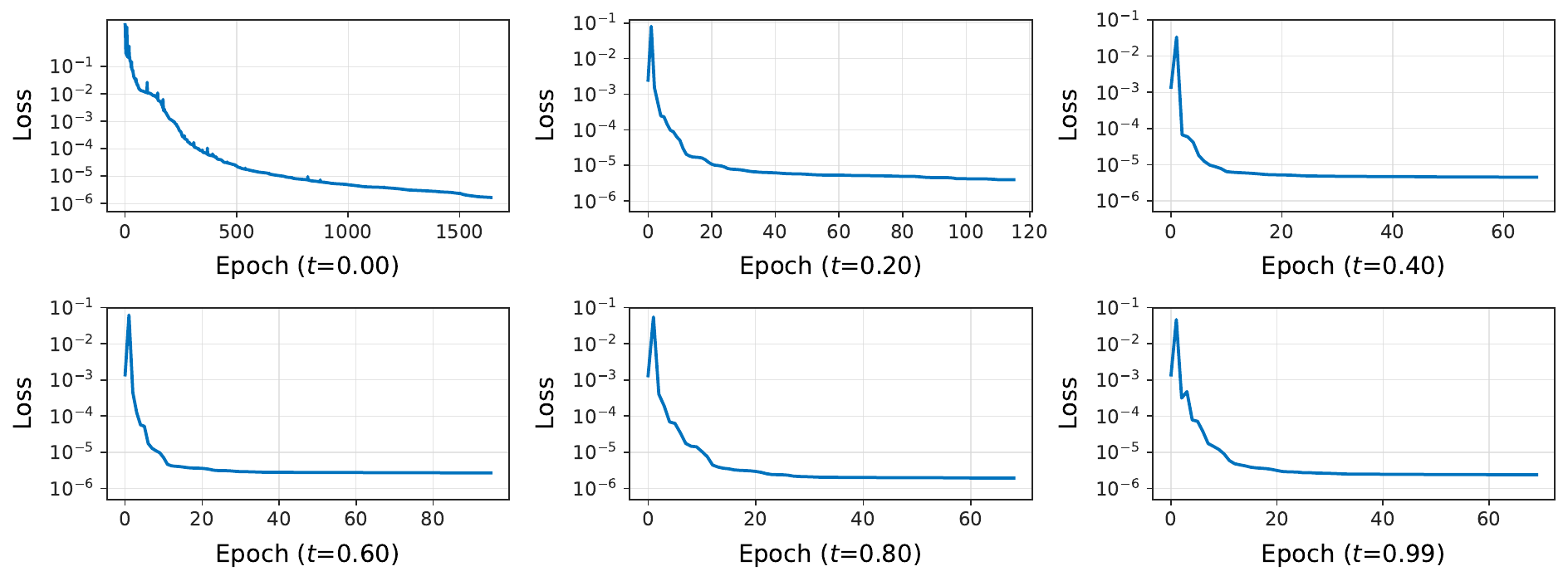}
			\caption{\Cref{example3}, training loss convergence curves at different time levels.}
			\label{Ex2D3-loss}
		\end{figure}
		\begin{figure}
			\centering
			\includegraphics[width=0.19\linewidth]{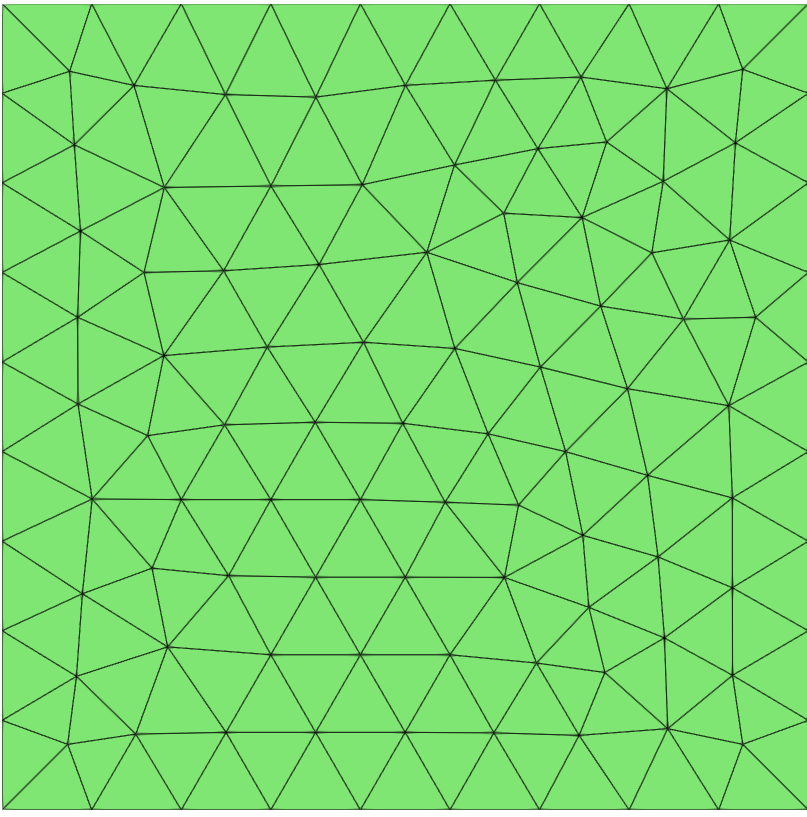}\hspace{0.5cm}
			\includegraphics[width=0.19\linewidth]{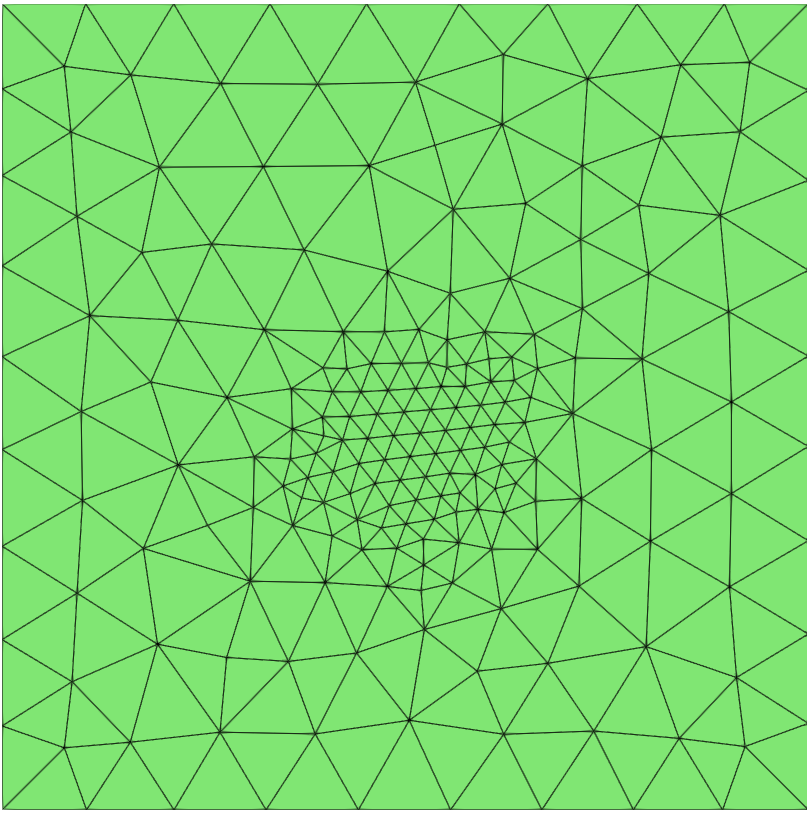}\hspace{0.5cm}
			\includegraphics[width=0.19\linewidth]{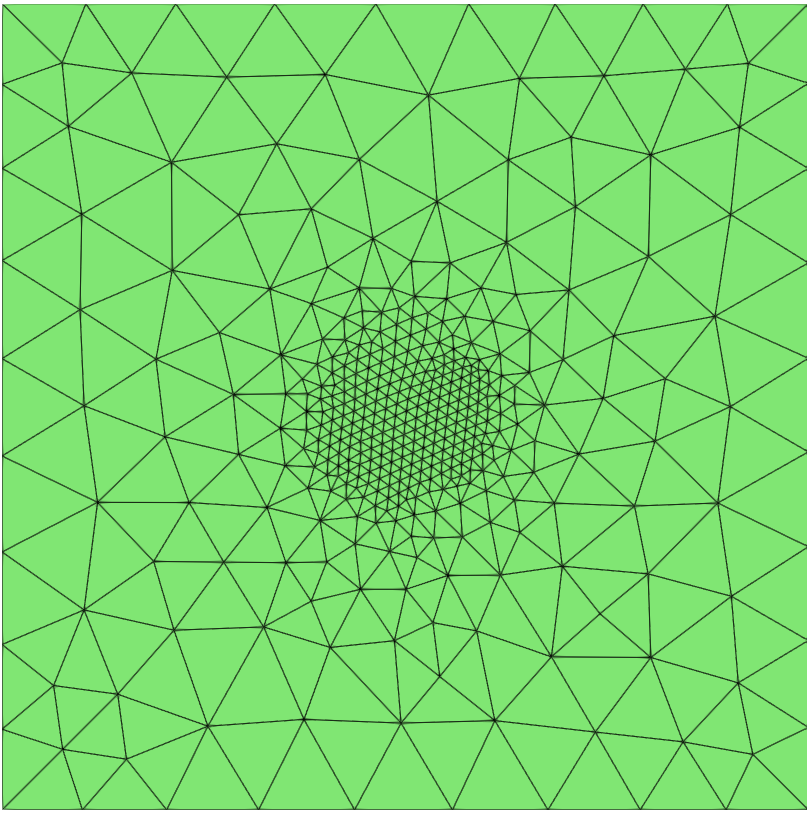}\\
			\includegraphics[width=0.19\linewidth]{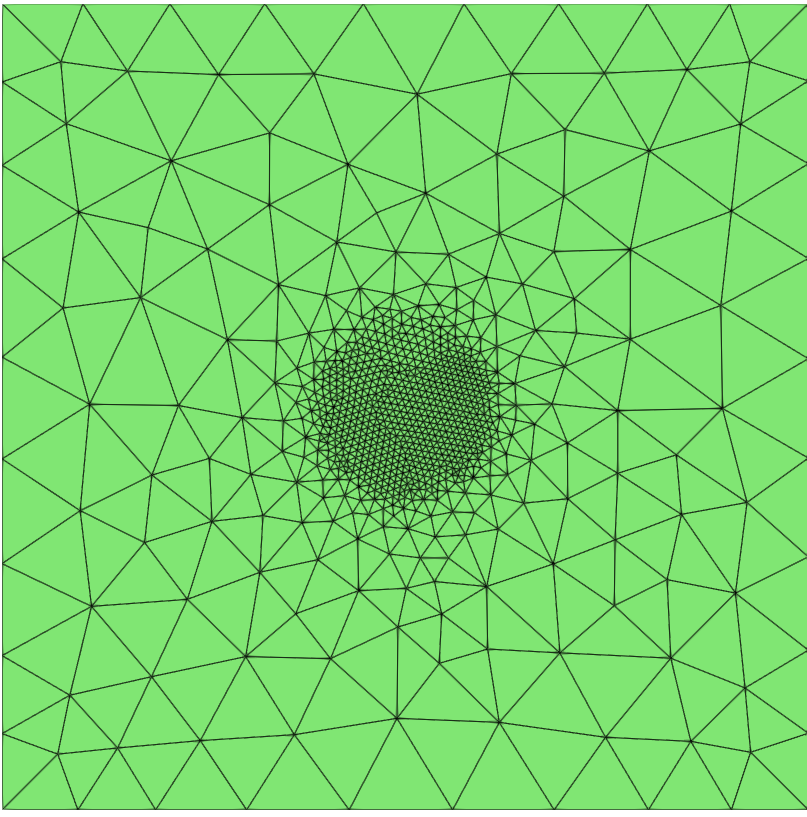}\hspace{0.5cm}
			\includegraphics[width=0.19\linewidth]{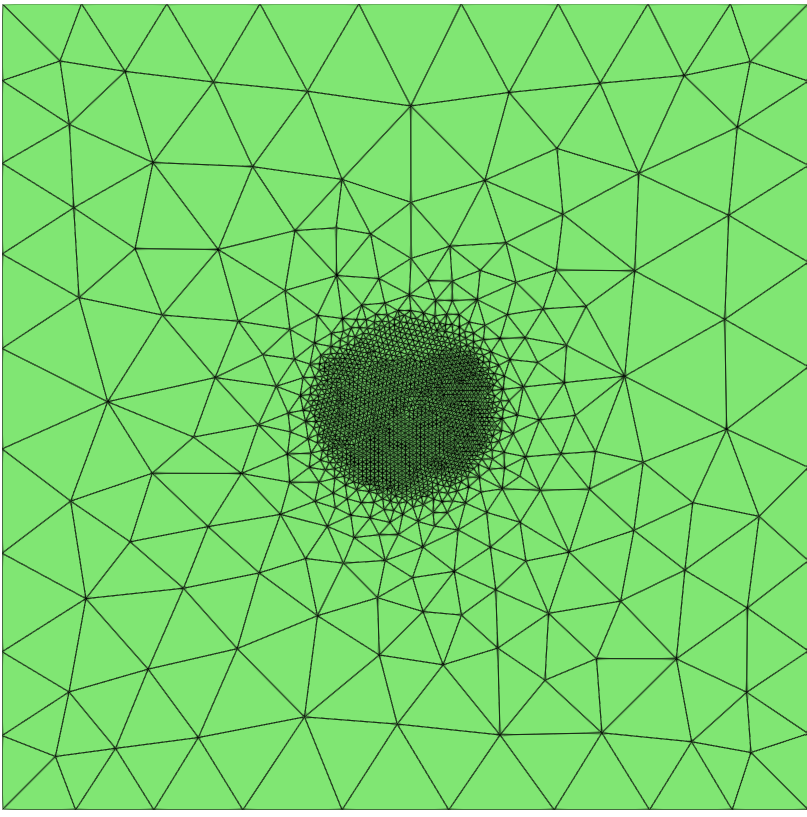}\hspace{0.5cm}
			\includegraphics[width=0.19\linewidth]{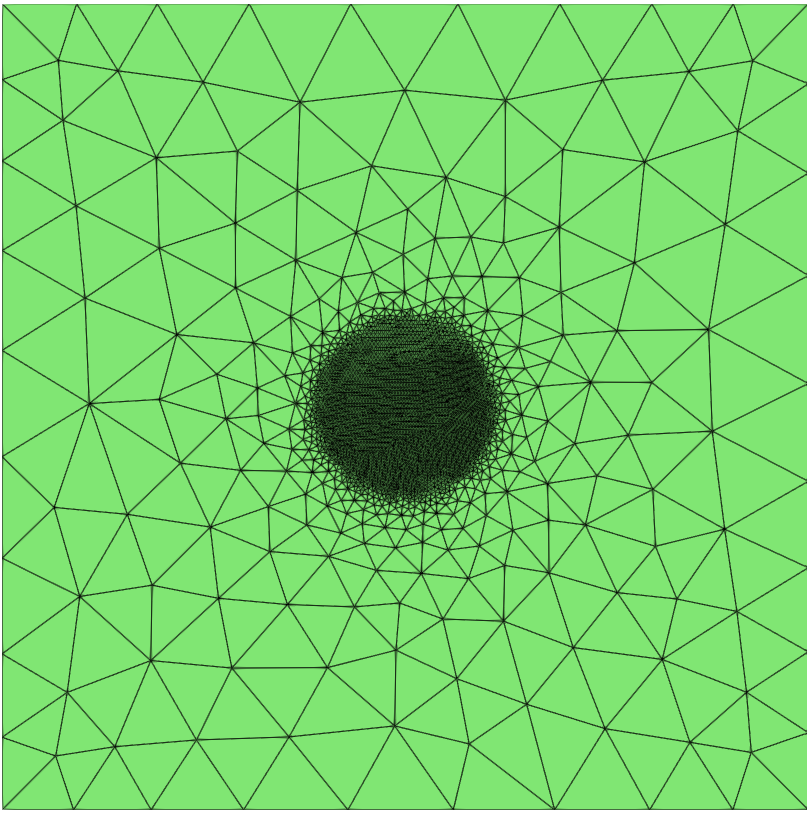}
			\caption{\Cref{example3}, initial mesh (top left) and its evolution through five adaptive refinements at $t=0.0$.}
			\label{Ex2D3-Time}
		\end{figure}
		\begin{figure}
			\centering
			\includegraphics[width=0.45\linewidth]{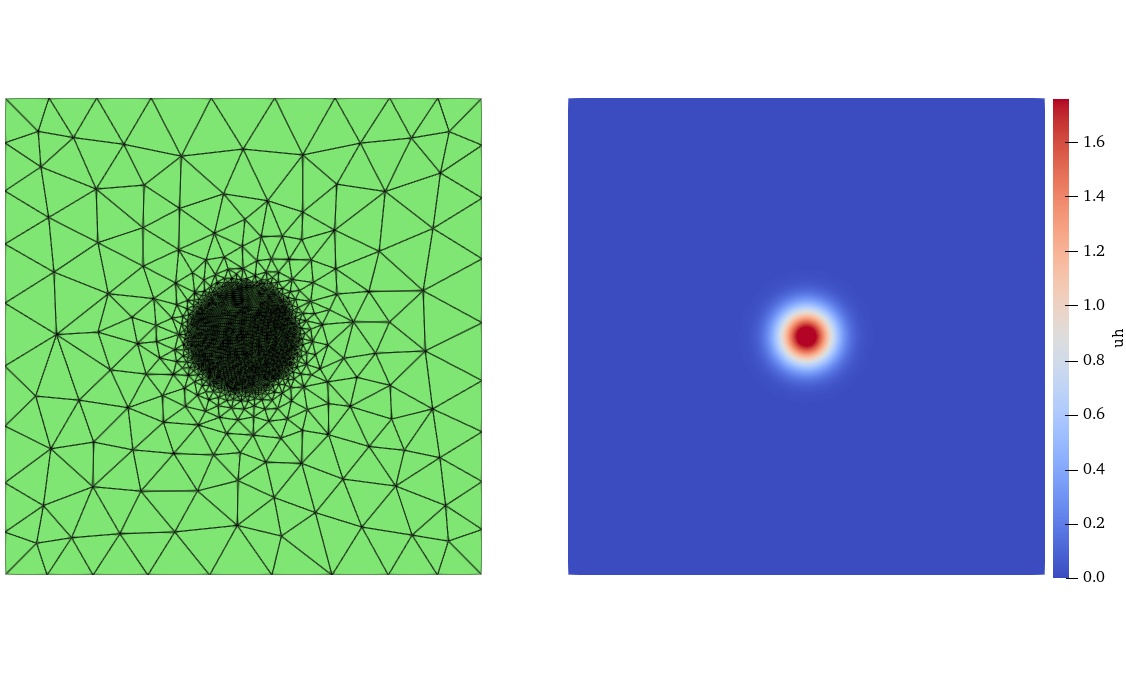}
			\includegraphics[width=0.45\linewidth]{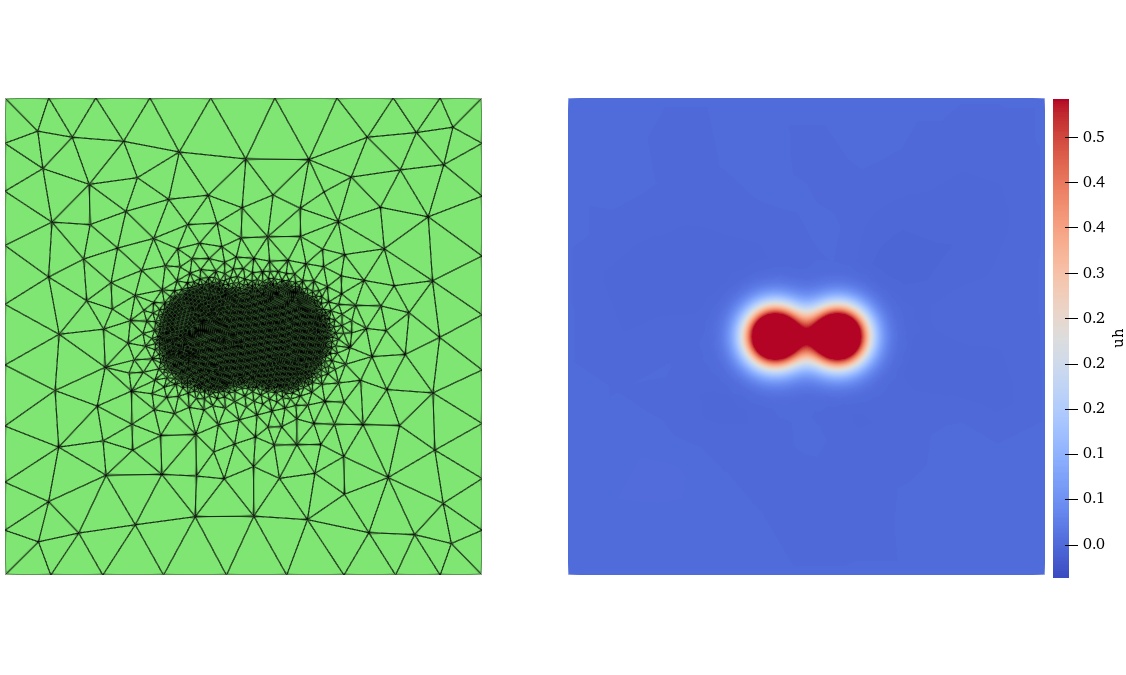}\vspace{-0.7cm}
			\includegraphics[width=0.45\linewidth]{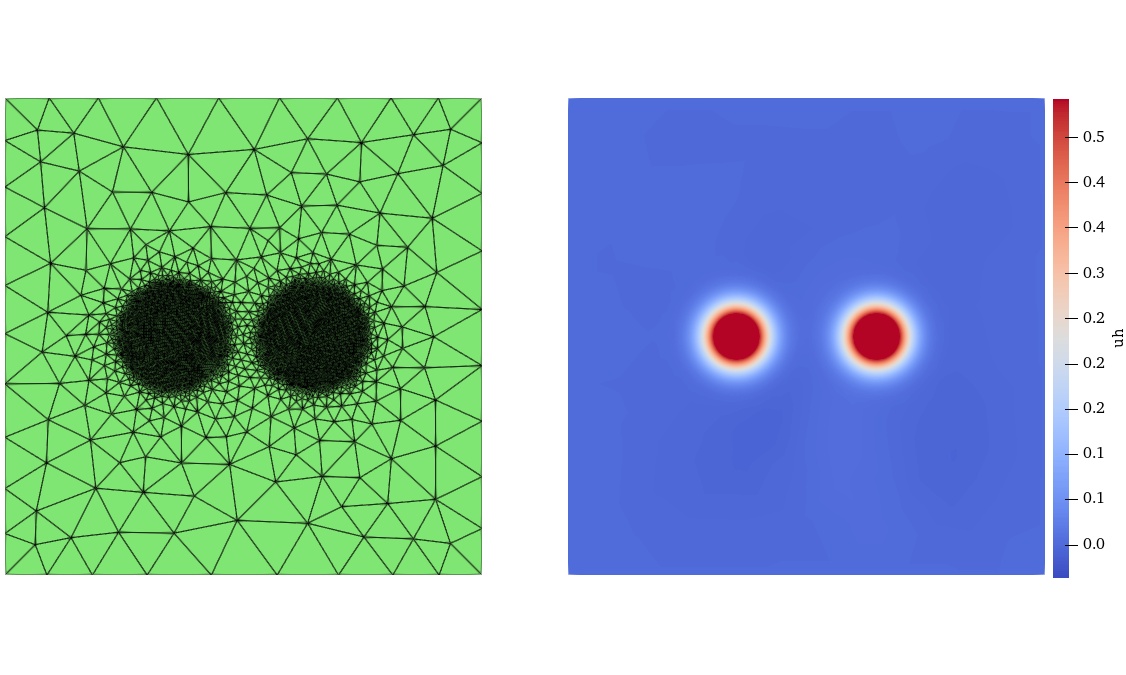}
			\includegraphics[width=0.45\linewidth]{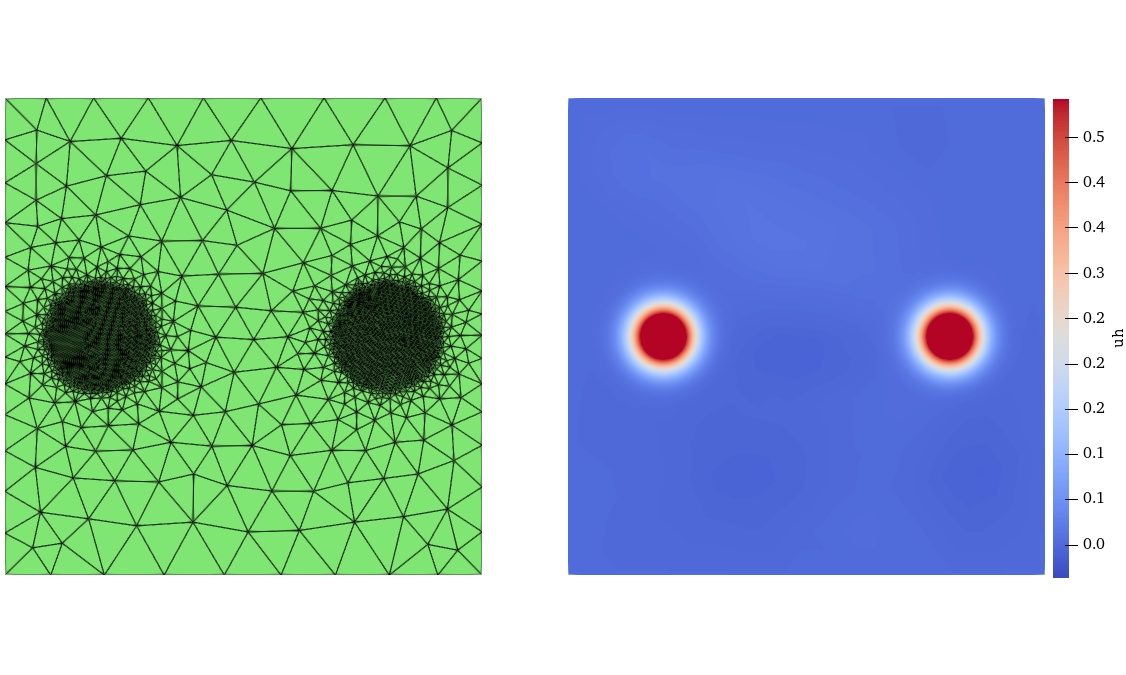}\vspace{-0.8cm}
			\caption{\Cref{example3}, snapshots of numerical solutions along with their corresponding adaptive meshes at $t=0.00, 0.25, 0.50$, and $t=0.75$.}
			\label{Ex2D3-Initial}
		\end{figure}
	\end{example}
	Consistent with \Cref{example1}, the computational settings remain unchanged with the tolerance $eTol=0.01$. 
	The evolution of the training loss across different time levels is shown in \Cref{Ex2D3-loss}. 
	To reach the desired MSE, the initial time level requires more than 500 training epochs. In contrast, fewer than 20 epochs are sufficient at all subsequent time levels to achieve a comparable level of accuracy. 
	\Cref{Ex2D3-Time} shows the adaptive process of the initial mesh over six iterations, with the NOV at each iteration being $118$, $206$, $401$, $813$, $1677$, and $3463$, respectively. 
	\Cref{Ex2D3-Initial} shows the evolution of the numerical solutions and the corresponding mesh at four different time steps. 
	The refinement process dynamically captures the splitting of the central peak into two distinct moving peaks, maintaining a coarser resolution in smoother regions while ensuring high resolution in areas with sharp gradients. 
	The right panel of \Cref{Ex2-order} shows that the adaptive algorithm reaches the quasi-optimal convergence rate in terms of the gradient error.

	\subsection{3D Example}
	Extending to three dimensions, the neural network is updated to the architecture $[3, 32, 32, 32, 32, 1]$ with the meaning as given earlier. 
	Owing to the lack of fully developed data structures for three-dimensional adaptive refinement in the \texttt{fealpy} framework, no direct three-dimensional comparisons are included in this study. 
	Nevertheless, this circumstance serves to highlight a salient advantage of the proposed algorithm: its adaptive mechanism is inherently dimension-agnostic. 
	Consequently, the entire implementation requires no modification when extending from 2D to 3D, underscoring the flexibility and robustness of our approach. 
	\begin{figure}
		\centering
		\includegraphics[width=0.3\linewidth, height=0.2\textheight]{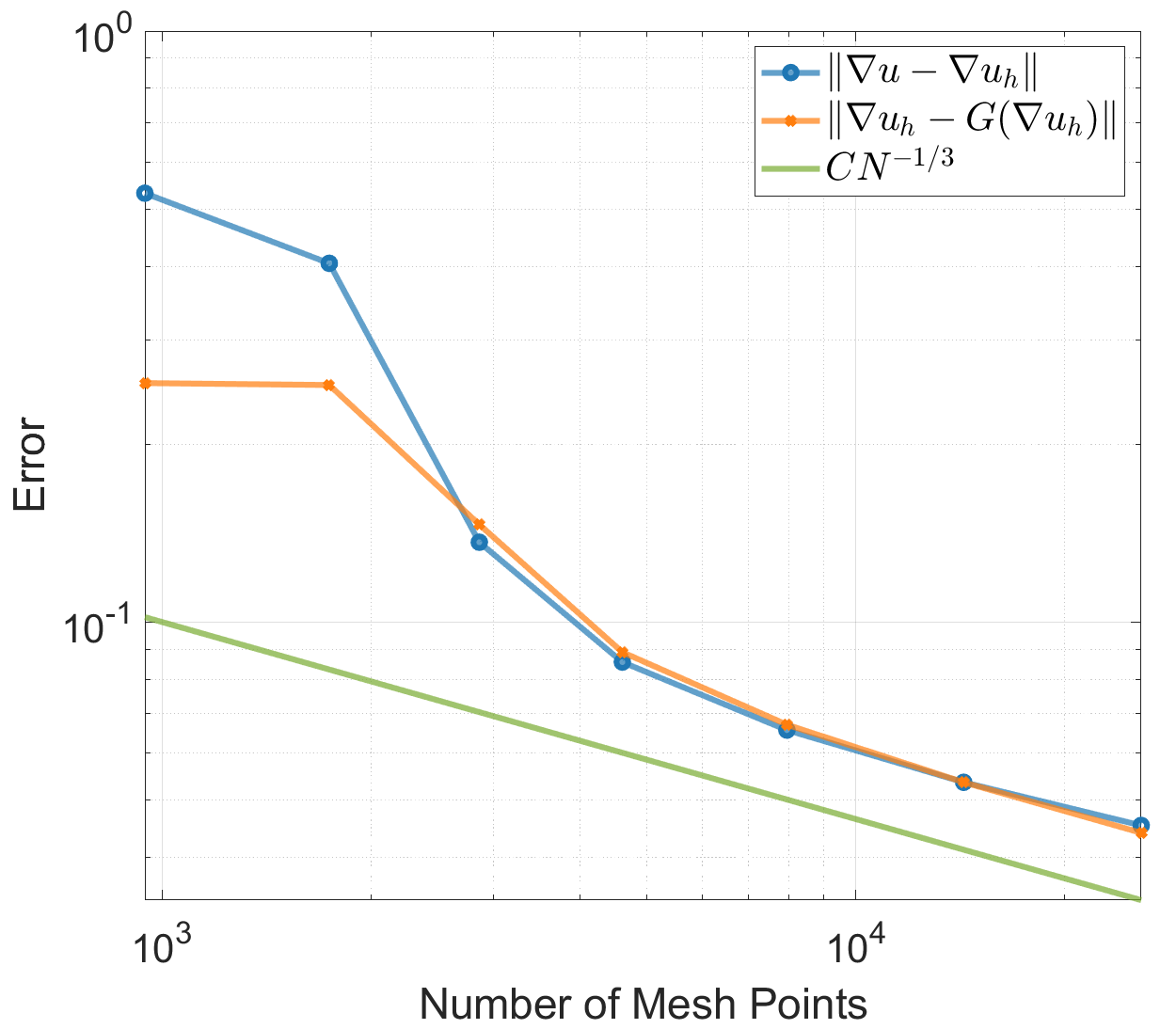}
		\includegraphics[width=0.3\linewidth, height=0.2\textheight]{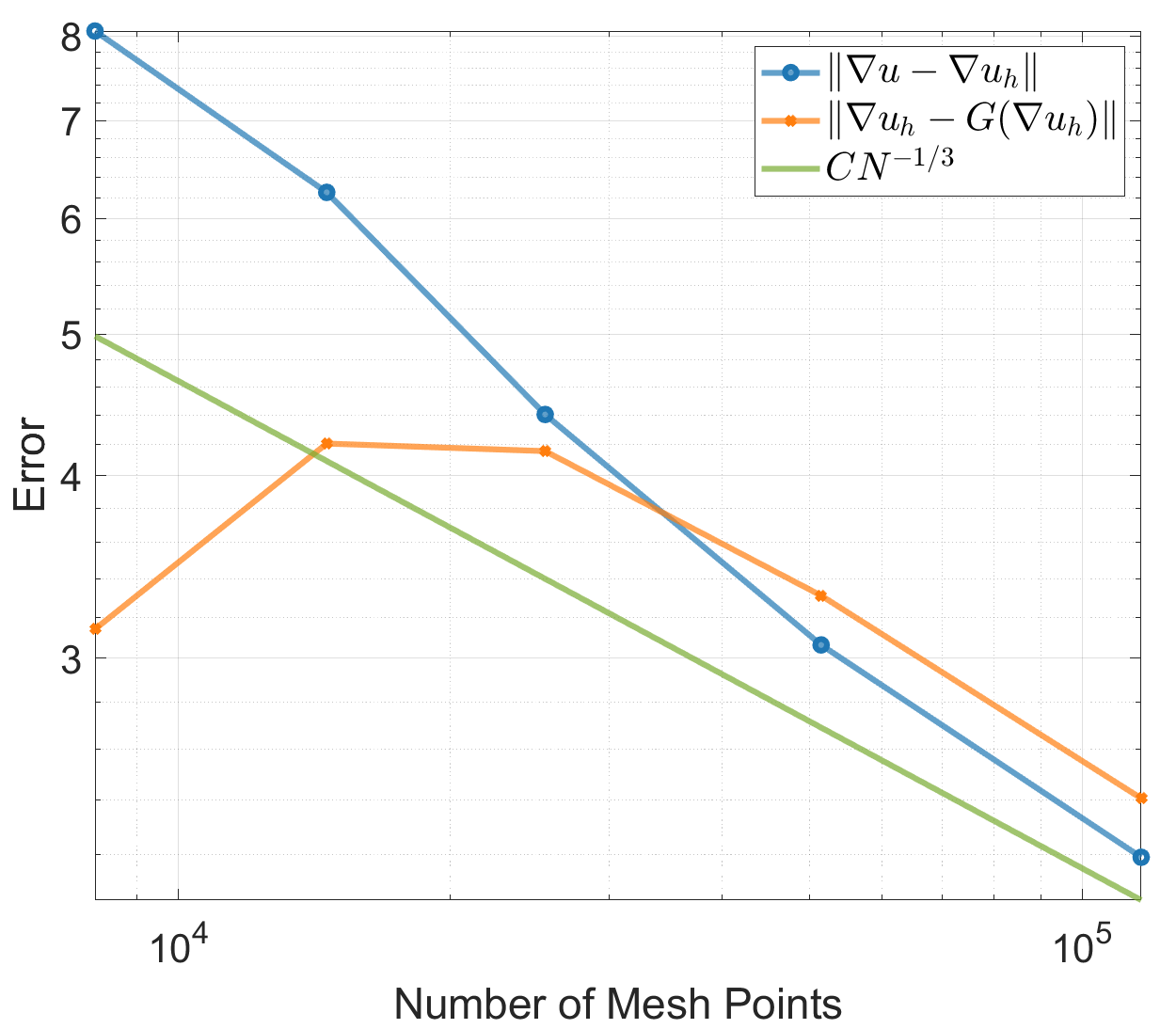}
		\includegraphics[width=0.3\linewidth, height=0.2\textheight]{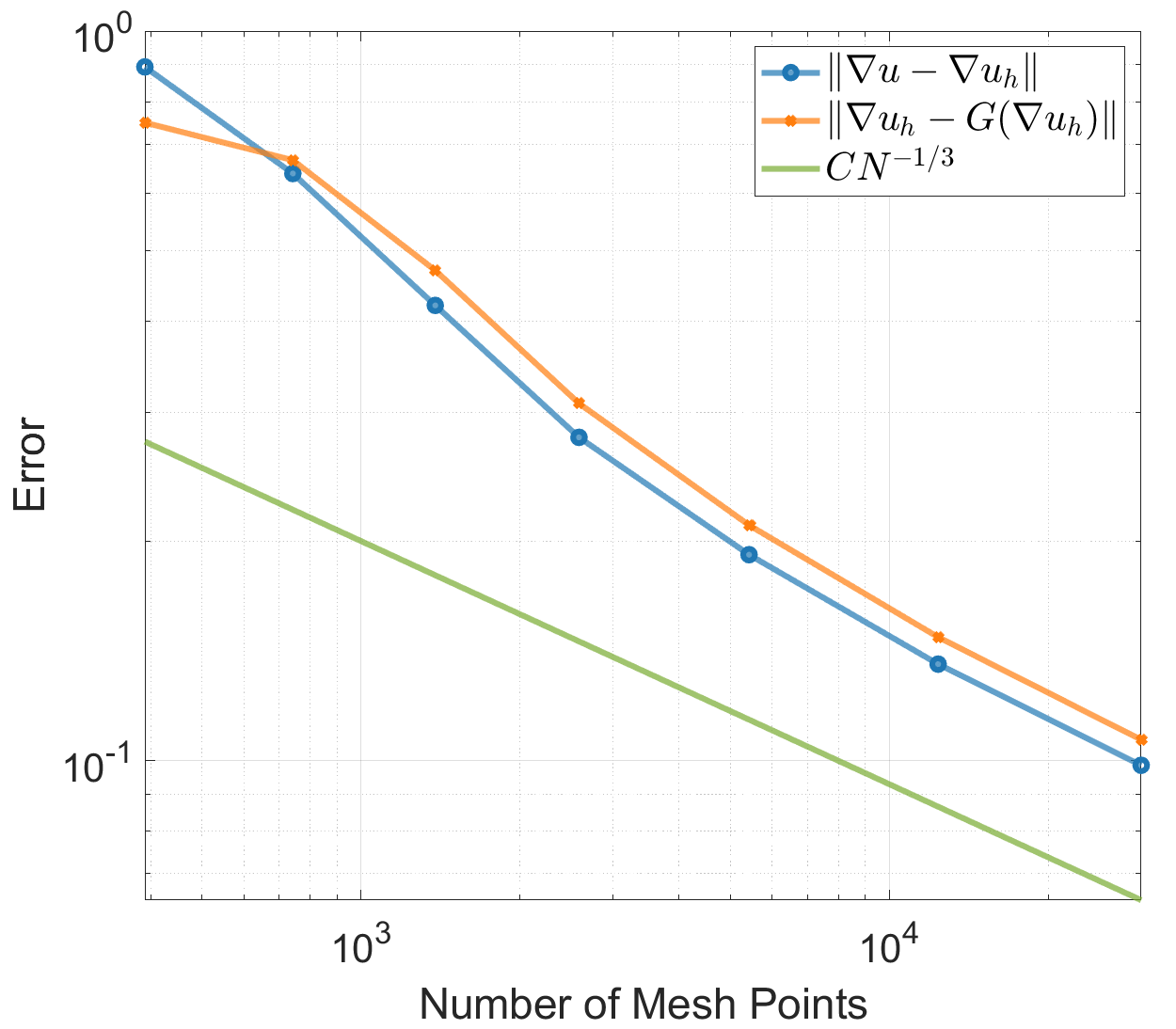}
		\caption{History of gradient error $\Vert\nabla u-\nabla u_h\Vert$ and recovered error estimator $\Vert\nabla u_h-G(\nabla u_h)\Vert$ at $T=1.0$ in \Cref{example4} (left), \Cref{example5} (middle) and \Cref{example6} (right).}
		\label{Ex3-order}
	\end{figure}
	
	\begin{example} \label{example4}
		(Rotation) Consider the 3D parabolic equation \eqref{model} in the domain $\Omega=[-1,1]^3$, with the source function $f$ such that the exact solution $u$ is given by:
		\[u(x,y,z,t) = \exp(-500(x-0.3\cos(2\pi t))^2)\exp(-500(y-0.3\sin(2\pi t))^2)\exp(-500z^2).\]
		\begin{figure}[!ht]
			\includegraphics[width=0.95\linewidth]{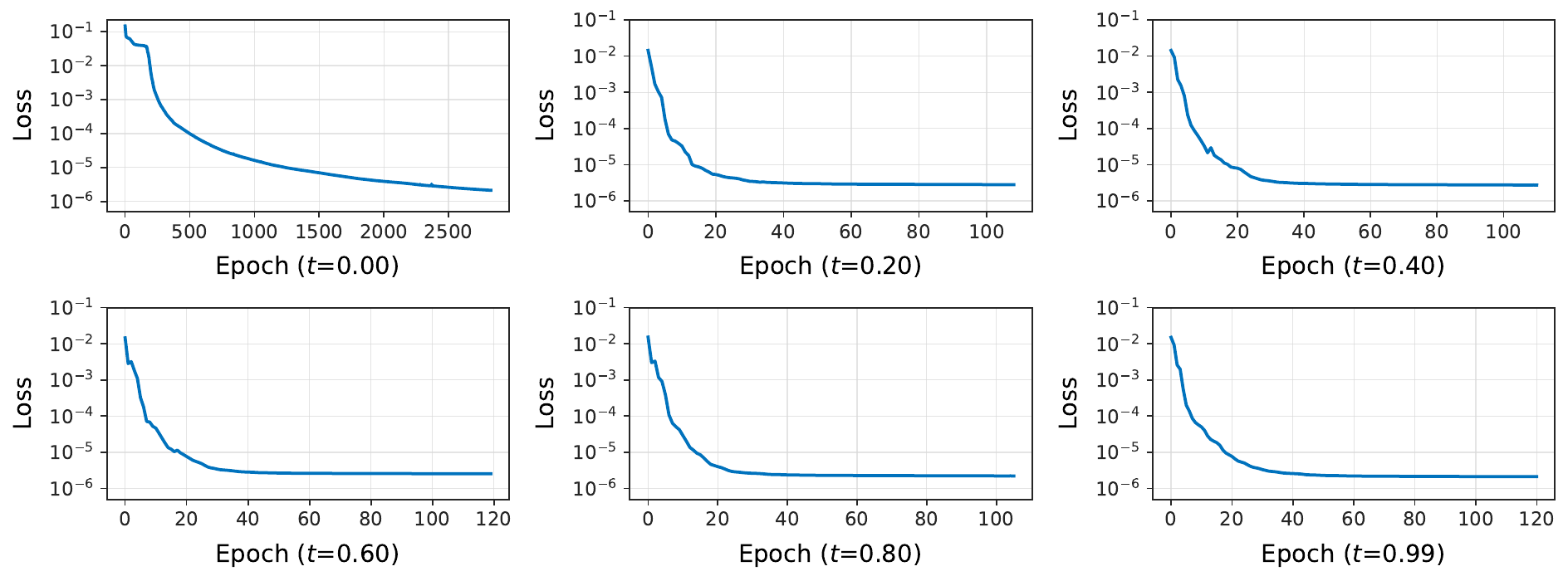}
			\caption{\Cref{example4}, training loss convergence curves at different time levels.}
			\label{Ex3D1-loss}
		\end{figure}
		\begin{figure}
			\centering
			\includegraphics[width=0.21\linewidth]{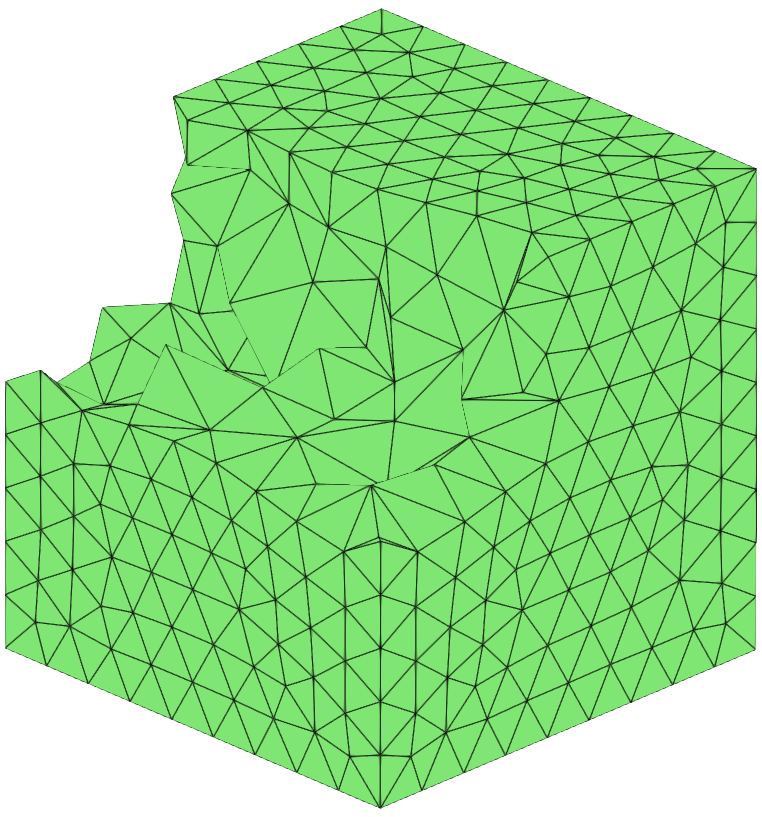}\hspace{0.8cm}
			\includegraphics[width=0.21\linewidth]{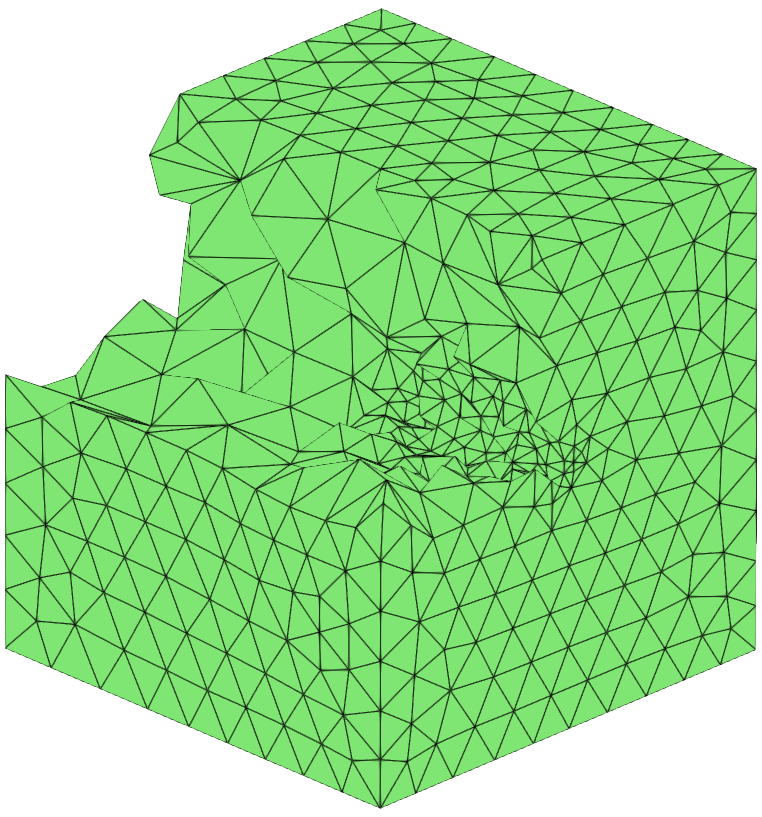}\hspace{0.8cm}
			\includegraphics[width=0.21\linewidth]{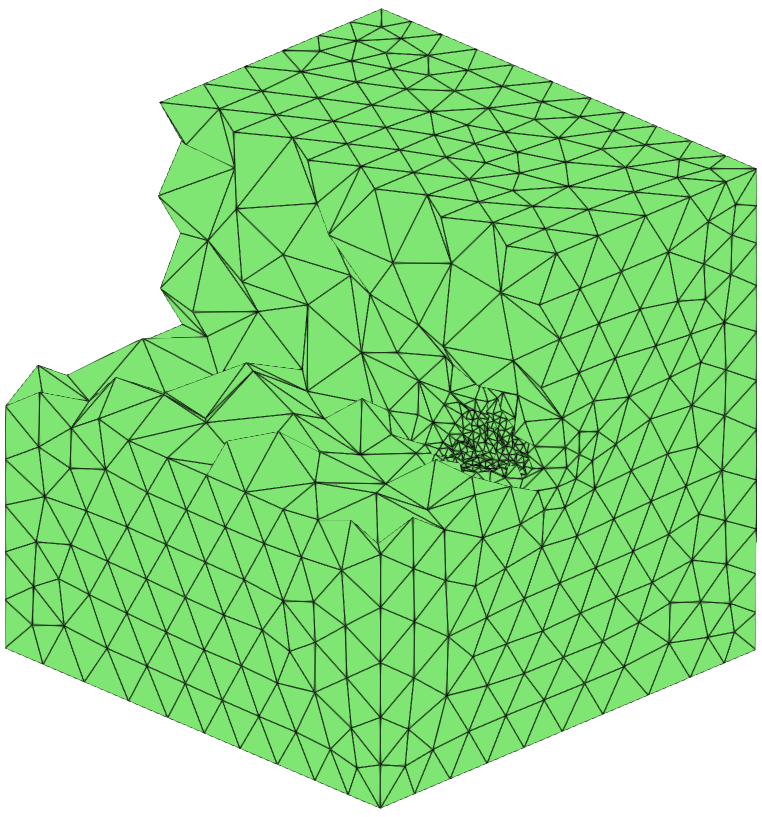}\\
			\includegraphics[width=0.21\linewidth]{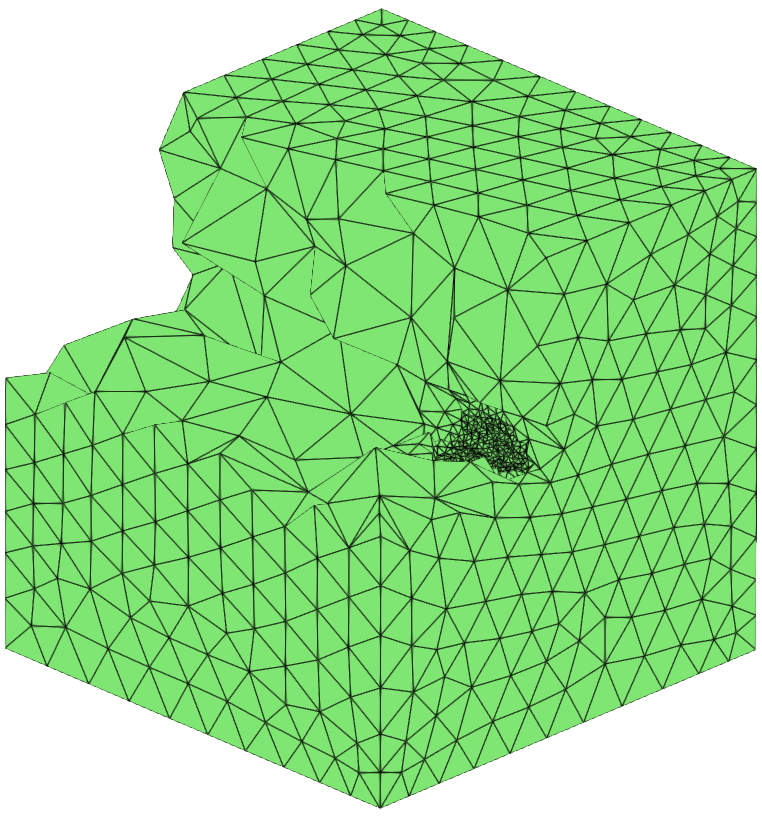}\hspace{0.8cm}
			\includegraphics[width=0.21\linewidth]{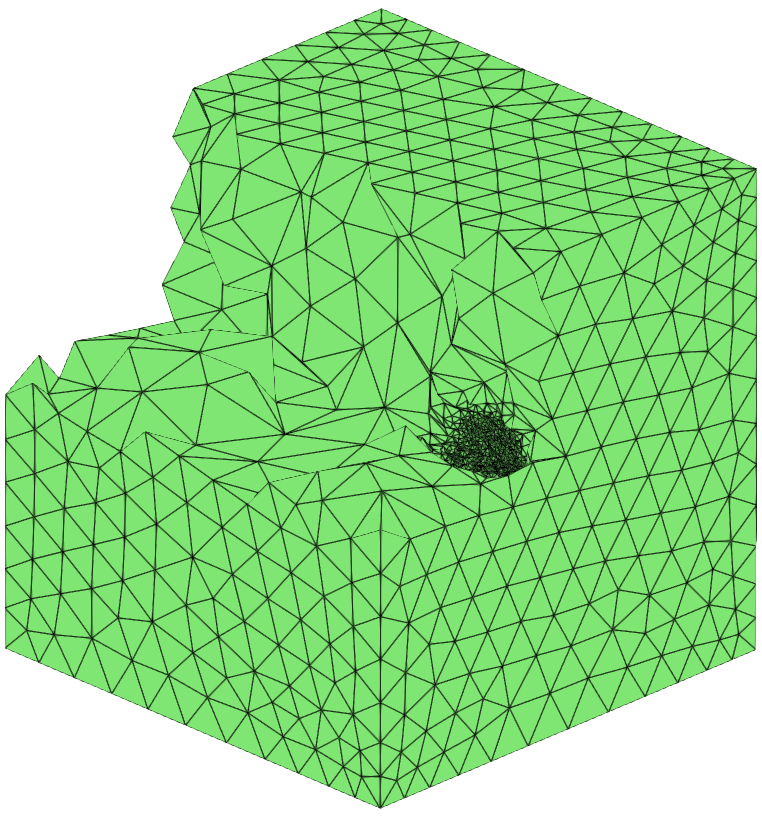}\hspace{0.8cm}
			\includegraphics[width=0.21\linewidth]{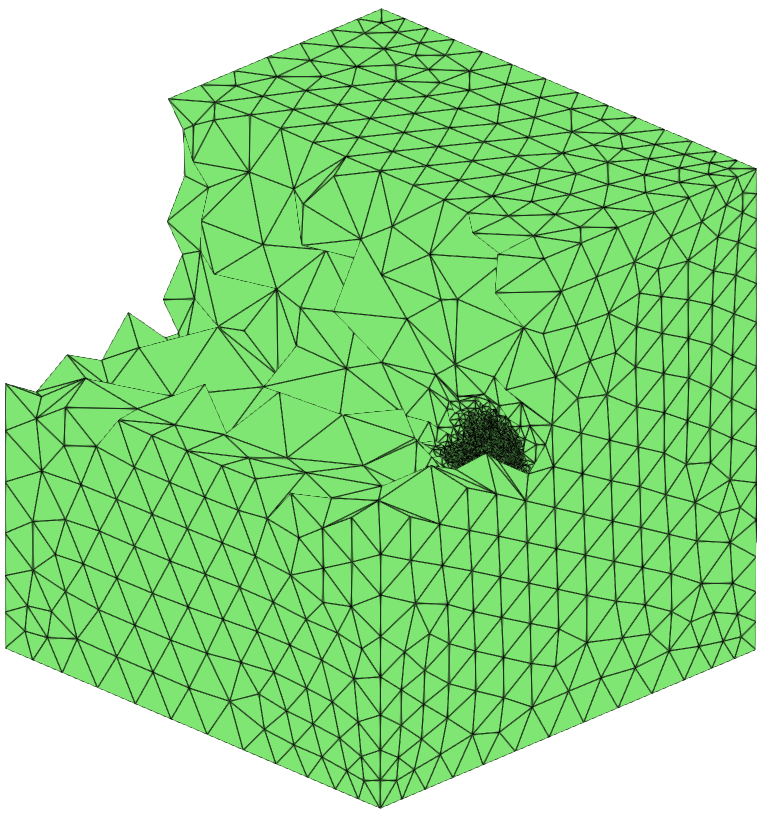}
			\caption{\Cref{example4}, initial mesh (top left) and its evolution through five adaptive refinements at $t=0.0$.}
			\label{Ex3D1-Time}
		\end{figure}
		\begin{figure}
			\centering
			\includegraphics[width=0.46\linewidth]{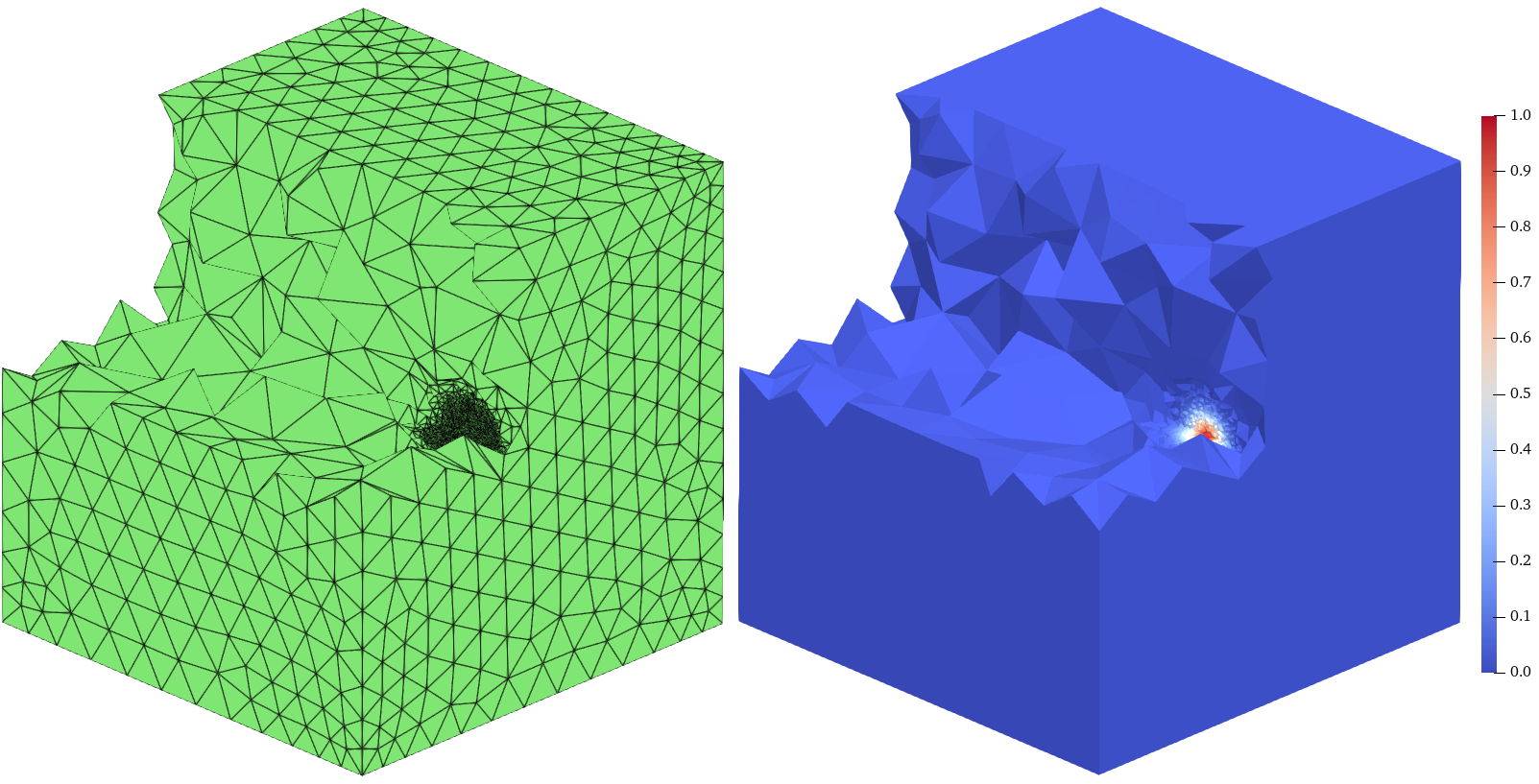}
			\includegraphics[width=0.46\linewidth]{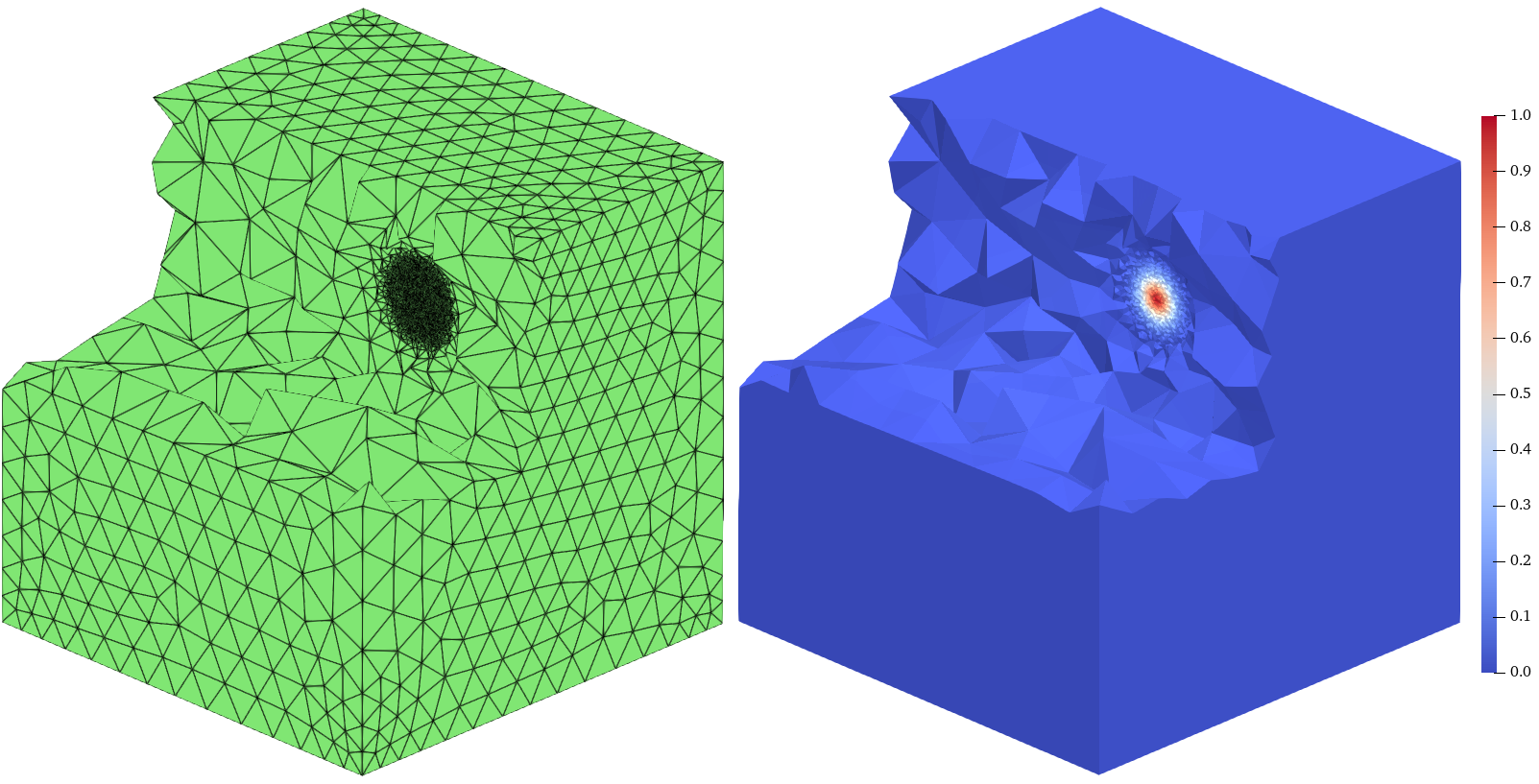}
			\includegraphics[width=0.46\linewidth]{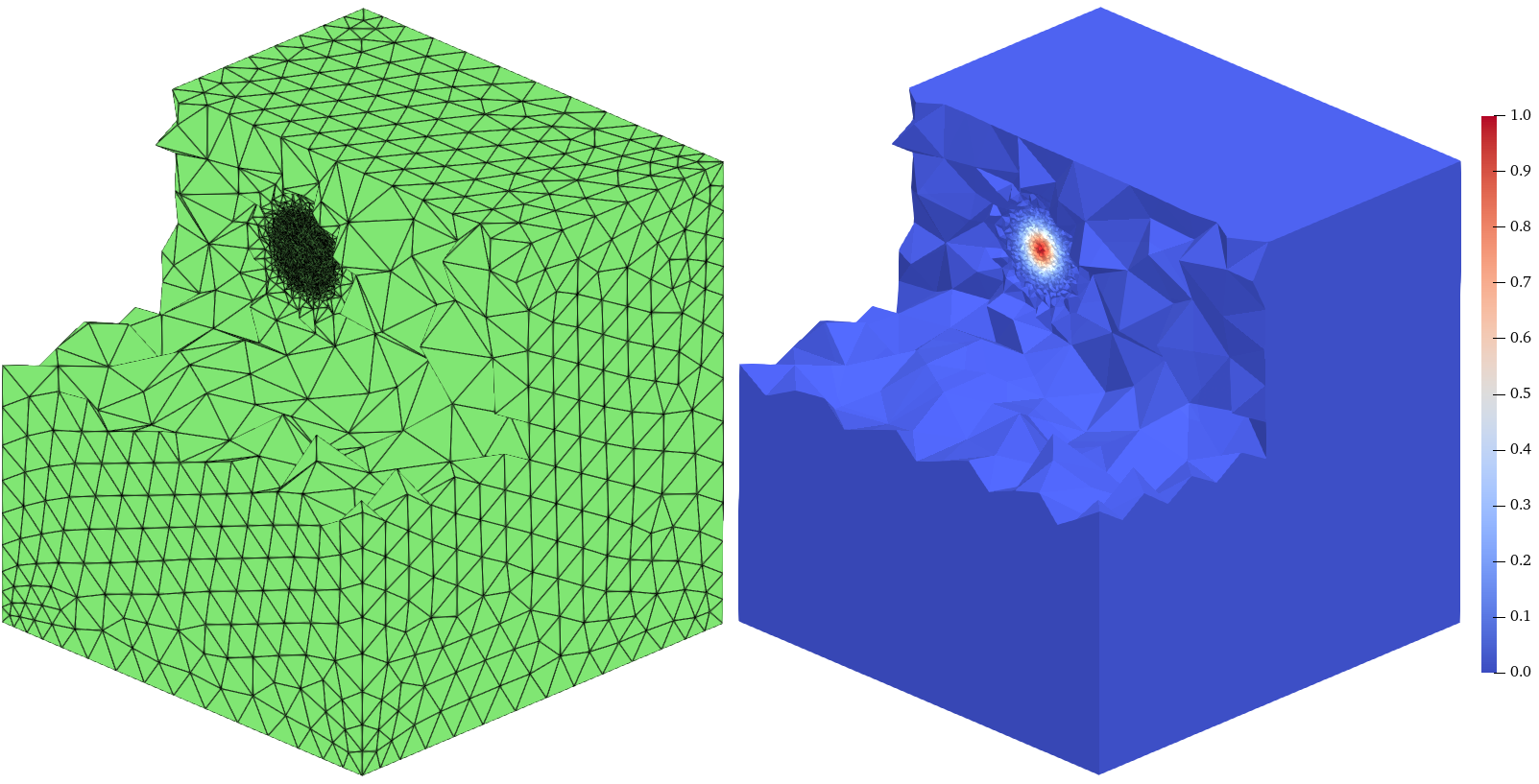}
			\includegraphics[width=0.46\linewidth]{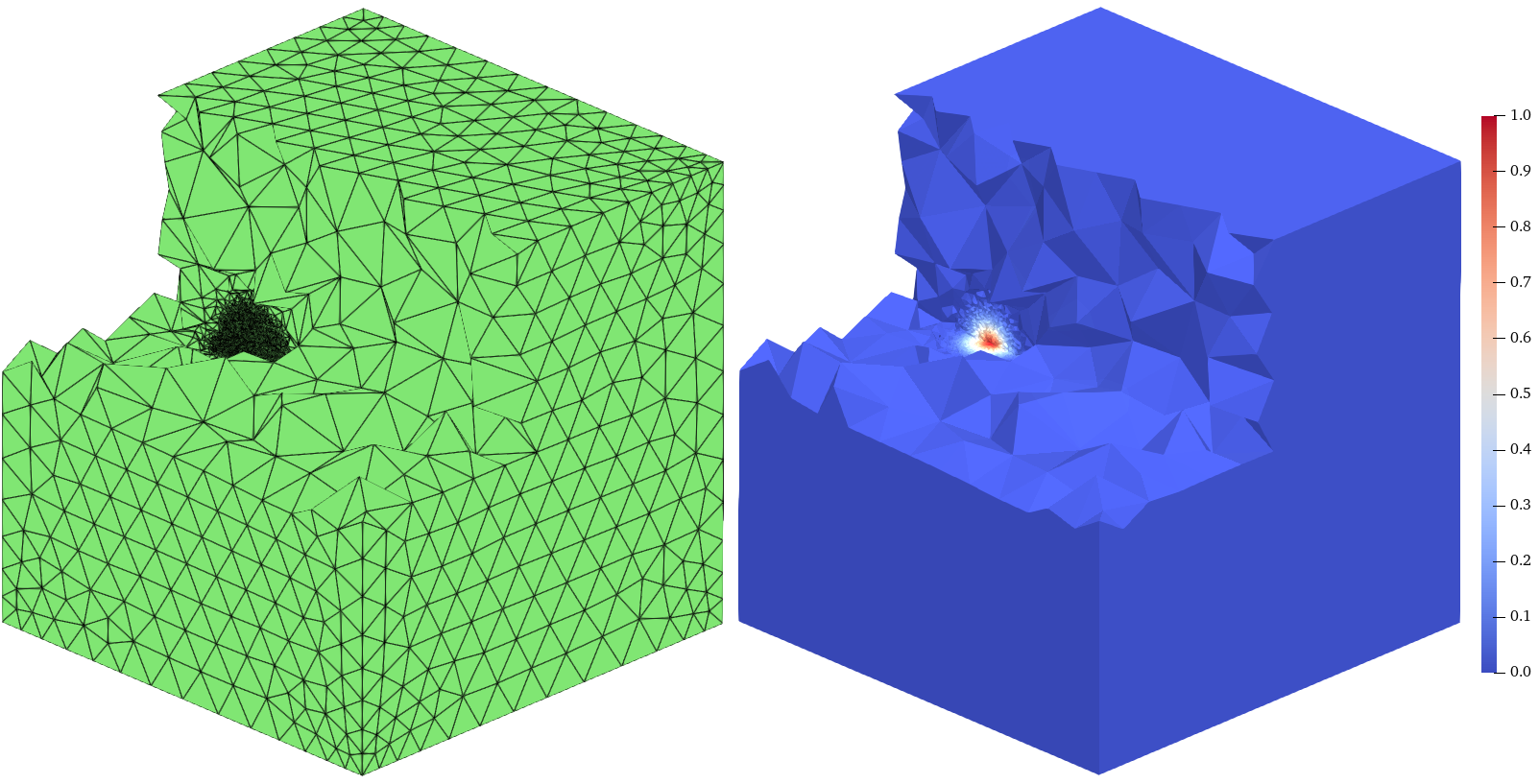}
			\caption{\Cref{example4}, snapshots of numerical solutions along with their corresponding adaptive meshes at $t=0.00, 0.17, 0.34$, and $t=0.50$.}
			\label{Ex3D1-Initial}
		\end{figure}
	\end{example}
	Apply \Cref{alg:ParabolicHR} with $eTol = 0.1$. 
	The evolution of the training loss across different time levels is shown in \Cref{Ex3D1-loss}. To reach the desired MSE, the initial time level requires more than 1000 training epochs. In contrast, fewer than 20 epochs are sufficient at all subsequent time levels to achieve a comparable level of accuracy. 
	\Cref{Ex3D1-Time} illustrates the mesh refinement process at initial time $t=0.0$, showing the transition from an initial uniform mesh to its evolution over five refinement steps. The NOV at each iteration are $944$, $1685$, $2819$, $4827$, $8733$, and $15908$, respectively. 
	\Cref{Ex3D1-Initial} presents the numerical results and corresponding meshes at four different time instances, demonstrating the adaptive method’s capability to efficiently handle time-dependent 3D problems with high precision. 
	The left panel of \Cref{Ex3-order} provides numerical evidence that the adaptive algorithm achieves the quasi-optimal convergence rate for the gradient error.

	\begin{example} \label{example5}
		(Diffusion) Consider the 3D parabolic equation \eqref{model} in the domain $\Omega=[-1,1]^3$, with the source function $f$ such that the exact solution $u$ is given by:
		\[u(x,y,z,t) = \exp(-5000(\sqrt{x^2+y^2+z^2}+0.3t-0.4)^2).\]
		
		\begin{figure}[!ht]
			\includegraphics[width=0.95\linewidth]{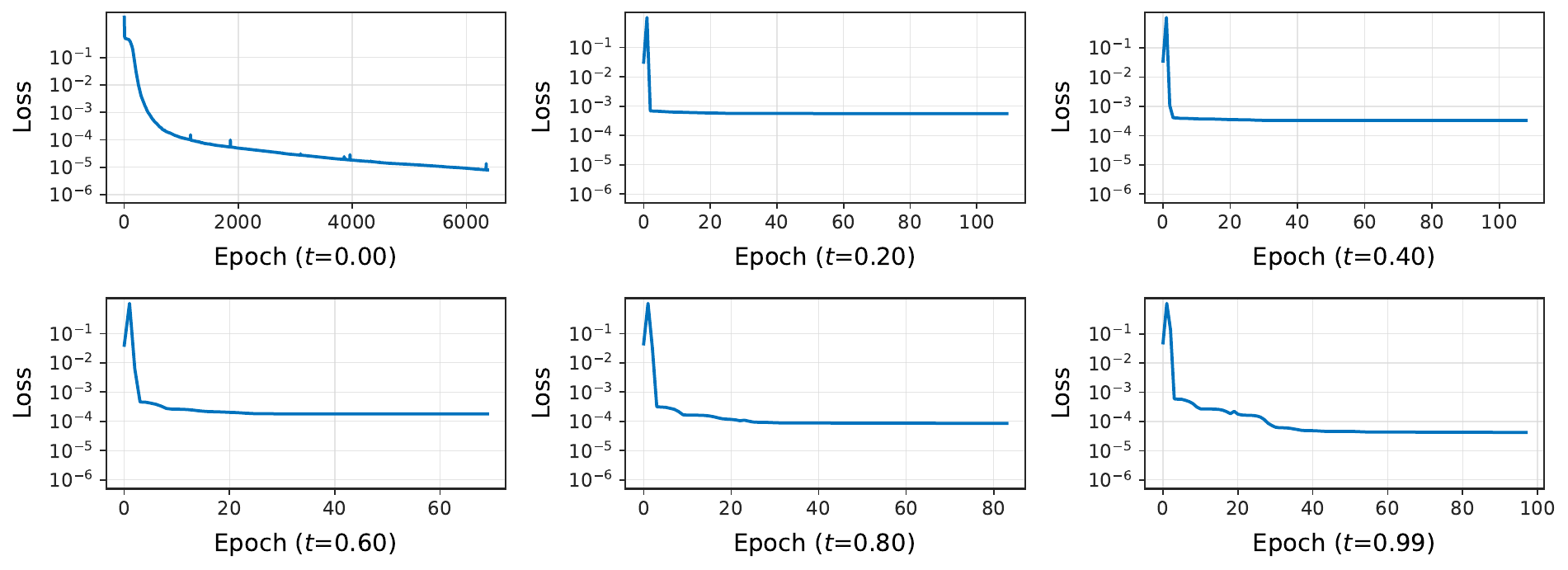}
			\caption{\Cref{example5}, training loss convergence curves at different time levels.}
			\label{Ex3D2-loss}
		\end{figure}
		\begin{figure}
			\begin{minipage}[c]{0.45\linewidth}
				\centering
				\includegraphics[width=0.48\linewidth]{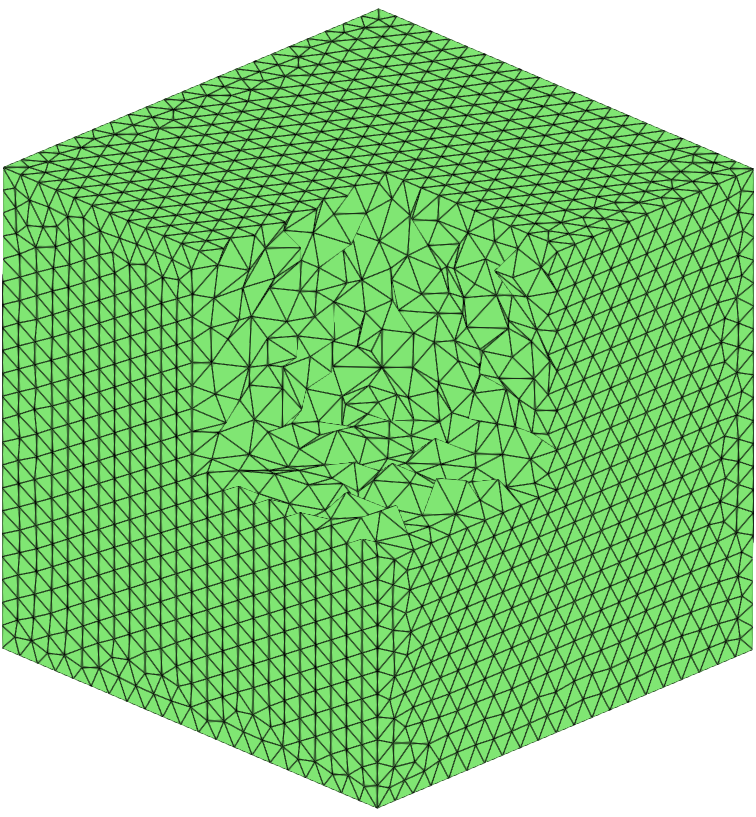}
			\end{minipage}%
			\begin{minipage}[c]{0.62\linewidth}
				\raggedright
				\includegraphics[width=0.35\linewidth]{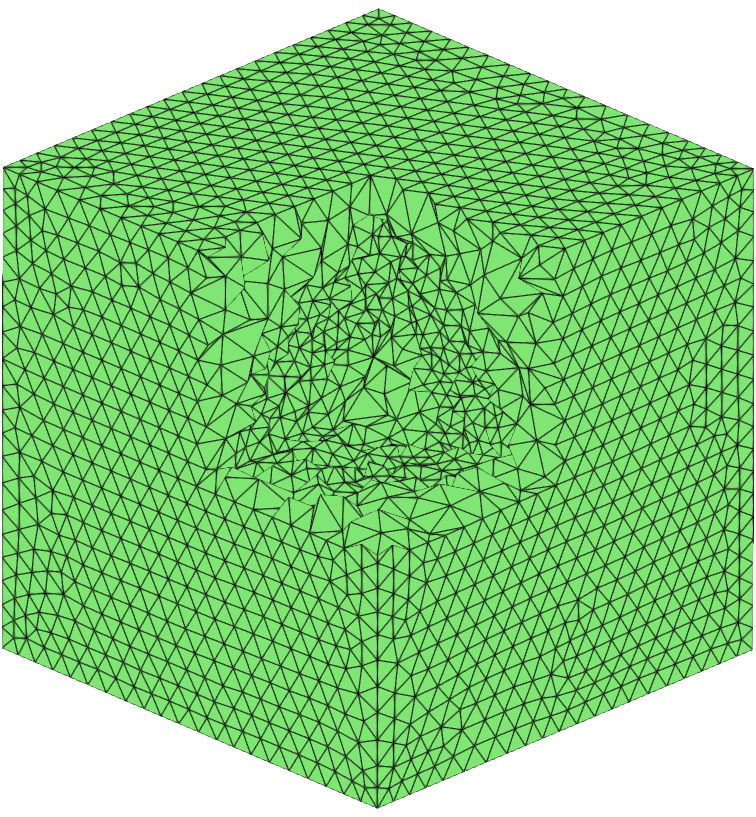}\hspace{0.5cm}
				\includegraphics[width=0.35\linewidth]{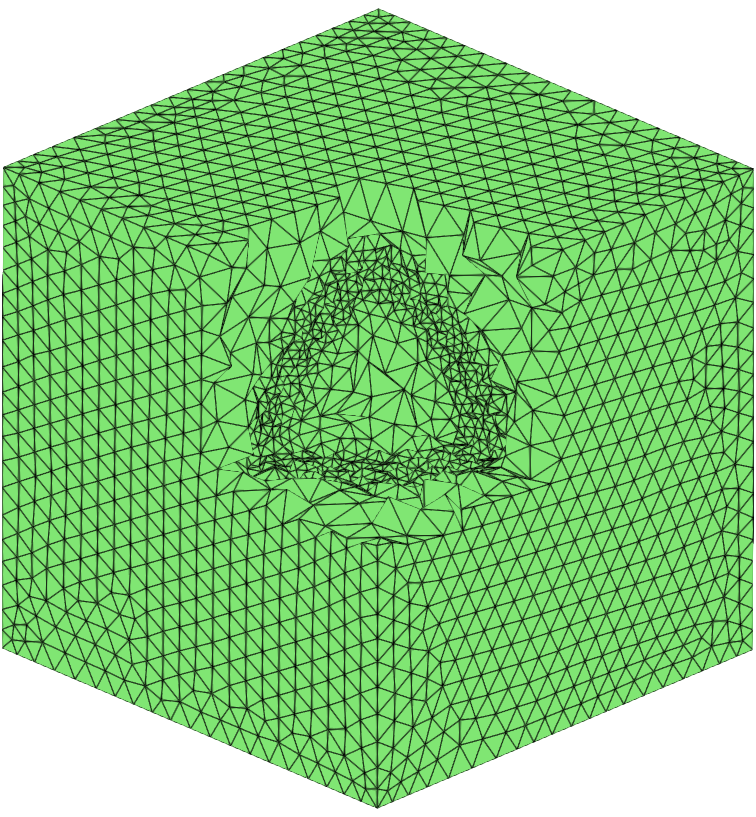}\\
				\includegraphics[width=0.35\linewidth]{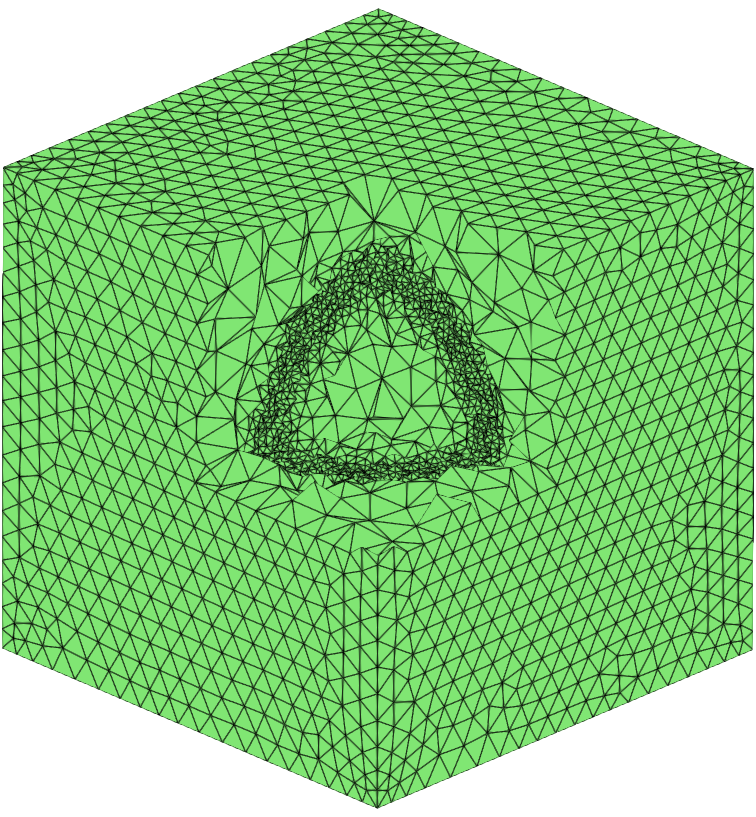}\hspace{0.5cm}
				\includegraphics[width=0.35\linewidth]{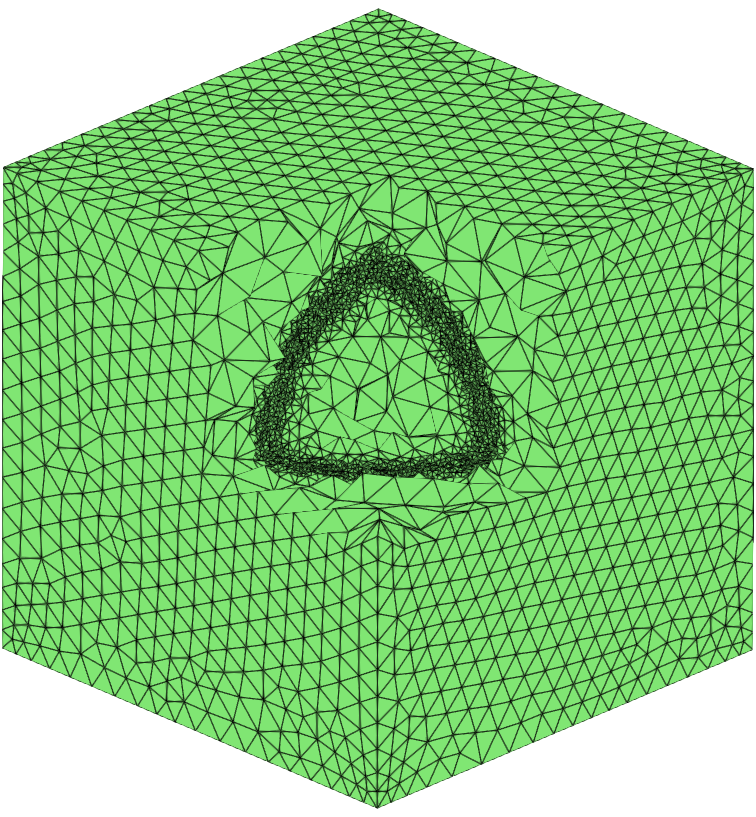}
			\end{minipage}
			\caption{\Cref{example5}, initial mesh (left) and its evolution through four adaptive refinements at $t=0.0$.}
			\label{Ex3D2-Time}
		\end{figure} 
		\begin{figure}
			\includegraphics[width=0.46\linewidth]{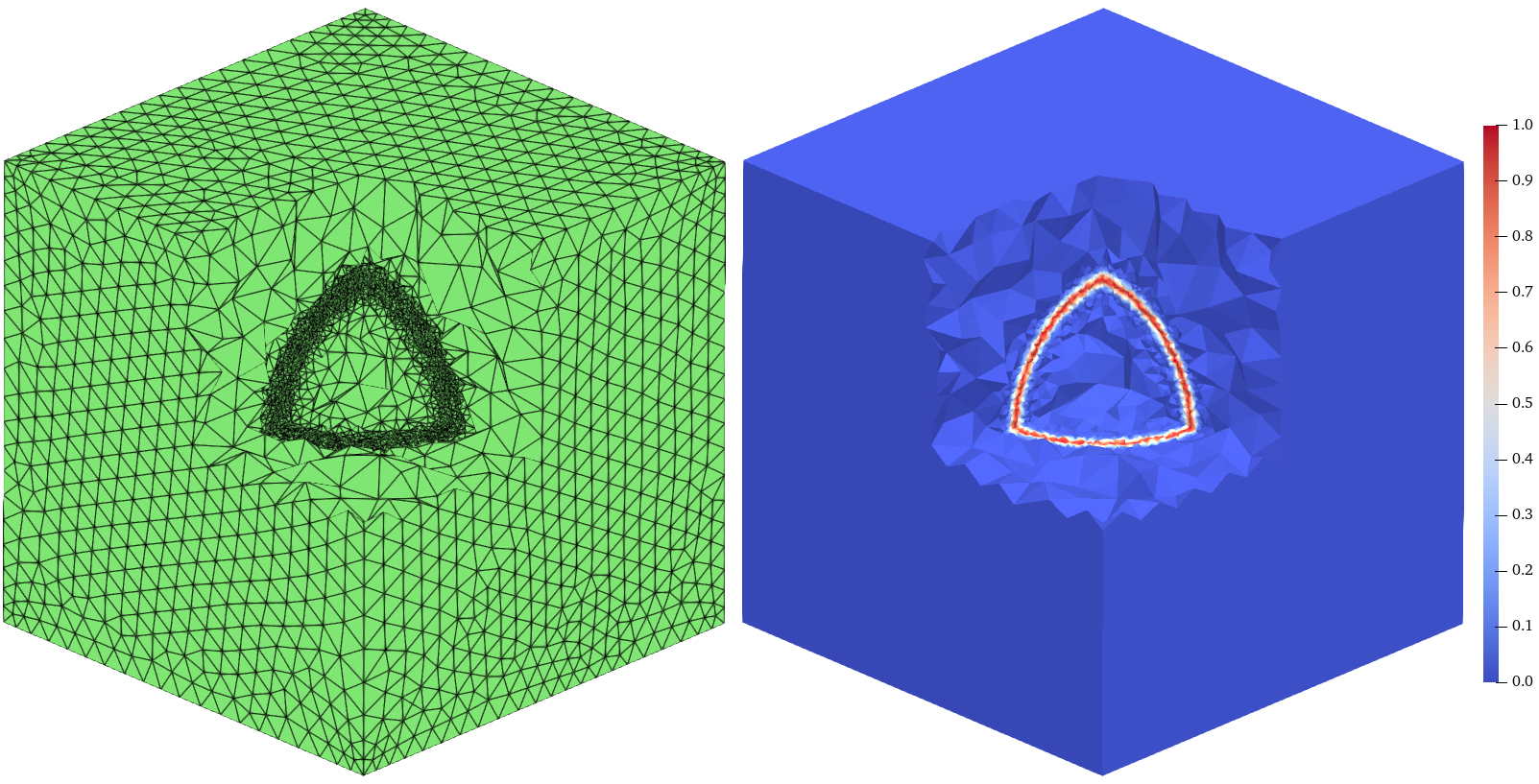}
			\includegraphics[width=0.46\linewidth]{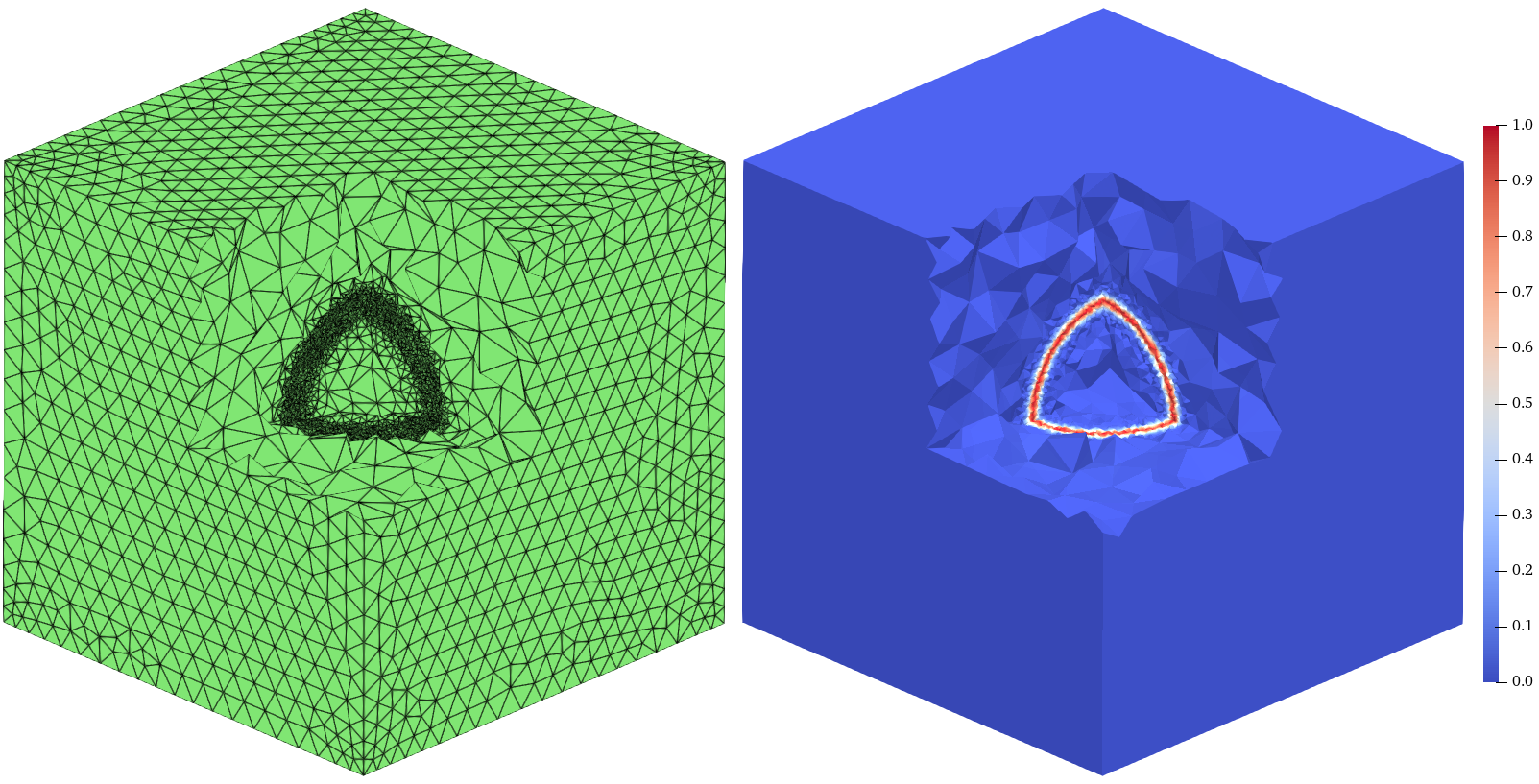}
			\includegraphics[width=0.46\linewidth]{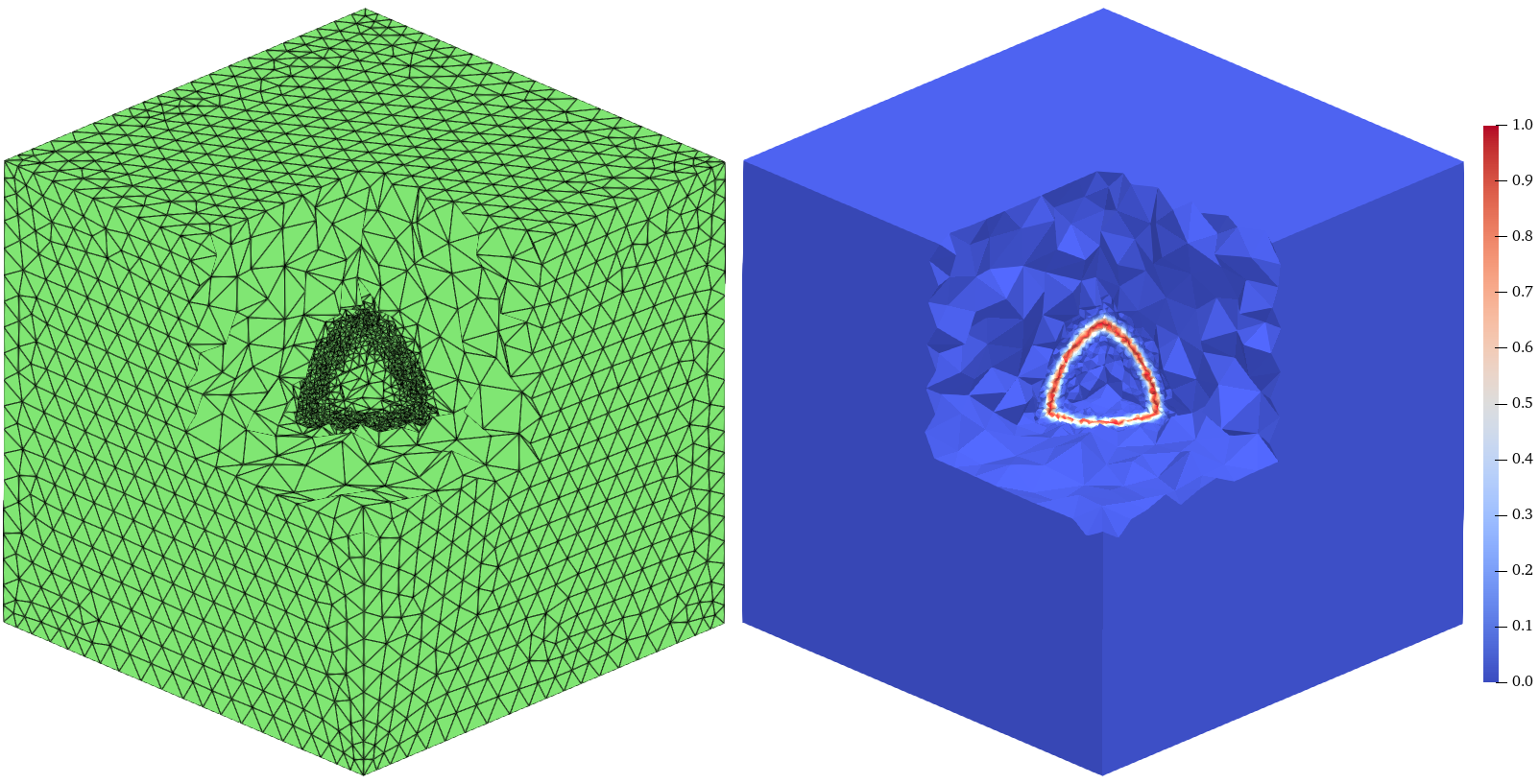}
			\includegraphics[width=0.46\linewidth]{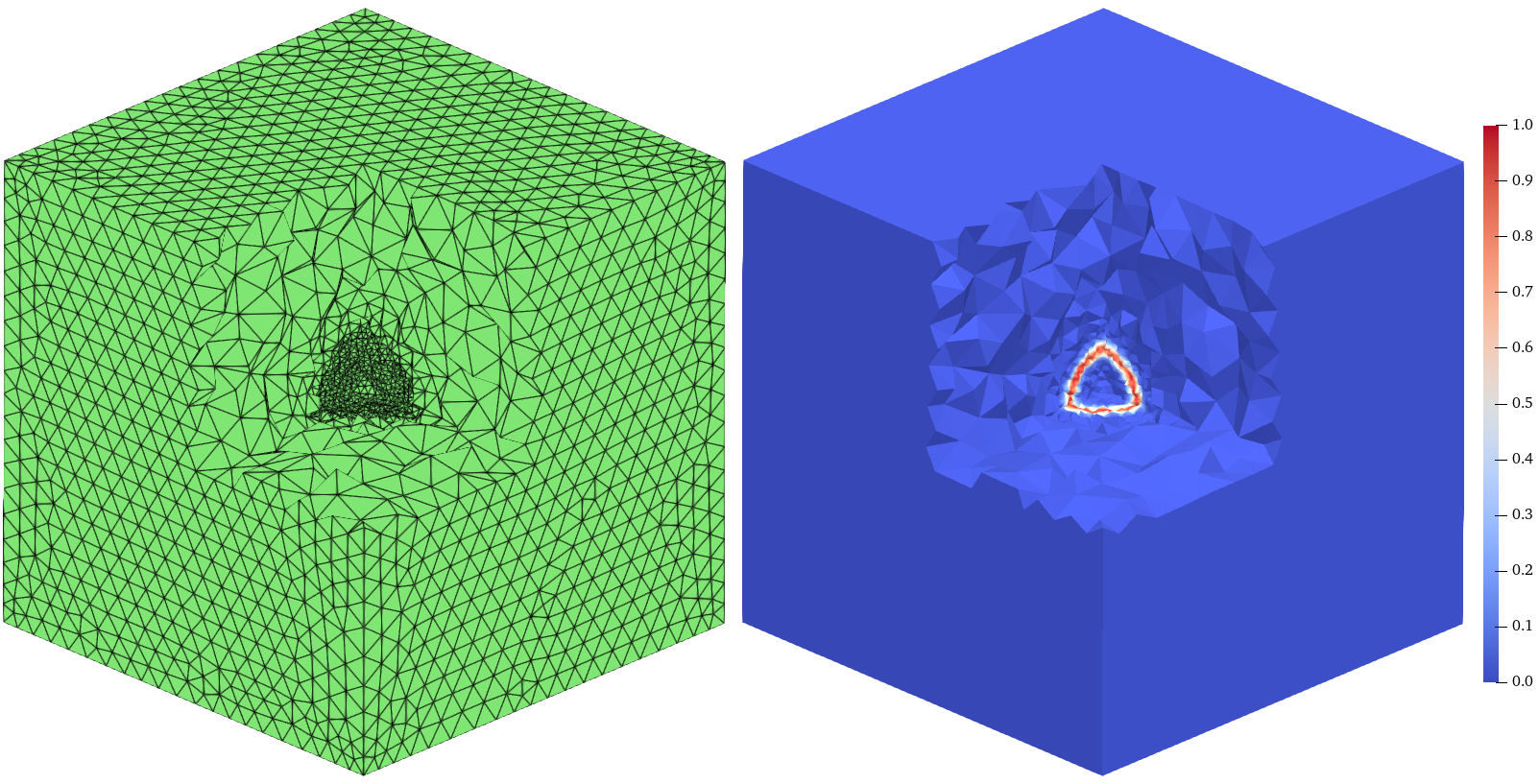}
			\caption{\Cref{example5}, snapshots of numerical solutions along with their corresponding adaptive meshes at $t=0.00, 0.25, 0.50$, and $t=0.75$.}
			\label{Ex3D2-Initial}
		\end{figure}
	\end{example}
	We set a relatively large mesh refinement tolerance, $eTol = 2.5$, to reduce the computational cost. 
	This choice leads to an earlier termination of the adaptive procedure, prior to the least-squares fitting step, while still allowing us to assess the effectiveness of both the adaptive algorithm and the neural network training strategy. 
	The evolution of the training loss across time levels is presented in \Cref{Ex3D2-loss}. 
	Unlike the previous examples, the loss in this case decreases only to the order of $10^{-4}$ for most time levels beyond the initial step. 
	This behavior is primarily due to the inherent noise in the training data resulting from the relatively large adaptive tolerance, which limits further convergence.
	Nevertheless, the training loss remains well below the prescribed tolerance, indicating that the achieved accuracy is sufficient to reliably advance the computation. 
	\Cref{Ex3D2-Time} illustrates the adaptive process of the initial mesh over four refinements, listing the NOV at each iteration as $8080$, $14502$, $27327$, $56739$, and $128527$, respectively. 
	\Cref{Ex3D2-Initial} presents the evolution of the numerical solutions and corresponding meshes at four distinct time points, with refinement specifically focusing on the contracting spherical shell.
	As illustrated in the middle panel of \Cref{Ex3-order}, the adaptive algorithm attains the quasi-optimal convergence rate for the gradient error.

	\begin{example} \label{example6}
		(Splitting) Consider the 3D parabolic equation \eqref{model} in the domain $\Omega=[-1,1]^3$, with the source function $f$ such that the exact solution $u$ is given by:
		\[u(x,y,z,t) = \exp(-300((x-0.3t)^2+y^2+z^2)) + \exp(-300((x+0.3t)^2+y^2+z^2)).\]
		
		\begin{figure}[!ht]
			\includegraphics[width=0.95\linewidth]{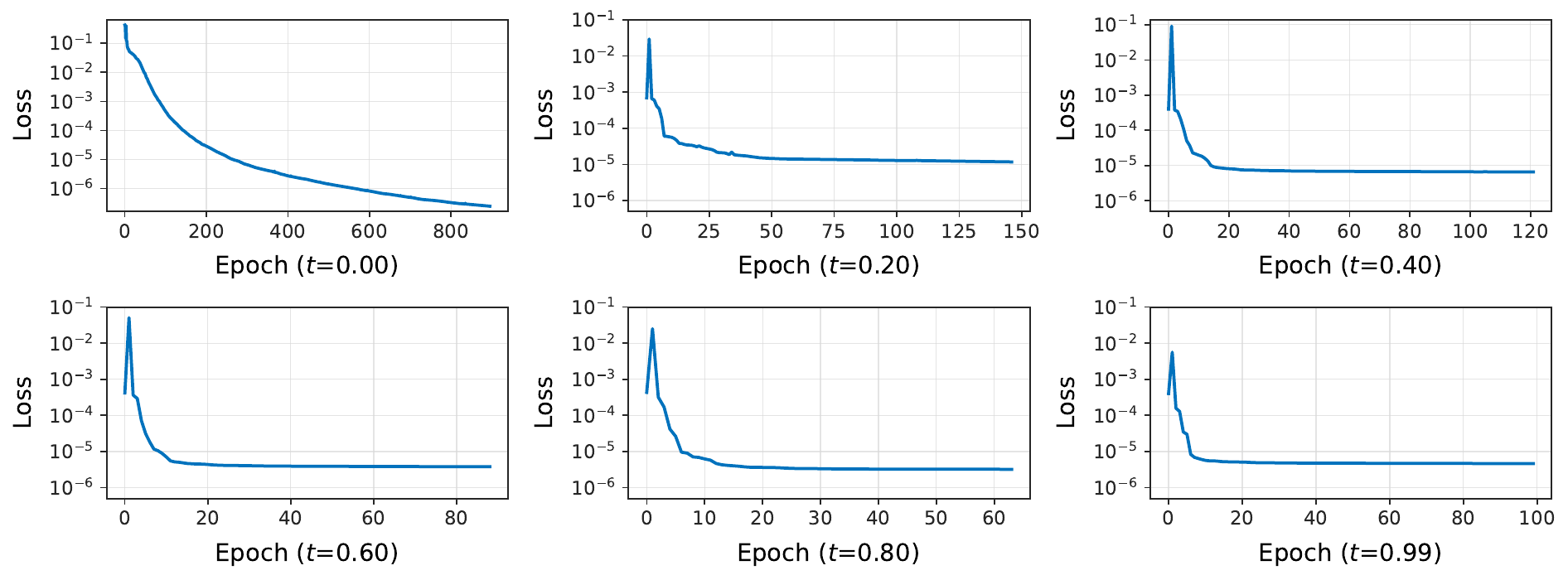}
			\caption{\Cref{example6}, training loss convergence curves at different time levels.}
			\label{Ex3D3-loss}
		\end{figure}
		\begin{figure}
			\begin{minipage}[c]{0.21\linewidth}
				\centering
				\includegraphics[width=\linewidth]{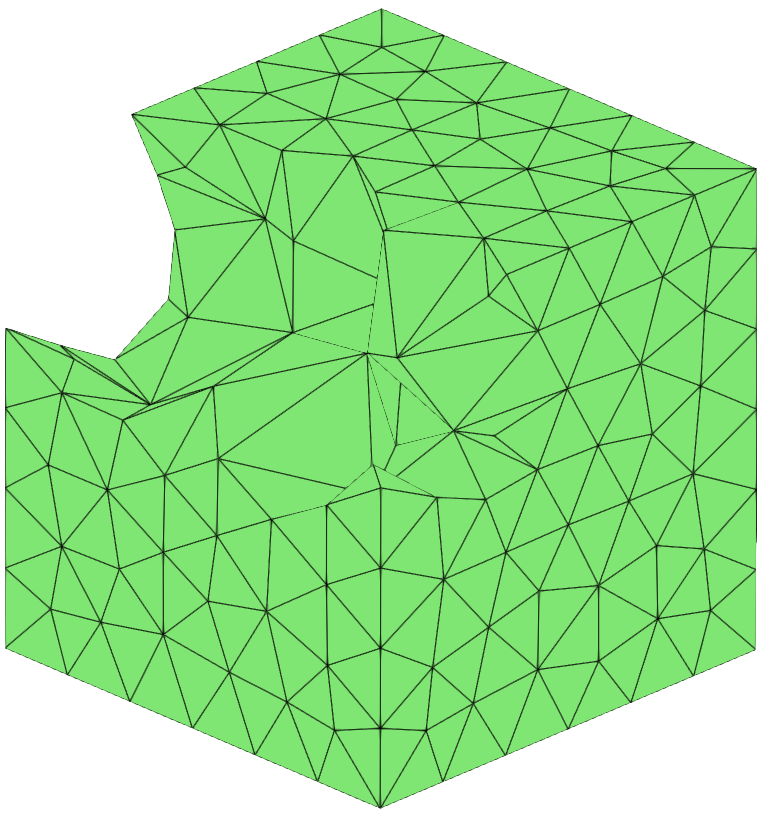}
			\end{minipage}%
			\begin{minipage}[c]{0.73\linewidth}
				\centering
				\includegraphics[width=0.3\linewidth]{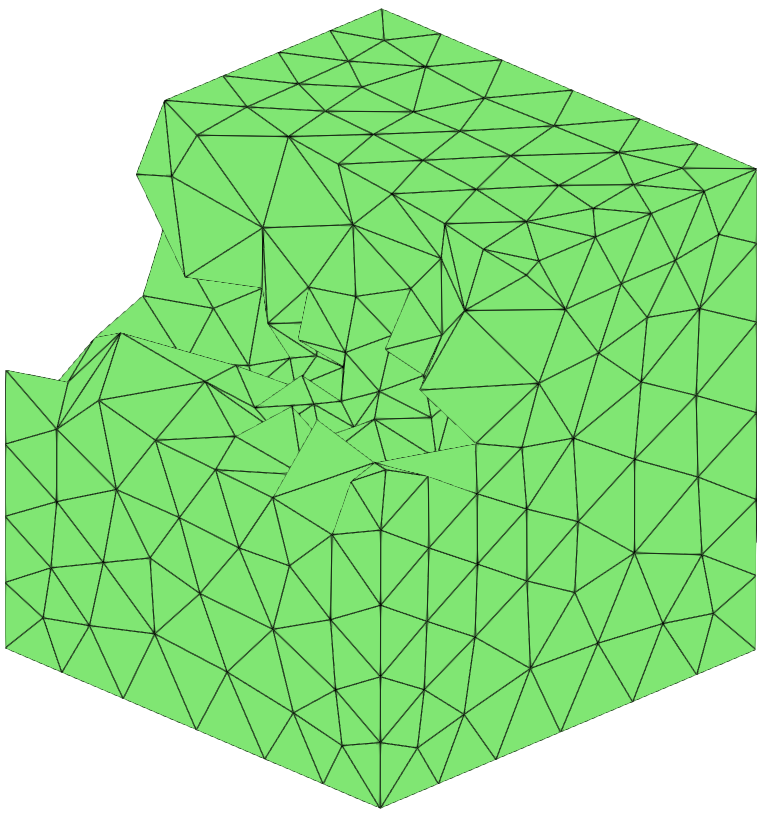}
				\includegraphics[width=0.3\linewidth]{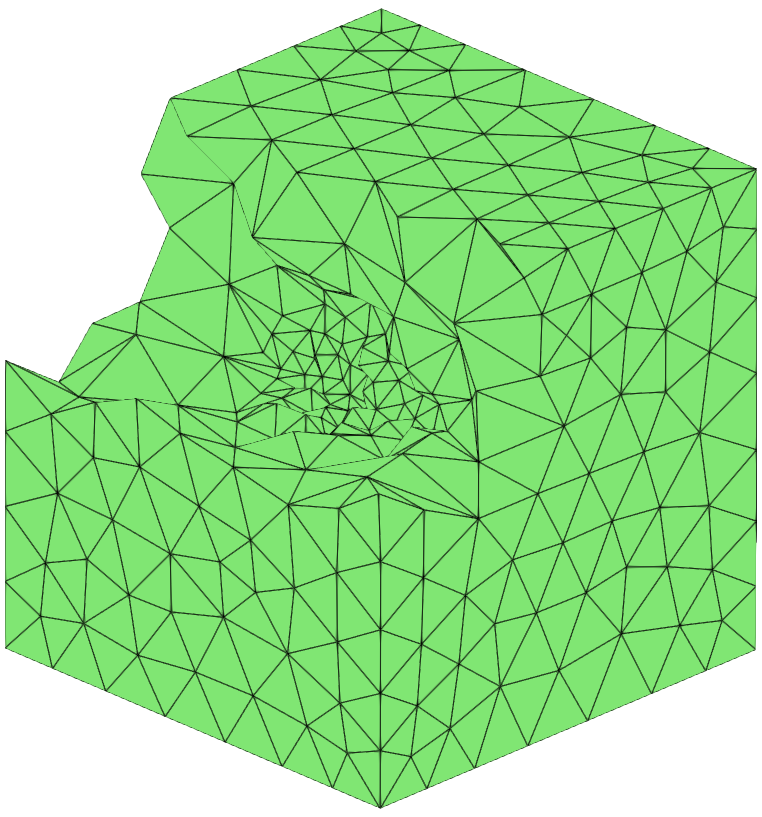}
				\includegraphics[width=0.3\linewidth]{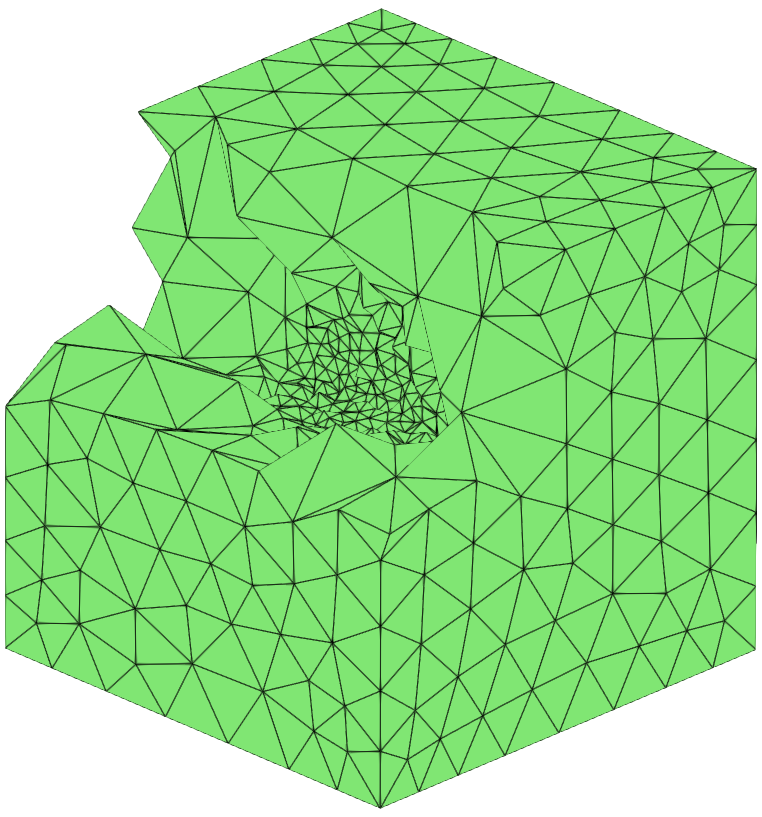} \\
				\includegraphics[width=0.3\linewidth]{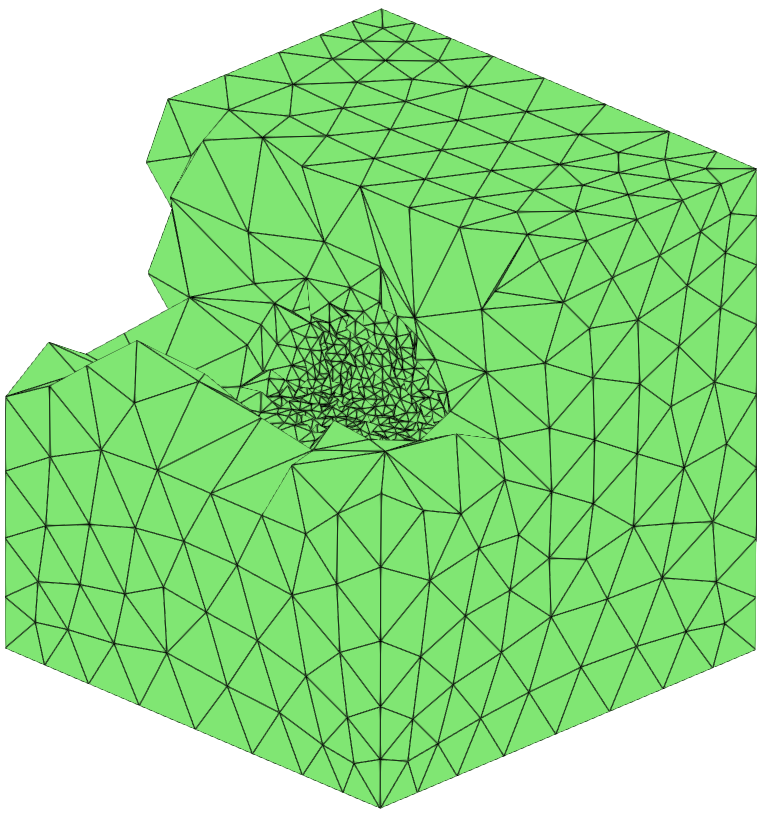}
				\includegraphics[width=0.3\linewidth]{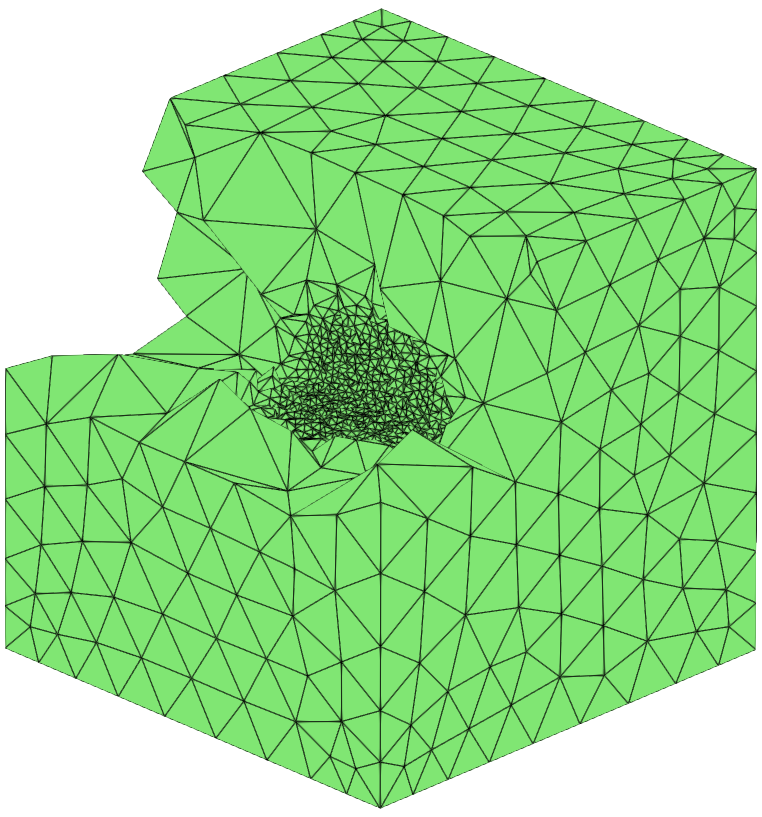}
				\includegraphics[width=0.3\linewidth]{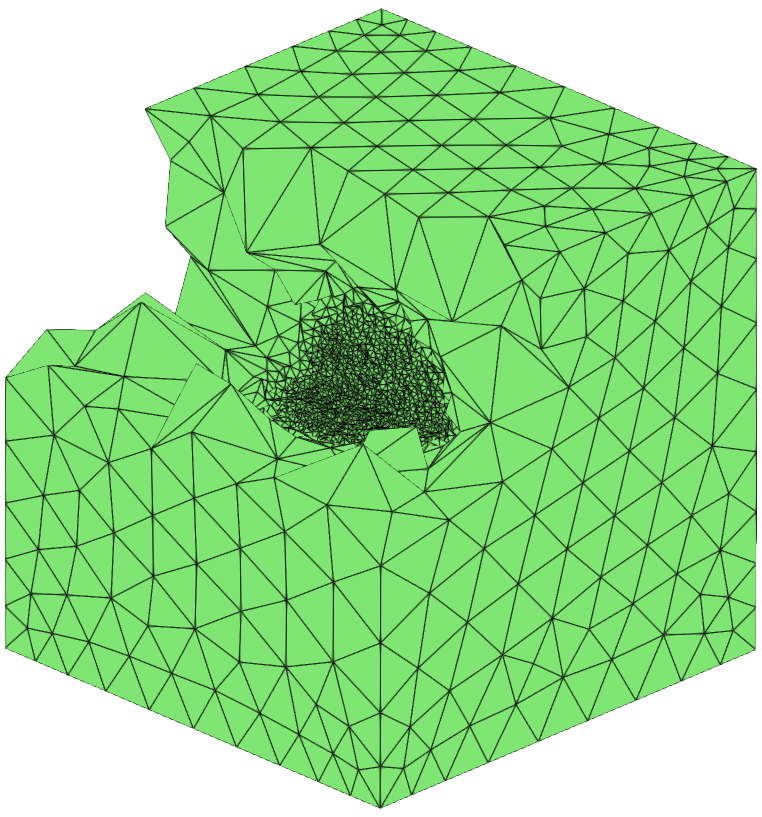}
			\end{minipage}
			\caption{\Cref{example6}, initial mesh (left) and its evolution through five adaptive refinements at $t=0.0$.}
			\label{Ex3D3-Time}
		\end{figure}   
		\begin{figure}
			\centering
			\includegraphics[width=0.47\linewidth]{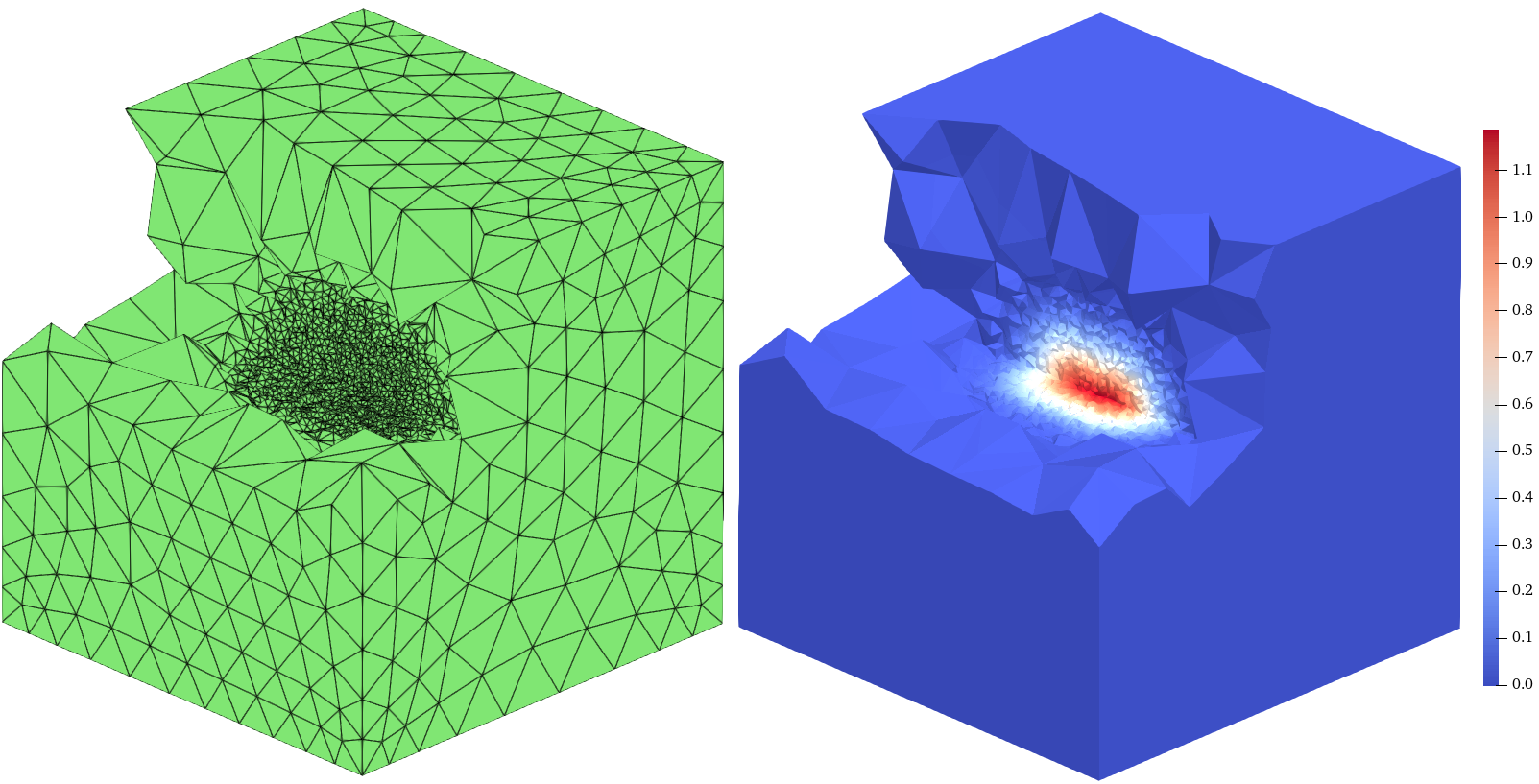}
			\includegraphics[width=0.47\linewidth]{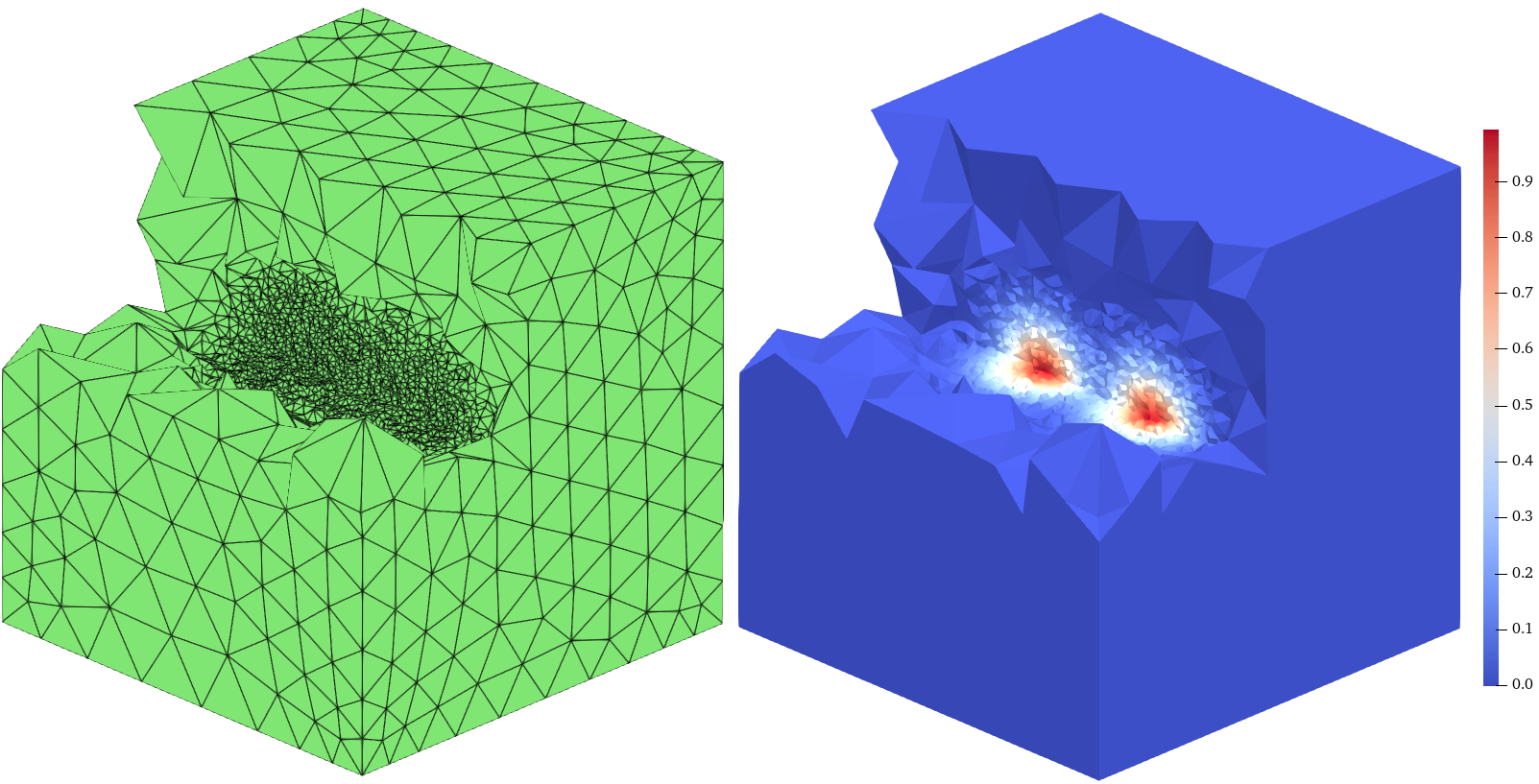}
			\includegraphics[width=0.47\linewidth]{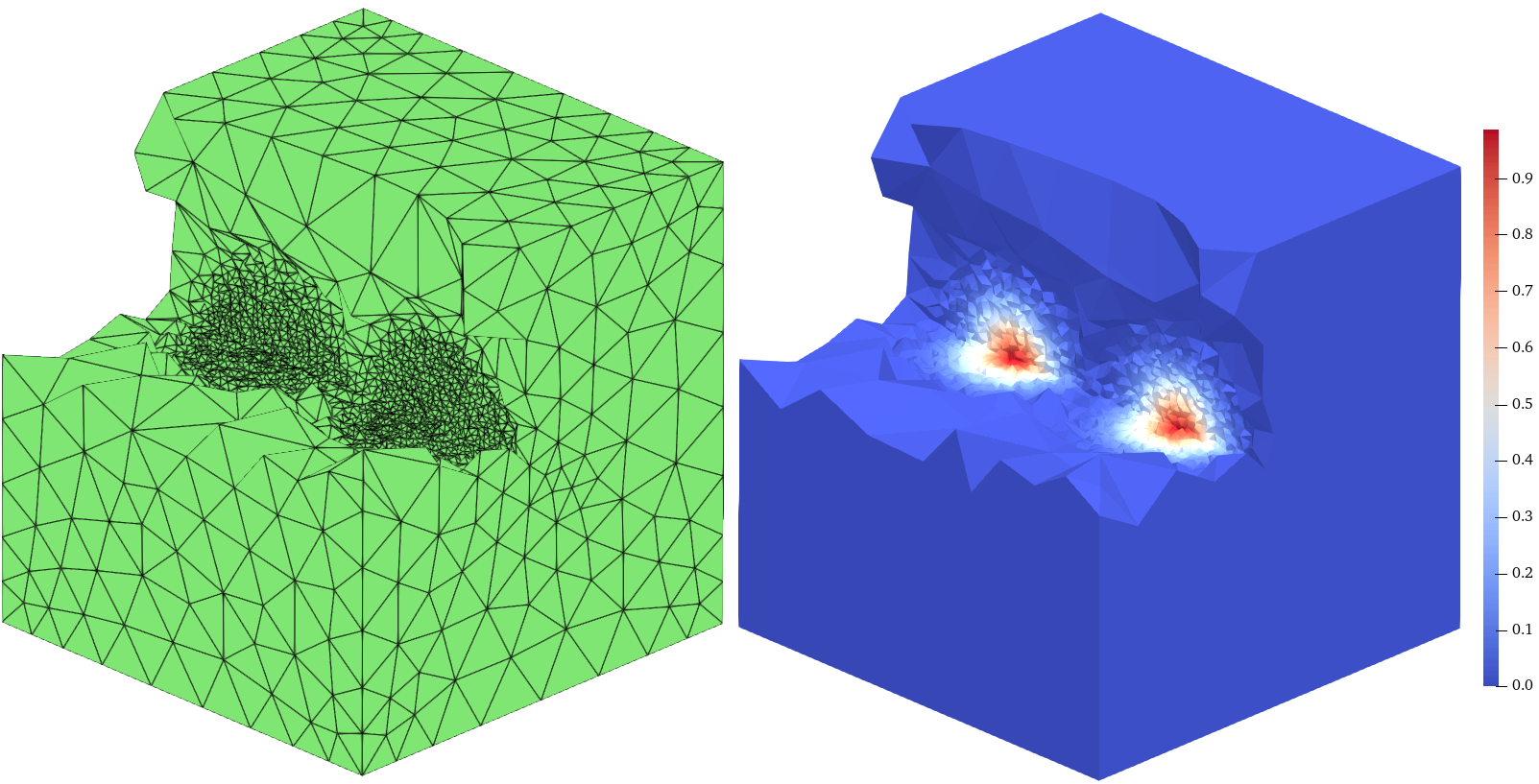}
			\includegraphics[width=0.47\linewidth]{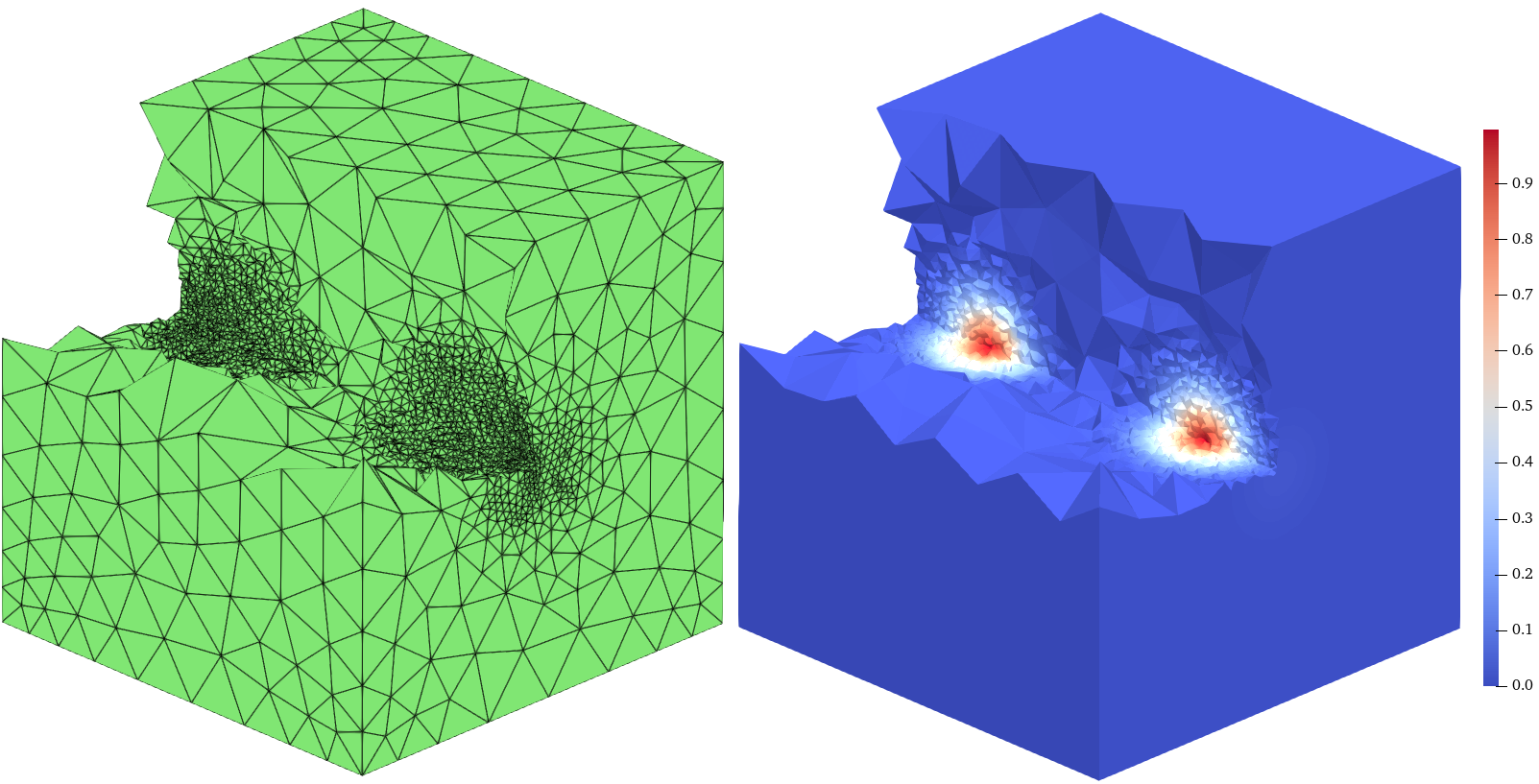}
			\caption{\Cref{example6}, snapshots of numerical solutions along with their corresponding adaptive meshes at $t=0.00, 0.25, 0.50$, and $t=0.75$.}
			\label{Ex3D3-Initial}
		\end{figure}
	\end{example}
	With the mesh refinement tolerance $eTol = 0.1$, \Cref{Ex3D3-Time} shows the adaptive process at $t=0$, terminating after seven steps and NOV at each iteration being $390$, $615$, $1062$, $2085$, $4606$, $10860$, and $27113$, respectively. 
	The mesh refinement process highlights that the adaptive algorithm achieves a highly efficient convergence procedure compared to conventional adaptive methods. 
	\Cref{Ex3D3-loss} presents the training loss at various time levels during the learning process. To reach the desired MSE, the initial time level requires more than 200 training epochs. In contrast, fewer than 20 epochs are sufficient at all subsequent time levels to achieve a comparable level of accuracy. 
	\Cref{Ex3D3-Initial} illustrates the evolution of the numerical solutions and corresponding meshes at four distinct time points, capturing the symmetric splitting of the initial peak into two separate regions moving in opposite directions over time.
	The quasi-optimal convergence rate of the adaptive algorithm for the gradient error is demonstrated in the right panel of \Cref{Ex3-order}.

	\subsection{Additional benchmark examples}
	As a further demonstration beyond the linear model equations on rectangular domains considered above, we apply the proposed method to benchmark tests involving nonlinear dynamics and non-rectangular geometry. 
	
	\begin{example}\label{example7}
		Consider the Allen–Cahn equation
		\[
		u_t - \Delta u= (u - u^3)/\varepsilon^2, \quad x \in \Omega=[-1,1]^2, \ t > 0,
		\]
		with a Dirichlet boundary condition
		$$
		u(x,y,t)|_{\partial \Omega} = -1.
		$$
		The initial condition corresponds to a circular interface of radius $0.25$ centered at $(0.5,0.5)$:
		\[
		u(x,y,0) = \tanh\left( \frac{0.25 - \sqrt{(x-0.5)^2 + (y-0.5)^2}}{\sqrt{2}\,\varepsilon} \right).
		\]
		The parameter $\varepsilon = 0.01$ determines the interface thickness, which scales as $O(\varepsilon)$; such a small value leads to a sharp transition layer and poses a challenge for numerical resolution. 
		The simulation proceeds to the final time $T = 0.03$ with a fixed time step $\tau = 0.0003$, and the adaptive algorithm employs an error tolerance $eTol = 0.5$ to guide local mesh refinement. 
		Time discretization and spatial treatment follow the same scheme as in the preceding examples, with the nonlinear term handled via Newton's method at each time step. 
		The neural network architecture and training configuration are identical to those employed in the 2D benchmark case. 
		The adaptive meshes and corresponding finite element approximations are plotted in \Cref{ExAC_mesh}. It shows that the mesh adaptation captures the sharp interface as it moves. 
		\begin{figure}
			\includegraphics[width=0.48\linewidth]{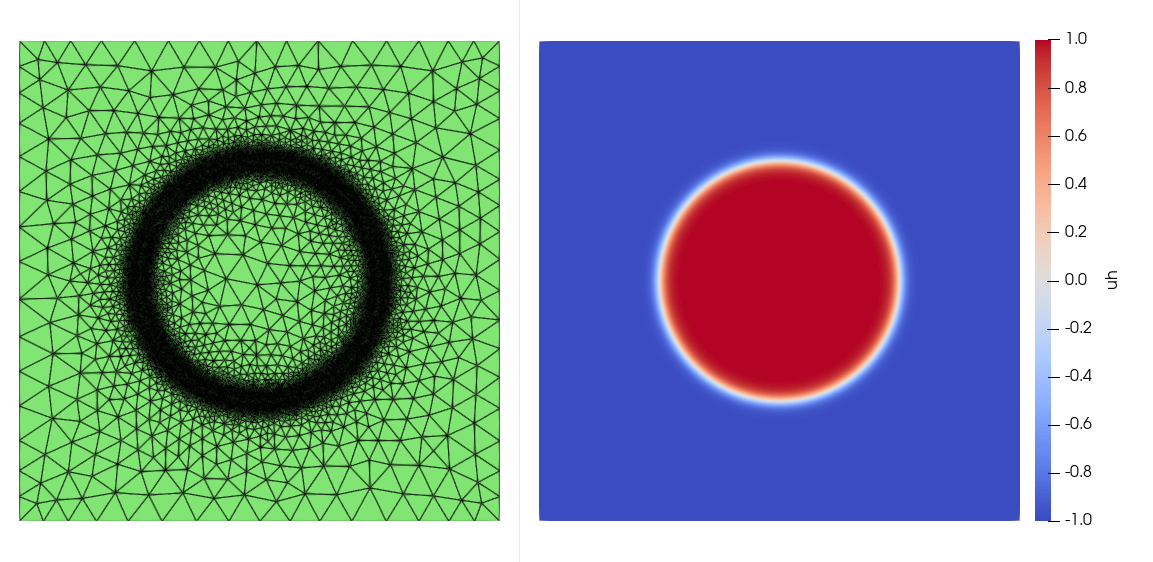}
			\includegraphics[width=0.48\linewidth]{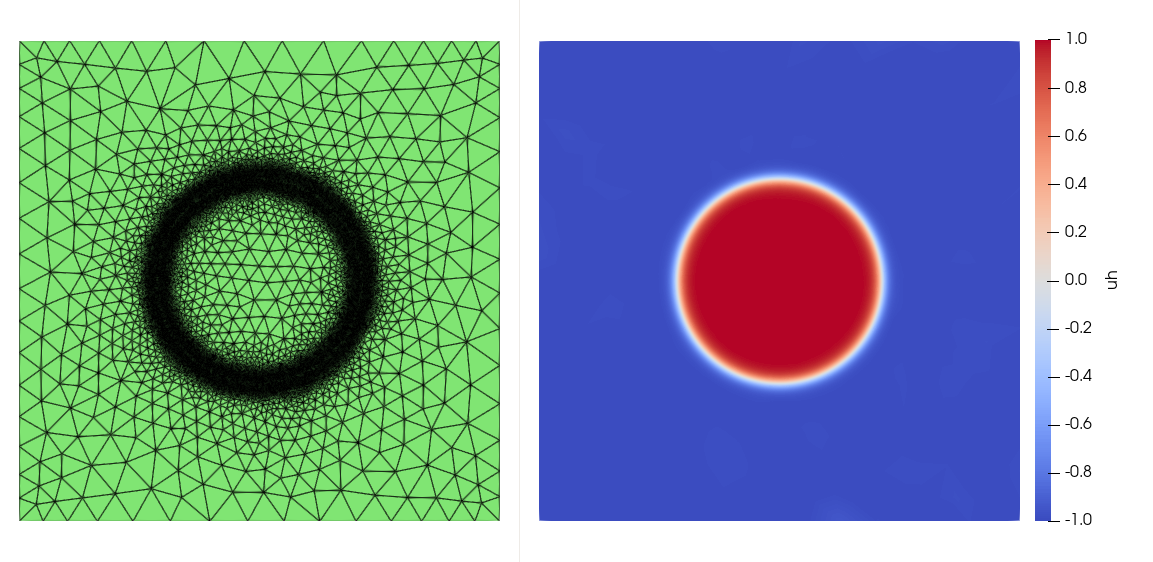}
			\includegraphics[width=0.48\linewidth]{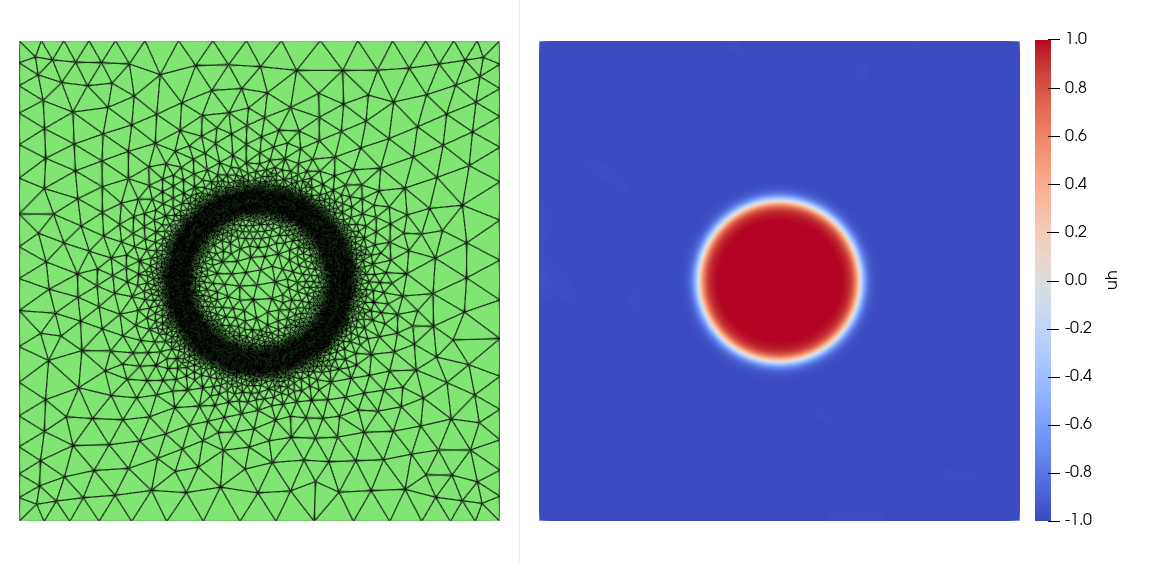}
			\includegraphics[width=0.48\linewidth]{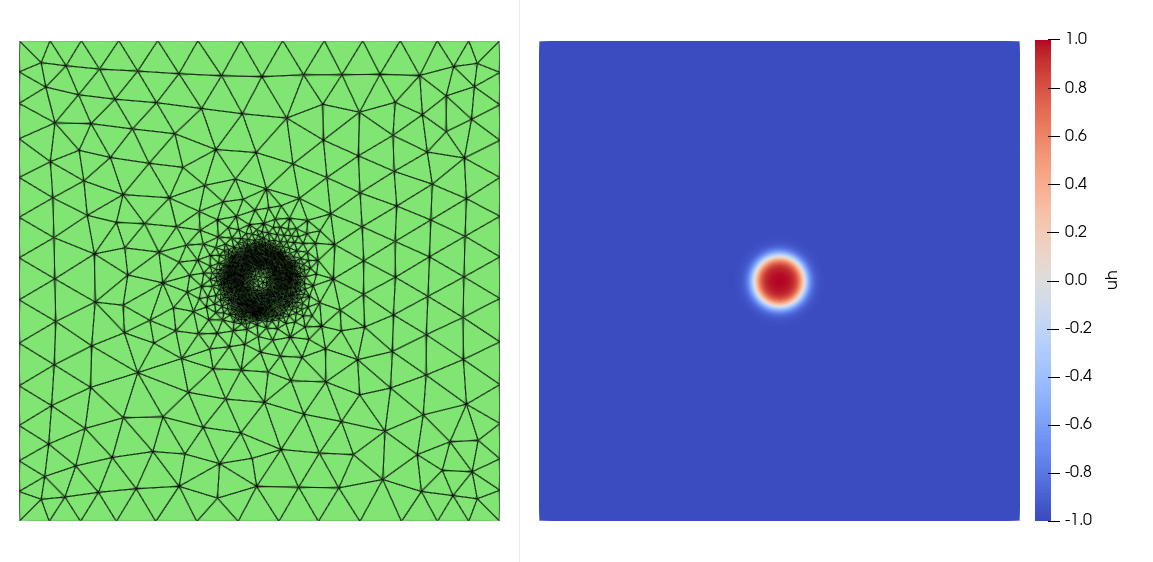}
			\caption{\Cref{example7}, Snapshots of numerical solutions along with their corresponding adaptive meshes at $t=0.00, 0.01, 0.02$, and $t=0.03$.}
			\label{ExAC_mesh}
		\end{figure}
	\end{example}
	
	\begin{example}\label{ex:lshaped}
		Consider the parabolic problem 
		\[ u_t-\Delta u=f(x,y,t), \quad (x,y)\in\Omega,\quad t>0, \] 
		where $f$ is a time-periodic localized heat source defined by
		\[
		f(x,y,t)= 
		\begin{cases} 
			60, & 0\leq \operatorname{mod}(t,0.2)\leq 0.04,\quad 
			-0.6\leq x\leq -0.4,\quad -0.6\leq y\leq -0.4,\\
			40, & 0.02\leq \operatorname{mod}(t,0.2)\leq 0.08,\quad 
			-0.6\leq x\leq -0.4,\quad 0.4\leq y\leq 0.6,\\
			80, & 0.10\leq \operatorname{mod}(t,0.2)\leq 0.18,\quad 
			0.4\leq x\leq 0.6,\quad 0.4\leq y\leq 0.6,\\
			0, & \text{otherwise}. 
		\end{cases} 
		\]
		on the L-shaped domain $\Omega=[-1,1]^2\setminus\big([0,1]\times[-1,0]\big)$, subject to the homogeneous Dirichlet boundary condition $u(x,y,t)|_{\partial\Omega}=0$ and the initial condition $u(x,y,0)=0$. 
		The corresponding distance function $d(\mathbf{x})$ and the continuous extension $\tilde{g}(\mathbf{x},t_{n-1})$ used in the neural network construction are provided in \Cref{app:func}. 
		The time step size, neural network architecture, and training configuration follow the settings specified at the beginning of this section. 
		\Cref{ExLshape_mesh} displays the adaptive meshes and the corresponding numerical solutions at $t = 0.04, 0.14, 0.21$ and $t = 0.24$.
		The temporal variation of the localized heat source induces significant changes in the solution profile throughout the simulation. The results demonstrate that the adaptive algorithm effectively concentrates mesh refinement in regions where the solution exhibits rapid variation, thereby accurately resolving the evolving solution features on this geometrically complex domain.
		\begin{figure}
			\includegraphics[width=0.45\linewidth]{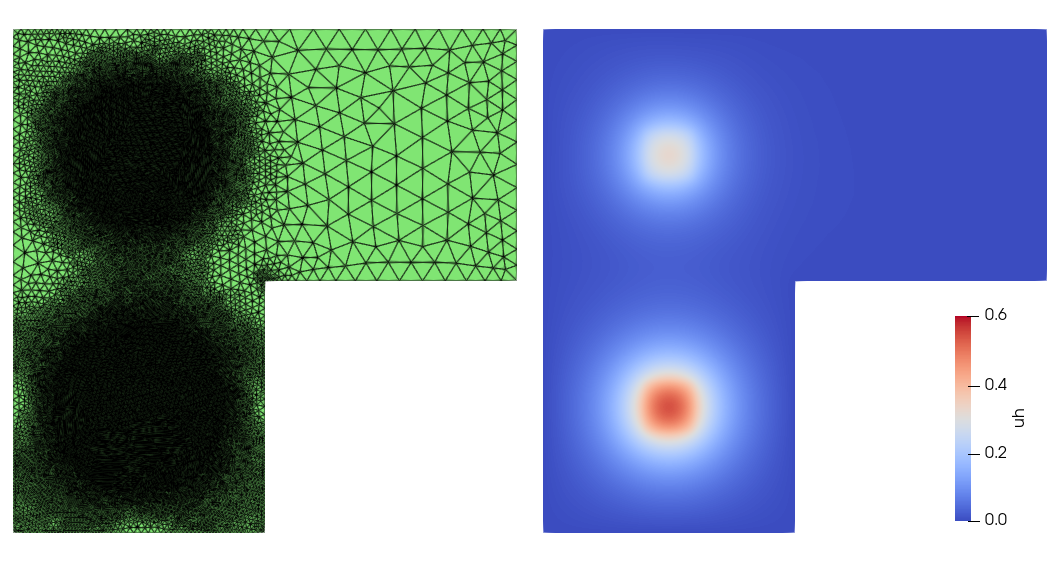}\hspace{0.5cm}
			\includegraphics[width=0.45\linewidth]{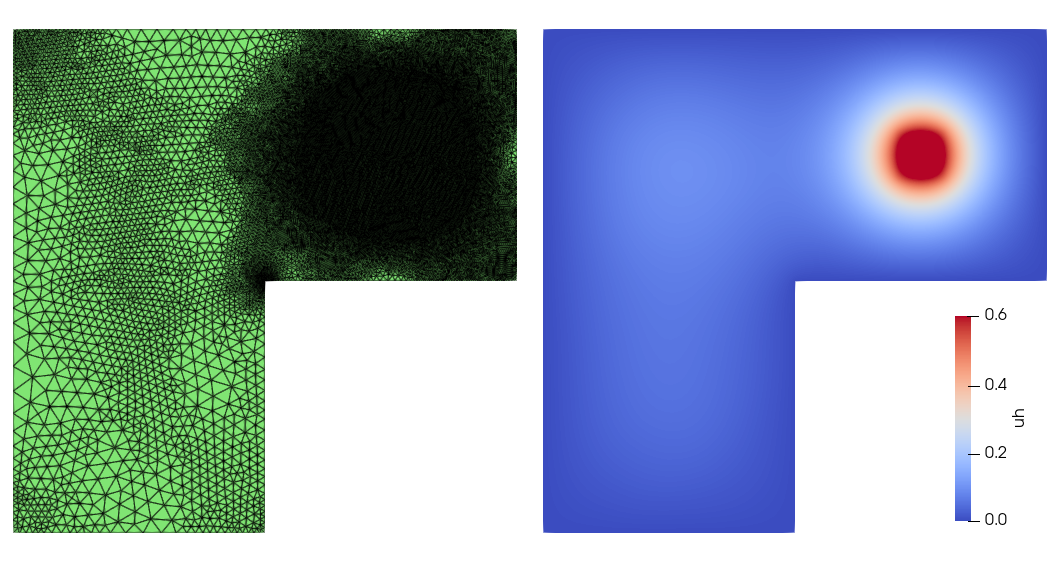}\\
			\includegraphics[width=0.45\linewidth]{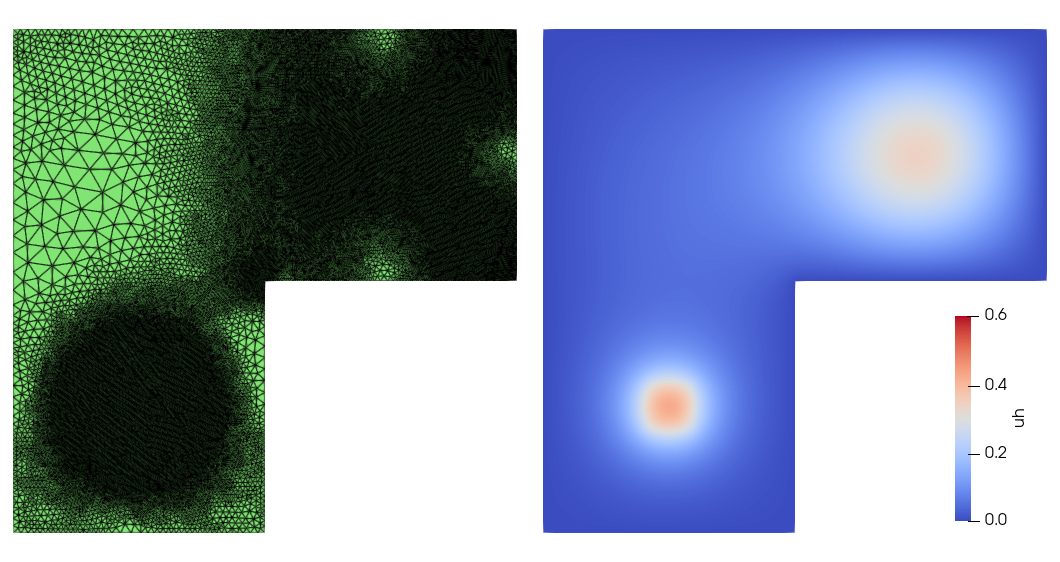}\hspace{0.5cm}
			\includegraphics[width=0.45\linewidth]{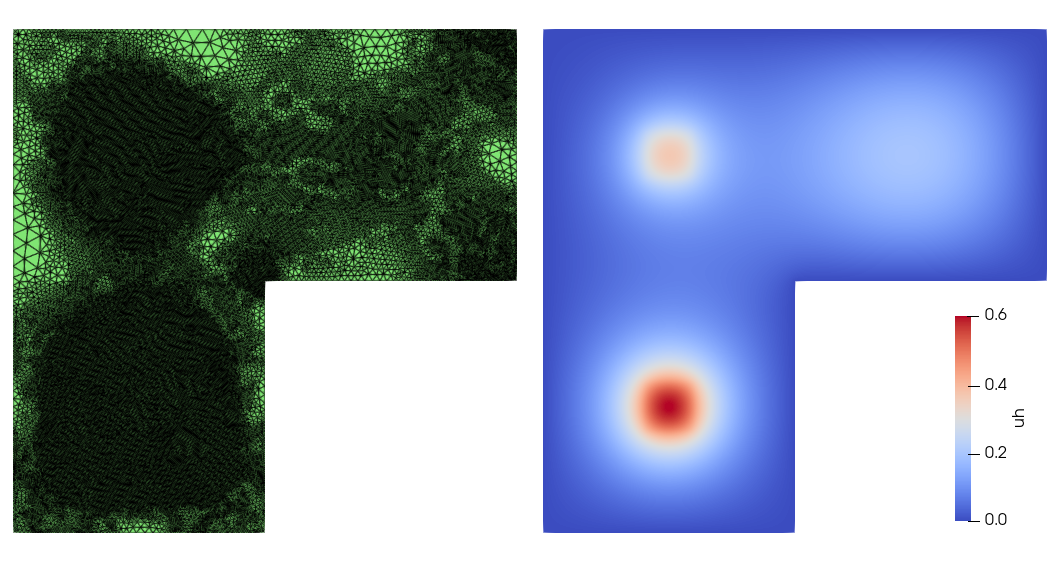}
			\caption{\Cref{ex:lshaped}, Snapshots of numerical solutions along with their corresponding adaptive meshes at $t=0.04$, $0.14$, $0.21$, and $t=0.24$.}
			\label{ExLshape_mesh}
		\end{figure}
	\end{example}

	\section*{Conclusion}
	In this paper, we proposed a neural network-enhanced $hr$-adaptive finite element method for parabolic equations that overcomes key limitations of conventional adaptive approaches, including excessive refinement iterations, costly interpolation between non-nested meshes, and the accumulation of unnecessary degrees of freedom across time steps. By employing a mesh-free neural network surrogate to replace the traditional interpolation operator, together with a controllable mesh size field strategy and a mesh-reset mechanism, the algorithm significantly reduces computational overhead while maintaining high accuracy. The combined error estimator ensures effective tracking of evolving singularities, and the adaptive process terminates within $7$ iterations per time step. Numerical experiments in both $2D$ and $3D$ demonstrate the efficiency, robustness, and quasi-optimal convergence behavior of the proposed method, highlighting the potential of integrating neural networks as specialized modules within classical finite element frameworks.
	
	In this work, we assume and numerically verify that the chosen network capacity is sufficient for the model problems under consideration, which are already challenging for conventional AFEM. This provides a reasonable starting point for the development of a complete and practical algorithmic framework.
	Nevertheless, several important questions remain open for future investigation, particularly the identification of minimal or near-optimal network architectures that can achieve the desired approximation accuracy with reduced computational complexity.
	
	A potential limitation of the current approach arises when the complexity of the solution increases substantially during the simulation, or when highly localized and multiscale features emerge. In such situations, a fixed network architecture may no longer provide sufficient approximation power. This challenge is closely related to the well-known spectral bias of neural networks, which can make the representation of increasingly complex high-frequency components more difficult.
	Although such behavior was not observed in the numerical experiments presented in this work, it may become relevant for more demanding applications. 
	One possible remedy is adaptive network enrichment. 
	For example, the network capacity could be increased whenever a suitable error indicator exceeds a prescribed tolerance over a sustained period. The parameters associated with the existing network could be inherited from the previously trained model, while newly introduced parameters could be initialized independently.
	Another possible strategy is to employ alternative neural network architectures. In our follow-up work on neural network-based interpolation, we compared several neural network architectures for representing finite element solutions, including feedforward neural networks and extreme learning machines \cite{hao2025datatransfer}. Our results indicate that certain architectures can achieve better approximation performance while maintaining low computational complexity. Consequently, selecting a more suitable network architecture may further improve the efficiency and robustness of the proposed framework.
	The design of reliable enrichment criteria, together with a systematic investigation of adaptive network-capacity strategies and alternative neural network architectures for problems with strongly localized or multiscale features, lies beyond the scope of the present work and will be explored in future research.

	\section*{Acknowledgments}
	The authors express sincere thanks to the referees for their useful comments and suggestions, which led to improvements of the presentation. 
	The authors thank Chunyu Chen for discussions about the implementation of the adaptive algorithm. 
	
	Huang and Yi were supported by the National Key R \& D Program of China (2024YFA1012600) and NSFC Project (12431014).
	Yin was supported by the University of Texas at El Paso Startup Award.
	
	\appendix
	\section{Functions in neural network \eqref{neuralnetwork}}\label{app:func}
	
	
	In this section, we provide the distance functions $d(\mathbf{x})$ and the 
	continuous extension $\tilde{g}(\mathbf{x},t_{n-1})$ of the boundary data 
	used in the neural network \eqref{neuralnetwork} for the examples in 
	\Cref{numExp}.

	\textbf{2D squared domain.}
	For the two-dimensional examples on $\Omega=[-1,1]^2$, we take
	\[
	d(\mathbf{x}) = (x^2-1)(y^2-1).
	\]  
	The continuous extension $\tilde{g}(\mathbf{x},t_{n-1})$ in 
	\eqref{neuralnetwork} is defined by
	\begin{align*}
		\tilde{g}(\mathbf{x},t_{n-1}) = &\ \frac{1-x}{2}u(-1,y,t_{n-1}) + \frac{1+x}{2}u(1,y,t_{n-1}) + \frac{1-y}{2}u(x,-1,t_{n-1}) + \frac{1+y}{2}u(x,1,t_{n-1}) \nonumber \\
		&- \frac{(1-x)(1-y)}{4}u(-1,-1,t_{n-1}) - \frac{(1-x)(1+y)}{4}u(-1,1,t_{n-1}) \nonumber \\
		&- \frac{(1+x)(1-y)}{4}u(1,-1,t_{n-1}) - \frac{(1+x)(1+y)}{4}u(1,1,t_{n-1}).
	\end{align*}
	
	\textbf{2D L-shaped domain.} 
	For the two-dimensional L-shaped domain $\Omega = [-1,1]^2 \setminus \big([0,1]\times[-1,0]\big)$, we choose 
	\begin{align*} 
		d_1(x,y) &= (1-x)(1+x)y(1-y),\\ d_2(x,y) &= (1+x)(-x)(1+y)(1-y), \end{align*} 
	and define 
	\begin{align*} 
		d(x,y) = d_1(x,y)+d_2(x,y) +\sqrt{d_1^2(x,y)+d_2^2(x,y)}. \end{align*} 
	Equivalently, 
	\begin{align*} 
		d(x,y) &= (1-x)(1+x)y(1-y) +(1+x)(-x)(1+y)(1-y) \\ &\quad + \sqrt{ \left[(1-x)(1+x)y(1-y)\right]^2 +\left[(1+x)(-x)(1+y)(1-y)\right]^2 }. 
	\end{align*}
	The continuous extension $\tilde{g}(x,y,t_{n-1})$ is given by 
	\begin{align*} 
		\tilde{g}(x,y,t_{n-1}) &= \frac{x(x-1)}{2}\,u(-1,y,t_{n-1}) + (1-x^2)\,u(0,y,t_{n-1}) + \frac{x(x+1)}{2}\,u(1,y,t_{n-1}) \\ 
		&\quad + \frac{y(y-1)}{2} \Bigg[ u(x,-1,t_{n-1}) - \frac{x(x-1)}{2}\,u(-1,-1,t_{n-1}) \\ 
		&\qquad\qquad - (1-x^2)\,u(0,-1,t_{n-1}) - \frac{x(x+1)}{2}\,u(1,-1,t_{n-1}) \Bigg] \\ 
		&\quad + (1-y^2) \Bigg[ u(x,0,t_{n-1}) - \frac{x(x-1)}{2}\,u(-1,0,t_{n-1}) \\ &\qquad\qquad - (1-x^2)\,u(0,0,t_{n-1}) - \frac{x(x+1)}{2}\,u(1,0,t_{n-1}) \Bigg]. 
	\end{align*}

	\textbf{3D domain.}
	For the three-dimensional examples on $\Omega=[-1,1]^3$, we choose 
	\begin{align*}
		d(\mathbf{x}) = (1-x^2)(1-y^2)(1-z^2).
	\end{align*}
	The continuous extension $\tilde{g}(\mathbf{x},t_{n-1})$ is given by
	\begin{align*}
		\tilde{g}(\mathbf{x},t_{n-1}) 
		&= \frac{1-x}{2}\,u(-1,y,z,t_{n-1}) + \frac{1+x}{2}\,u(1,y,z,t_{n-1}) \\
		&\quad + \frac{1-y}{2}\,u(x,-1,z,t_{n-1}) + \frac{1+y}{2}\,u(x,1,z,t_{n-1}) \\
		&\quad + \frac{1-z}{2}\,u(x,y,-1,t_{n-1}) + \frac{1+z}{2}\,u(x,y,1,t_{n-1}) \\
		&\quad - \frac{(1-x)(1-y)}{4}\,u(-1,-1,z,t_{n-1}) - \frac{(1-x)(1+y)}{4}\,u(-1,1,z,t_{n-1}) \\
		&\quad - \frac{(1+x)(1-y)}{4}\,u(1,-1,z,t_{n-1}) - \frac{(1+x)(1+y)}{4}\,u(1,1,z,t_{n-1}) \\
		&\quad - \frac{(1-x)(1-z)}{4}\,u(-1,y,-1,t_{n-1}) - \frac{(1-x)(1+z)}{4}\,u(-1,y,1,t_{n-1}) \\
		&\quad - \frac{(1+x)(1-z)}{4}\,u(1,y,-1,t_{n-1}) - \frac{(1+x)(1+z)}{4}\,u(1,y,1,t_{n-1}) \\
		&\quad - \frac{(1-y)(1-z)}{4}\,u(x,-1,-1,t_{n-1}) - \frac{(1-y)(1+z)}{4}\,u(x,-1,1,t_{n-1}) \\
		&\quad - \frac{(1+y)(1-z)}{4}\,u(x,1,-1,t_{n-1}) - \frac{(1+y)(1+z)}{4}\,u(x,1,1,t_{n-1}) \\
		&\quad + \frac{(1-x)(1-y)(1-z)}{8}\,u(-1,-1,-1,t_{n-1}) + \frac{(1-x)(1-y)(1+z)}{8}\,u(-1,-1,1,t_{n-1}) \\
		&\quad + \frac{(1-x)(1+y)(1-z)}{8}\,u(-1,1,-1,t_{n-1}) + \frac{(1-x)(1+y)(1+z)}{8}\,u(-1,1,1,t_{n-1}) \\
		&\quad + \frac{(1+x)(1-y)(1-z)}{8}\,u(1,-1,-1,t_{n-1}) + \frac{(1+x)(1-y)(1+z)}{8}\,u(1,-1,1,t_{n-1}) \\
		&\quad + \frac{(1+x)(1+y)(1-z)}{8}\,u(1,1,-1,t_{n-1}) + \frac{(1+x)(1+y)(1+z)}{8}\,u(1,1,1,t_{n-1}).
	\end{align*}

	\bibliographystyle{plain}   
	\bibliography{references}
	
\end{document}